# Regularization of Discrete Ill-Conditioned Problems Done Right - I

Ibrahima Dione

*Université de Moncton, Département de mathématiques et de statistique,*
*18 Antonine-Maillet Ave, Moncton, NB, E1A 3E9, Canada,*
*ibrahima.dione@umoncton.ca, ORCID: 0000-0003-4146-8753*

**Abstract**

When solving *rank-deficient* or *discrete ill-posed* problems by regularization methods, the choice of the regularization parameter is crucial. It is also of interest, the regularization norm used in the selection of the solution. In this work, we propose a stabilization of the existing regularization methods to address the delicate task of choosing this parameter. The analysis we carried out is independent of the chosen regularization norm. Under an unperturbed data least squares problem and of a maximal rank matrix, the stabilized-regularized method we propose provides the minimal norm solution whatever the chosen positive regularization parameter. And under a perturbed data least squares problem, this approach provides increasingly accurate and stable approximations of the minimal norm solution with respect to a refined mesh and a huge regularization parameter (over-regularization). We also investigate standard rank-deficient and ill-posed numerical examples corroborating the theoretical analysis, where the accuracy and the stability of the proposed approach is widely discussed.

## 1. Introduction

Many difficulties in the solving of linear algebraic systems, by direct or iterative numerical methods, are sometimes linked to ill-conditioning. A problem is said to be well-conditioned (or stable) if small variations in the data induce small variations in the solution. Otherwise, it is said to be ill-conditioned or unstable. The word "*conditioning*" is thus used to designate one or more constants that precisely measure the dependence of the solution on the data ([4] - Chap. 4).
The conditioning of the matrix of a linear system measures the influence of rounding errors on the solution of a given problem. A linear system facing high conditioning is called an *ill-conditioned problem*, meaning that its solution is potentially very sensitive to perturbations of its data.

Two classes of problems give rise to ill-conditioned problems: *rank-deficient problems* [24] and *discrete ill-posed problems* [6, 21]. These problems are omnipresent in applied sciences and may arise from a variety of settings including data analysis problems, machine learning algorithms [39, 38], image debluring problems [26], approximations of the Fredholm integral equation of the first kind (describing many mathematical models) [19], etc. Solving these problems have been the topic of numerous investigations, mainly through the following least squares formulation (see [24] and references within)

$$\min_{\boldsymbol{v}\in\mathbb{R}^n} \frac{1}{2} \|\mathbf{A}\boldsymbol{v} - \boldsymbol{b}\|_m^2 , \tag{1}$$

where $\mathbf{A} \in \mathcal{M}_{m,n}(\mathbb{R})$ is a severely ill-conditioned and possibly singular matrix with $m$ rows and $n$ columns, $\boldsymbol{b}$ is a vector of $\mathbb{R}^n$, and $\|\cdot\|_m$ (resp. $\|\cdot\|_n$) denotes the euclidean vector norm in $\mathbb{R}^m$ (resp. in $\mathbb{R}^n$). The least squares problem (1) is also equivalent to the following normal equation problem (see [11] - Theorem 7.4-4 and [1] - Lemma 7.2.1)

$$\text{Find } \boldsymbol{u} \in \mathbb{R}^n \text{ such that, } \mathbf{A}^t\mathbf{A}\boldsymbol{u} = \mathbf{A}^t\boldsymbol{b}, \tag{2}$$

where $\mathbf{A}^t$ is the transpose of the matrix $\mathbf{A}$. The normal equation (2) is sometimes used in the solving of ill-conditioned problems, despite that this is less recommanded for stability reasons ([37] - Chap. 1.3, [7] - Chap. 1.4).

One of the simplest and effectively used method in the solving of ill-conditioned problems is the Tikhonov regularization method ([24] - Chap. 5; [31] - Part II; [15] - Chap. 3, 4 and 5). This method is based on the idea that, in order to avoid numerical difficulties linked to the ill-conditioning, we need to incorporate about the desired solutions further informations of the type $\mathbf{L}\boldsymbol{u} = \mathbf{g}$. Here, the regularization matrix $\mathbf{L}$ of dimension $p \times n$ is not necessary the identity matrix $\mathbf{I}$, and $\mathbf{g}$ is a given $p$-vector. For incorporating these informations in the problem (1), the following minimization problem is thus considered

$$\text{Find } \boldsymbol{u}^{\star}_{\text{tik},\gamma} \in \mathbb{R}^n \text{ such that, } J_\gamma(\boldsymbol{u}^{\star}_{\text{tik},\gamma}) := \min_{\boldsymbol{v}\in\mathbb{R}^n} J_\gamma(\boldsymbol{v}), \tag{3}$$

where the quadratic regularization functional $J_\gamma(\cdot)$ is defined as follows

$$J_\gamma(\boldsymbol{v}) := \frac{1}{2}\left(\|\mathbf{A}\boldsymbol{v} - \boldsymbol{b}\|_m^2 + \gamma\,\|\mathbf{L}\boldsymbol{v} - \mathbf{g}\|_*^2\right),\ \forall\, \gamma > 0, \tag{4}$$

and where $\|\cdot\|_*^2$ denotes a chosen vector norm in $\mathbb{R}^p$, and $\gamma$ is the regularisation parameter which controls the weight given to the term $\|\mathbf{L}\boldsymbol{v} - \mathbf{g}\|_*^2$, relative to the minimisation of the residual $\|\mathbf{A}\boldsymbol{v} - \boldsymbol{b}\|_m^2$. When the euclidean norm $\|\cdot\|_p$ is used in place of $\|\cdot\|_*$, we obtain from (3) the following problem:

$$\text{Find } \boldsymbol{u}^{\star}_{\text{tik},\gamma} \in \mathbb{R}^n \text{ such that, } \left(\mathbf{A}^t\mathbf{A} + \gamma\mathbf{L}^t\mathbf{L}\right)\boldsymbol{u}^{\star}_{\text{tik},\gamma} = \mathbf{A}^t\boldsymbol{b} + \gamma\mathbf{L}^t\mathbf{g}. \tag{5}$$

**Remark 1.** *Under the choice* $\mathbf{L} = \mathbf{I}_{n\times n}$ *with* $\mathbf{g} = \mathbf{0}$, *it is etablished in [17] (Proposition 19.1) that the regularized solution* $\boldsymbol{u}^{\star}_{tik,\gamma}$ *converges to the minimum norm solution given by ([1] - Proposition 7.2.1)*

$$\boldsymbol{u}^{\star} = \mathbf{A}^{\dagger}\boldsymbol{b}, \tag{6}$$

*when the penalty parameter* $\gamma$ *goes to zero, where* $\mathbf{A}^{\dagger}$ *is the pseudoinverse matrix of* $\mathbf{A}$. *Whereas in [11] (Chapitre 8 -Theorem 8.6-3), it is established the same convergence but when* $\gamma$ *goes to* $+\infty$.

The principle of the regularization method is to replace the original problem (1) by the better conditioned and well-posed problem (3) in order to reduce the impacts of the data noises and rounding errors, and to produce a solution closed to the exact (noise-free) solution. Even if the regularization method provides accurate solutions, it also presents some drawbacks whose the good choice of the regularization parameter $\gamma$. The proper choice of this parameter is in the center of the solving process and has been the subject of many studies (see [23, 25, 28, 29, 32]). In fact, no regularization approach is completed without a method for the determination of the parameter providing an accurate solution. Several well-known methods, namely the *discrepancy principle* [30, 16], the *L-curve criterion* [25, 22], *Generalized cross-validation* [18, 40], the *Quasi-optimality criterion* [30, 16], are proposed for the determination of this parameter, and this still remains a very important research topic.
However, by rewriting the problem (3) under the following form

$$\text{Find } \boldsymbol{u}^{\star}_{\text{tik},\gamma} \in \mathbb{R}^n \text{ such that, } J_\gamma(\boldsymbol{u}^{\star}_{\text{tik},\gamma}) := \min_{\boldsymbol{v}\in\mathbb{R}^n} \frac{1}{2}\left\|\begin{bmatrix}\mathbf{A}\\ \sqrt{\gamma}\,\mathbf{L}\end{bmatrix}\boldsymbol{v} - \begin{bmatrix}\boldsymbol{b}\\ \sqrt{\gamma}\mathbf{g}\end{bmatrix}\right\|^2, \tag{7}$$

we clearly see that the regularization method, as it is applied for discrete ill-conditioned problems, does not act directly on the matrix $\mathbf{A}$ but completes the system by adding some additional informations. Thus, this method provides what is brought through these informations (more regular solution), but nothing reminds the system to also keep a small residual. This will come from the stabilization of the whole system in addition to the regularization, which is what we propose here.

In this paper, we propose a simple and practical approach that remedies to the choice of the regularization parameter. We introduce a stabilized-regularized method whose the basic idea is, when solving the *least squares problem* (1), we need to control at the same time the residual (through the *normal equation* (2)) and the solution or one of its derivative (through the type of informations $\mathbf{L}\boldsymbol{u} = \mathbf{g}$). In the solving of the stabilized-regularized equation proposed here, we use a direct method through the Matlab command \ (see [35] - Chap. 5.8, for more details on this command) and obtain the results:

• Under unperturbed data, the stabilized-regularized solution obtained $\boldsymbol{u}_\gamma^\star$, replicates the right hand term $\boldsymbol{b}$ whatever the chosen positive regularization parameter $\gamma > 0$ and with no proprieties on the matrix $\mathbf{A}$

$$\mathbf{A}\boldsymbol{u}^\star = \mathbf{A}\boldsymbol{u}_\gamma^\star, \; \forall\, \gamma > 0.$$

Moreover, with a maximal rank matrix $\mathbf{A}$, the solution $\boldsymbol{u}_\gamma^\star$ coincides with the minimal norm solution

$$\boldsymbol{u}^\star = \boldsymbol{u}_\gamma^\star, \; \forall\, \gamma > 0,$$

meaning that no regularization error is made. These theoretical results are coherent with the obtained numerical results.

• However, under a perturbed right hand term $\tilde{\boldsymbol{b}}$ this approach provides the following estimate between the stabilized-regularized solution $\widetilde{\boldsymbol{u}}_\gamma^\star$ and the minimal norm solution $\boldsymbol{u}^\star$

$$\left\|\mathbf{A}\boldsymbol{u}^\star - \mathbf{A}\widetilde{\boldsymbol{u}}_\gamma^\star\right\|_m \leq \left\|\left(\mathbf{I} + \gamma\mathbf{A}\mathbf{A}^t\right)\left(\boldsymbol{b} - \tilde{\boldsymbol{b}}\right)\right\|_m, \; \forall\, \gamma > 0.$$

And if the matrix $\mathbf{A}$ is of maximal rank, we hence obtain the estimate

$$\left\|\boldsymbol{u}^\star - \widetilde{\boldsymbol{u}}_\gamma^\star\right\|_n \leq \sqrt{\frac{1}{\lambda_n}}\left\|\left(\mathbf{I} + \gamma\mathbf{A}\mathbf{A}^t\right)\left(\boldsymbol{b} - \tilde{\boldsymbol{b}}\right)\right\|_m, \; \forall\, \gamma > 0,$$

where $\lambda_n > 0$ is the smallest eigenvalue of the matrix $\mathbf{A}^t\mathbf{A}$. The numerical results obtained here, provided increasingly accurate and stable approximations $\widetilde{\boldsymbol{u}}_\gamma^\star$ of the minimal norm solution as long as the mesh is refined and more regularization is added (over-regularization).

This work is organized as follows: we introduce in Section 2, the stabilized-regularized method where we give an insight on the accuracy of this approach through its Singular Value Decomposition (SVD). We also present in Section 3, the convergence analysis of this approach without and with perturbed data. Several numerical examples are presented in Section 4, to illustrate the theoretical results and specially to investigate the cases of perturbed data. Through these examples, we also investigate the accuracy and stability of this appoach.

## 2. Stabilization of the regularized discrete ill-conditioned problem

### *2.1.* Stabilization of the regularized problem

Based in the remark that regularization of discrete ill-conditioned problems are not stable (that is why the choice of the parameter $\gamma$ is needed), the approach we propose here is formulated from the problem (3) but with the additional worry to obtain stability. Hence, in addition to the type of prior information $\mathbf{L}\boldsymbol{u} = \mathbf{g}$ that we have already considered, we also introduce the normal equation in (2) as a constraint of the least squares problem (1). We thus interest in the following affine constrained minimization problem

$$\text{Find } \boldsymbol{u}, \text{ solution of the problem } \begin{cases} \min\limits_{\boldsymbol{v}\in\mathbb{R}^n} \dfrac{1}{2}\left\|\mathbf{A}\boldsymbol{v} - \boldsymbol{b}\right\|_m^2 \\ \text{Under the constraints:} \\ \left(\mathbf{A}^t\mathbf{A}\right)\boldsymbol{v} = \mathbf{A}^t\boldsymbol{b} \\ \qquad \mathbf{L}\boldsymbol{v} = \mathbf{g} \end{cases} \tag{8}$$

**Remark 2.** *What we mean from the formulation (8) is, when solving a discrete ill-conditioned problem under the least squares approach (1), even if we have in our disposal some a priori information of the form* $\mathbf{L}\boldsymbol{u} = \mathbf{g}$ *for regularization purpose, we also need to take into account the normal equation (2) as supplement informations about the solution for stability purpose. The constrained normal least squares problem obtained, namely problem (8), is thus the main difference between this study and all the works on the least squares problem.*

For a better handling of the problem (8), we introduce the quadratic functional $\mathcal{J}(\cdot)$ as follows

$$\begin{aligned}\mathcal{J}(\boldsymbol{v}) &= \frac{1}{2}\,\|\mathbf{A}\boldsymbol{v} - \boldsymbol{b}\|_m^2 - \frac{1}{2}\|\boldsymbol{b}\|_m^2 \\ &= \frac{1}{2}\,\langle \mathbf{A}\boldsymbol{v}, \mathbf{A}\boldsymbol{v}\rangle_m - \langle \boldsymbol{b}, \mathbf{A}\boldsymbol{v}\rangle_m \\ &= \frac{1}{2}\,\langle \mathbf{A}^t\mathbf{A}\boldsymbol{v}, \boldsymbol{v}\rangle_n - \langle \mathbf{A}^t\boldsymbol{b}, \boldsymbol{v}\rangle_n\,,\ \forall\, \boldsymbol{v} \in \mathbb{R}^n, \end{aligned} \tag{9}$$

where $\langle\,,\,\rangle_m$ and $\langle\,,\,\rangle_n$ refer to the scalar product associated to the euclidean norms $\|\cdot\|_m$ and $\|\cdot\|_n$ in $\mathbb{R}^m$ and $\mathbb{R}^n$, respectively. We also define the bilinear and linear forms $\mathcal{A}(\cdot,\cdot)$ and $\mathcal{L}(\cdot)$ by

$$\mathcal{A}(\boldsymbol{u}, \boldsymbol{v}) = \langle \mathbf{A}\boldsymbol{u},\ \mathbf{A}\boldsymbol{v}\rangle_m\,,\ \forall\, \boldsymbol{u}, \boldsymbol{v} \in \mathbb{R}^n, \tag{10}$$

$$\mathcal{L}(\boldsymbol{v}) = \langle \mathbf{A}^t\boldsymbol{b},\ \boldsymbol{v}\rangle_n\,,\ \forall\, \boldsymbol{v} \in \mathbb{R}^n. \tag{11}$$

The quadratic functional introduced in (9) can thus be rewritten under the following form

$$\mathcal{J}(\boldsymbol{v}) = \frac{1}{2}\mathcal{A}(\boldsymbol{v}, \boldsymbol{v}) - \mathcal{L}(\boldsymbol{v}),\ \forall\, \boldsymbol{v} \in \mathbb{R}^n, \tag{12}$$

and the normal least squares problem in (8) may be formulated under the minimization problem:

$$\text{Find } \boldsymbol{u}\ \in \boldsymbol{K} \text{ such that, } \mathcal{J}(\boldsymbol{u}) := \min_{\boldsymbol{v}\in\boldsymbol{K}} \mathcal{J}(\boldsymbol{v}), \tag{13}$$

where $\boldsymbol{K}$ is a closed convex subset of $\mathbb{R}^n$ defined as follows

$$\boldsymbol{K} = \left\{\boldsymbol{v} \in \mathbb{R}^n \mid \mathbf{A}^t\mathbf{A}\boldsymbol{v} = \mathbf{A}^t\boldsymbol{b},\ \text{and}\ \ \mathbf{L}\boldsymbol{v} = \mathbf{g}\right\}.$$

Under certain circumstances (see [17] - Theorem 13-4), the minimization problem (13) is equivalent to the following variational inequality on the convex subset $\boldsymbol{K} \subset \mathbb{R}^n$:

$$\text{Find } \boldsymbol{u} \in \boldsymbol{K} \text{ such that, } \mathcal{A}(\boldsymbol{u}, \boldsymbol{v} - \boldsymbol{u}) \geq \mathcal{L}(\boldsymbol{v} - \boldsymbol{u}),\ \boldsymbol{v}\ \in \boldsymbol{K}. \tag{14}$$

The solutions of the least squares problem we are looking for through the constrained problem (13), are characterized by the variational inequality (14) on the convex subset $\boldsymbol{K}$. However, the minimal norm solution defined in (6), despite being assumed belonging in $\boldsymbol{K}$, also verifies the following.

**Theorem 1.** *The minimum norm solution* $\boldsymbol{u}^\star$ *defined in (6), verifies the following equation*

$$\mathcal{A}(\boldsymbol{u}^\star, \boldsymbol{v}) = \mathcal{L}(\boldsymbol{v}),\ \forall\, \boldsymbol{v} \in \mathbb{R}^n. \tag{15}$$

*Conversely, if the matrix* $\mathbf{A}$ *is of maximal rank* $r = \min(m, n)$*, then the minimum norm solution* $\boldsymbol{u}^\star$ *is the unique vector verifying (15).*

*Proof.* *Since the matrix* $\mathbf{A}\mathbf{A}^\dagger$ *is symetric (i.e.* $\mathbf{A}\mathbf{A}^\dagger = \left(\mathbf{A}\mathbf{A}^\dagger\right)^t = (\mathbf{A}^\dagger)^t\mathbf{A}^t$*, where* $\mathbf{A}^\dagger$ *is the pseudoinverse matrix of* $\mathbf{A}$*) and taking into account that* $\boldsymbol{u}^\star = \mathbf{A}^\dagger\boldsymbol{b}$*, we obtain from definition (10) the following equility*

$$\begin{aligned}\mathcal{A}(\boldsymbol{u}^\star, \boldsymbol{v}) = \left\langle \mathbf{A}\mathbf{A}^\dagger\boldsymbol{b},\ \mathbf{A}\boldsymbol{v}\right\rangle_m &= \left\langle (\mathbf{A}^\dagger)^t\mathbf{A}^t\boldsymbol{b},\ \mathbf{A}\boldsymbol{v}\right\rangle_m \\ &= \left\langle \mathbf{A}^t(\mathbf{A}^\dagger)^t\mathbf{A}^t\boldsymbol{b},\ \boldsymbol{v}\right\rangle_n \\ &= \left\langle (\mathbf{A}\mathbf{A}^\dagger\mathbf{A})^t\boldsymbol{b},\ \boldsymbol{v}\right\rangle_n,\ \forall\, \boldsymbol{v} \in \mathbb{R}^n. \end{aligned} \tag{16}$$

*Using the MoorePenrose condition $\mathbf{A}\mathbf{A}^{\dagger}\mathbf{A} = \mathbf{A}$, we obtain from (16) the desired result (15).*

*Conversely, any vector $\boldsymbol{u}$ satisfying the equality (15), also verifies using definitions (10) and (11), the following equation*

$$\left\langle \mathbf{A}^t\mathbf{A}\boldsymbol{u} - \mathbf{A}^t\boldsymbol{b},\, \boldsymbol{v}\right\rangle_n = 0,\ \forall\, \boldsymbol{v} \in \mathbb{R}^n. \tag{17}$$

*It is thus clear from (17) that $\mathbf{A}^t\mathbf{A}\boldsymbol{u} - \mathbf{A}^t\boldsymbol{b} = 0$, which provides $\boldsymbol{u} = \left(\mathbf{A}^t\mathbf{A}\right)^{-1}\mathbf{A}^t\boldsymbol{b}$ since the $n \times n$ matrix $\mathbf{A}^t\mathbf{A}$ has an inverse (since the matrix $\mathbf{A}$ is of maximal rank). Otherwise, since the pseudoinverse of $\mathbf{A}$ is given by $\mathbf{A}^{\dagger} = \left(\mathbf{A}^t\mathbf{A}\right)^{-1}\mathbf{A}^t$ (see [1] - Section 2.7), then $\boldsymbol{u} = \boldsymbol{u}^{\star}$ which ends the proof.* □

Instead of working in the constrained set $\boldsymbol{K}$, we relaxe these constraints by applying the regularization method to the minimization problem (13). We thus interest in the following *stabilized-regularized least squares problem*:

$$\text{Find } \boldsymbol{u}_{\gamma}^{\star} \in \mathbb{R}^n \text{ such that, } \mathcal{J}_{\gamma}(\boldsymbol{u}_{\gamma}^{\star}) := \min_{\boldsymbol{v}\in\mathbb{R}^n} \mathcal{J}_{\gamma}(\boldsymbol{v}), \tag{18}$$

where the regularized functional $\mathcal{J}_{\gamma}(\cdot)$ is defined by

$$\mathcal{J}_{\gamma}(\boldsymbol{v}) := \mathcal{J}(\boldsymbol{v}) + \frac{\gamma}{2}\left\|\mathbf{A}^t\mathbf{A}\boldsymbol{v} - \mathbf{A}^t\boldsymbol{b}\right\|_n^2 + \frac{\gamma}{2}\|\mathbf{L}\boldsymbol{v} - \mathbf{g}\|_*^2,\ \forall\, \gamma > 0. \tag{19}$$

**Remark 3.** *The stabilized-regularized problem (18) is interpreted (see [21] - Chap. 4.4) as follows:*

- *The first right-hand term $\mathcal{J}(\boldsymbol{v})$ in (19) measures the better fitting to the data, such that the solution $\boldsymbol{u}_{\gamma}^{\star}$ predicts the right-hand side $\boldsymbol{b}$.*

- *The third right-hand term $\|\mathbf{L}\boldsymbol{v} - \mathbf{g}\|_*^2$ in (19) controls the solution (or it's first or second derivative), when it is dominated by high-frequency components with large amplitudes (which can be damaging when the data are noisy). This term provides regularity to the solution.*

- *As for the second right-hand term $\left\|\mathbf{A}^t\mathbf{A}\boldsymbol{v} - \mathbf{A}^t\boldsymbol{b}\right\|_n^2$ in (19), it controls the residual through the normal equation. By trying to give more regularity to the solution from the term $\|\mathbf{L}\boldsymbol{v} - \mathbf{g}\|_*^2$ and by taking the parameter $\gamma$ larger, we run the risk of less fitting the data. This is the reason why we propose to add this term.*

- *The balance between these three objectives, namely to measure good-fit, to provide more regularity to the solution and to control the residual, is more and more archieved as large as is the parameter $\gamma$.*

We may rewrite the stabilized-regularized problem (18) under the following form

$$\text{Find } \boldsymbol{u}_{\gamma}^{\star} \in \mathbb{R}^n \text{ such that, } \mathcal{J}_{\gamma}(\boldsymbol{u}_{\gamma}^{\star}) := \min_{\boldsymbol{v}\in\mathbb{R}^n} \frac{1}{2}\left\|\begin{bmatrix}\mathbf{A}\\ \gamma^{\frac{1}{2}}\mathbf{A}^t\mathbf{A}\\ \gamma^{\frac{1}{2}}\mathbf{L}\end{bmatrix}\boldsymbol{v} - \begin{bmatrix}\boldsymbol{b}\\ \gamma^{\frac{1}{2}}\mathbf{A}^t\boldsymbol{b}\\ \gamma^{\frac{1}{2}}\mathbf{g}\end{bmatrix}\right\|^2 - \frac{1}{2}\|\boldsymbol{b}\|_m^2, \tag{20}$$

where the first term in the right hand side of (20) is defined as follows

$$\left\|\begin{bmatrix}\mathbf{A}\\ \gamma^{\frac{1}{2}}\mathbf{A}^t\mathbf{A}\\ \gamma^{\frac{1}{2}}\mathbf{L}\end{bmatrix}\boldsymbol{v} - \begin{bmatrix}\boldsymbol{b}\\ \gamma^{\frac{1}{2}}\mathbf{A}^t\boldsymbol{b}\\ \gamma^{\frac{1}{2}}\mathbf{g}\end{bmatrix}\right\|^2 = \|\mathbf{A}\boldsymbol{v} - \boldsymbol{b}\|_m^2 + \gamma\left\|\mathbf{A}^t\mathbf{A}\boldsymbol{v} - \mathbf{A}^t\boldsymbol{b}\right\|_n^2 + \gamma\|\mathbf{L}\boldsymbol{v} - \mathbf{g}\|_*^2.$$

**Remark 4.** *We can clearly see from the formulation (20) that the stabilized-regularized problem (18) has an influence on both the matrix $\mathbf{A}$ through the term $\gamma^{\frac{1}{2}}\mathbf{A}^t\mathbf{A}$ (as a stabilizer) and the overall system by the term $\gamma^{\frac{1}{2}}\mathbf{L}$ (as a regulator): The desired goal is thus achieved!*
*In other words, the central idea of the method that we propose here could be summarized in these terms: A satisfactory approximation by regularization of an ill-conditioned problem requires both stabilization (of the system) and regularization (of the solution).*

In order to ensure a unique solution to the problem (20), we may assume that $\mathbf{rank}(\mathbf{L}) = p \leq n \leq m$ and that $\mathbf{rank}\begin{pmatrix} \mathbf{A}^t\mathbf{A} \\ \mathbf{L} \end{pmatrix} = n$, i.e., when the null spaces of $\mathbf{A}^t\mathbf{A}$ and $\mathbf{L}$ intersect trivially

$$\mathcal{N}(\mathbf{A}^t\mathbf{A}) \cap \mathcal{N}(\mathbf{L}) = \emptyset.$$

When the norm $\|\cdot\|_*$ in (19) is the euclidean vector norm $\|\cdot\|_p$ in $\mathbb{R}^p$, the minimization problem (18) is equivalent to the following variational problem:

$$\text{Find } \boldsymbol{u}_\gamma^\star \in \mathbb{R}^n \text{ such that, } \mathcal{A}_\gamma(\boldsymbol{u}_\gamma^\star, \boldsymbol{v}) = \mathcal{L}_\gamma(\boldsymbol{v}),\ \forall\, \boldsymbol{v} \in \mathbb{R}^n, \tag{21}$$

where the bilinear and the linear forms $\mathcal{A}_\gamma(\cdot,\cdot)$ and $\mathcal{L}_\gamma(\cdot)$ are given as follows

$$\begin{aligned} \mathcal{A}_\gamma(\boldsymbol{u}_\gamma^\star, \boldsymbol{v}) &:= \mathcal{A}(\boldsymbol{u}_\gamma^\star, \boldsymbol{v}) + \gamma \left\langle \mathbf{A}^t\mathbf{A}\boldsymbol{u}_\gamma^\star, \mathbf{A}^t\mathbf{A}\boldsymbol{v} \right\rangle_n + \gamma \left\langle \mathbf{L}\boldsymbol{u}_\gamma^\star, \mathbf{L}\boldsymbol{v} \right\rangle_p, \\ \mathcal{L}_\gamma(\boldsymbol{v}) &:= \mathcal{L}(\boldsymbol{v}) + \gamma \left\langle \mathbf{A}^t\mathbf{A}\mathbf{A}^t\boldsymbol{b},\, \boldsymbol{v} \right\rangle_n + \gamma \left\langle \mathbf{L}^t\mathbf{g},\, \boldsymbol{v} \right\rangle_n. \end{aligned}$$

By the Lax-Milgram lemma (see [17] - Theorem 13.5), the problem (21) admits a unique solution when the functional $\mathcal{A}_\gamma(\cdot,\cdot)$ is a continuous and elliptic bilinear form. For implementation purposes, we can write the problem (21) under the following vectorial form

$$\Big( \left(\mathbf{I} + \gamma\mathbf{A}^t\mathbf{A}\right) \mathbf{A}^t\mathbf{A} + \gamma\mathbf{L}^t\mathbf{L} \Big) \boldsymbol{u}_\gamma^\star = \left(\mathbf{I} + \gamma\mathbf{A}^t\mathbf{A}\right) \mathbf{A}^t\boldsymbol{b} + \gamma\mathbf{L}^t\mathbf{g}. \tag{22}$$

**Remark 5.** *In the solving of the stabilized-regularized equation (22) by a direct method, we will use the Matlab command \ (see [35] - Chap. 5.8, for more details on this command).*

### *2.2.* Insight on the stabilized-regularized solution from Singular Value Decomposition

For some insight into the stabilized-regularized solution $\boldsymbol{u}_\gamma^\star$ of the problem (22), we use the Singular Value Decomposition (SVD) [5] of the matrix $\mathbf{A} \in \mathcal{M}_{m,n}(\mathbb{R})$, with $m \geq n$,

$$\mathbf{A} = \mathbf{U}\boldsymbol{\Sigma}\mathbf{V}^t,\ \ \mathbf{A}^t\mathbf{A} = \mathbf{V}\boldsymbol{\Sigma}^2\mathbf{V}^t,\ \text{ and }\ \mathbf{A}^t\mathbf{A}\mathbf{A}^t\mathbf{A} = \mathbf{V}\boldsymbol{\Sigma}^4\mathbf{V}^t, \tag{23}$$

where $\boldsymbol{\Sigma} \in \mathcal{M}_{n,n}(\mathbb{R})$ is a diagonal matrix with the singular values $\sigma_i$, $i = 1, \cdots, n$, satisfying

$$\sigma_1 \geq \sigma_2 \geq \cdots \geq \sigma_n, \quad \text{and } \ \boldsymbol{\Sigma} = \mathrm{diag}(\sigma_1, \sigma_2, \cdots, \sigma_n).$$

The matrices $\mathbf{U} \in \mathcal{M}_{m,n}(\mathbb{R})$ and $\mathbf{V} \in \mathcal{M}_{n,n}(\mathbb{R})$ are orthogonal, that is $\mathbf{U}^t\mathbf{U} = \mathbf{I}$, and $\mathbf{V}\mathbf{V}^t = \mathbf{V}^t\mathbf{V} = \mathbf{I}$ (see [21] - Chap 3.2). For simplicity purpose, we set in the formulation (22) that $\mathbf{g} = \mathbf{0}$ and the regularization matrix $\mathbf{L} = \mathbf{I}$. We thus use the SVD decomposition (23) to write the solution $\boldsymbol{u}_\gamma^\star$ under the filtered expansion

$$\begin{aligned} \boldsymbol{u}_\gamma^\star &= \left(\mathbf{A}^t\mathbf{A} + \gamma\mathbf{A}^t\mathbf{A}\mathbf{A}^t\mathbf{A} + \gamma\mathbf{I}\right)^{-1} \left(\mathbf{I} + \gamma\mathbf{A}^t\mathbf{A}\right) \mathbf{A}^t\boldsymbol{b} & (24) \\ &= \left(\mathbf{V}\boldsymbol{\Sigma}^2\mathbf{V}^t + \gamma\mathbf{V}\boldsymbol{\Sigma}^4\mathbf{V}^t + \gamma\mathbf{V}\mathbf{V}^t\right)^{-1} \left(\mathbf{V}\mathbf{V}^t + \gamma\mathbf{V}\boldsymbol{\Sigma}^2\mathbf{V}^t\right) \mathbf{V}\boldsymbol{\Sigma}\mathbf{U}^t\boldsymbol{b} \\ &= \mathbf{V}\left(\boldsymbol{\Sigma}^2 + \gamma\boldsymbol{\Sigma}^4 + \gamma\mathbf{I}\right)^{-1} \mathbf{V}^t\mathbf{V}\left(\mathbf{I} + \gamma\boldsymbol{\Sigma}^2\right) \boldsymbol{\Sigma}\mathbf{U}^t\boldsymbol{b} \\ &= \mathbf{V}\left(\boldsymbol{\Sigma}^2 + \gamma\boldsymbol{\Sigma}^4 + \gamma\mathbf{I}\right)^{-1} \left(\mathbf{I} + \gamma\boldsymbol{\Sigma}^2\right) \boldsymbol{\Sigma}\mathbf{U}^t\boldsymbol{b} \\ &= \sum_{i=1}^{n} \left( \phi_i^{[\gamma]} \frac{\boldsymbol{u}_i^t\boldsymbol{b}}{\sigma_i} \right) \boldsymbol{v}_i & (25) \end{aligned}$$

where $\boldsymbol{u}_i$ and $\boldsymbol{v}_i$ , $i = 1, \cdots, n$, are columns of the matrices $\mathbf{U}$ and $\mathbf{V}$, and $\phi_i^{[\gamma]}$, $i = 1, \cdots, n$, are the filter factors defined as follows

$$\phi_i^{[\gamma]} = \frac{\sigma_i^2 + \gamma\sigma_i^4}{\sigma_i^2 + \gamma\sigma_i^4 + \gamma}. \tag{26}$$

Furthermore, the SVD expansion of the standard Tikhonov regularization solution obtained from (5) is given by (see [5] - Chap. 2.7, or [21] - Chap. 4.4)

$$\boldsymbol{u}^{\star}_{\text{tik},\gamma} = \sum_{i=1}^{n} \left( \varphi_i^{[\gamma]} \frac{\boldsymbol{u}_i^t \boldsymbol{b}}{\sigma_i} \right) \boldsymbol{v}_i, \tag{27}$$

where the filter factors $\varphi_i^{[\gamma]}$, $i = 1, \cdots, n$, are given by

$$\varphi_i^{[\gamma]} = \frac{\sigma_i^2}{\sigma_i^2 + \gamma}. \tag{28}$$

The behavior of filter factors in (26) and in (28) are illustrated in Figure 1, for different values of the regularization parameter ($\gamma = 10^{-6}, 10^{-3}, 1, 10^3, 10^6, 10^{12}$), and by considering (for simplicity) $\sigma_i$ as a continuous variable in $[10^{-4}, 1]$. For small values of this parameter (namely $\gamma = 10^{-6}, 10^{-3}, 1$), we observe no difference between these filter factors. We see that for singular values $\sigma_i$ larger than the regularization parameter $\gamma$, the filter factors are close to one, meaning that the corresponding SVD components contribute to solutions $\boldsymbol{u}^{\star}_{\gamma}$ and $\boldsymbol{u}^{\star}_{\text{tik},\gamma}$ with almost full strength (see [21] - Chap 4.4). And for singular values much smaller than the parameter $\gamma$, the filter factors are small therefore the SVD components are just filtered (or damped) to suppress the noise. The parameter $\gamma$ is always taken, under the Tikhonov method, as follows $\sigma_n^2 < \gamma < \sigma_1^2$.
However huge differences appear when the regularization parameter is greater than 1 (namely $\gamma > \sigma_1^2$), that is when one over-regularizes. We see in Figure 1 filter factors from Tikhonov method going toward zero, when $\gamma$ becomes more and more large by successively taking values $\gamma = 10^3, 10^6, 10^{12}$. Whereas for the stabilized-regularized method, the situation is completly different. We observe filter factors corresponding to the large singular values are only slightly damped (which is desirable since these singular values are the most representative), while filter factors associated to small singular values are effectively filtered (since the filter factors $\phi_i^{[\gamma]}$ decay rapidly like $\sigma^4$). Moreover, we abtain from equation (26) the limit asymptotic behavior of the filter factors when $\gamma \to \infty$

$$\phi_i^{[\gamma]} \longrightarrow \phi_i^{[\infty]} = \frac{\sigma_i^4}{\sigma_i^4 + 1}, \; i = 1, \cdots, n, \;\; \text{for} \;\; \gamma \to \infty. \tag{29}$$

**Remark 6.** *The SVD expansion we reach by over-regularizing the stabilized-regularized method is*

$$\boldsymbol{u}^{\star}_{\infty} = \sum_{i=1}^{n} \left( \phi_i^{[\infty]} \frac{\boldsymbol{u}_i^t \boldsymbol{b}}{\sigma_i} \right) \boldsymbol{v}_i. \tag{30}$$

*This solution should be the more accurate solution one may hope when the mesh is refined, that is whenever the dimensions $(m, n)$ become large enough (see Figure 9 and Figure 11).*

As for a brief look at the statistical aspects of the stabilized-regularized method, we consider the perturbed right hand side $\widetilde{\boldsymbol{b}} = \boldsymbol{b} + \boldsymbol{e}$, where errors $\boldsymbol{e} \in \mathbb{R}^m$ are Gaussian white noises, i.e. the elements of $\boldsymbol{e}$ are from Gaussian distribution with zero mean and standard deviation denoted by $\eta$. The covariance matrix of $\boldsymbol{e}$ is thus given by

$$\mathbb{C}ov(\boldsymbol{e}) = \eta^2 \mathbf{I}.$$

We compute the covariance matrix of the stabilized-regularized solution $\widetilde{\boldsymbol{u}}^{\star}_{\gamma}$ from (24) and (23)

$$\mathbb{C}ov(\widetilde{\boldsymbol{u}}^{\star}_{\gamma}) = \eta^2 \sum_{i=1}^{n} \left( \phi_i^{[\gamma]} \right)^2 \sigma_i^{-2} \boldsymbol{v}_i \boldsymbol{v}_i^t.$$

The norm of this covariance matrix is hence given by

$$\begin{aligned} \|\mathbb{C}ov(\widetilde{\boldsymbol{u}}^{\star}_{\gamma})\|_1 = \|\mathbb{C}ov(\widetilde{\boldsymbol{u}}^{\star}_{\gamma})\|_\infty &= \eta^2 \max_{i=1,\cdots,n} \left| \left( \phi_i^{[\gamma]} \right)^2 \sigma_i^{-2} \right| \\ &= \eta^2 \max_{i=1,\cdots,n} \left| \frac{\left(1 + \gamma \sigma_i^2\right)^2 \sigma_i^2}{\left(\sigma_i^2 + \gamma \sigma_i^4 + \gamma\right)^2} \right|. \end{aligned} \tag{31}$$

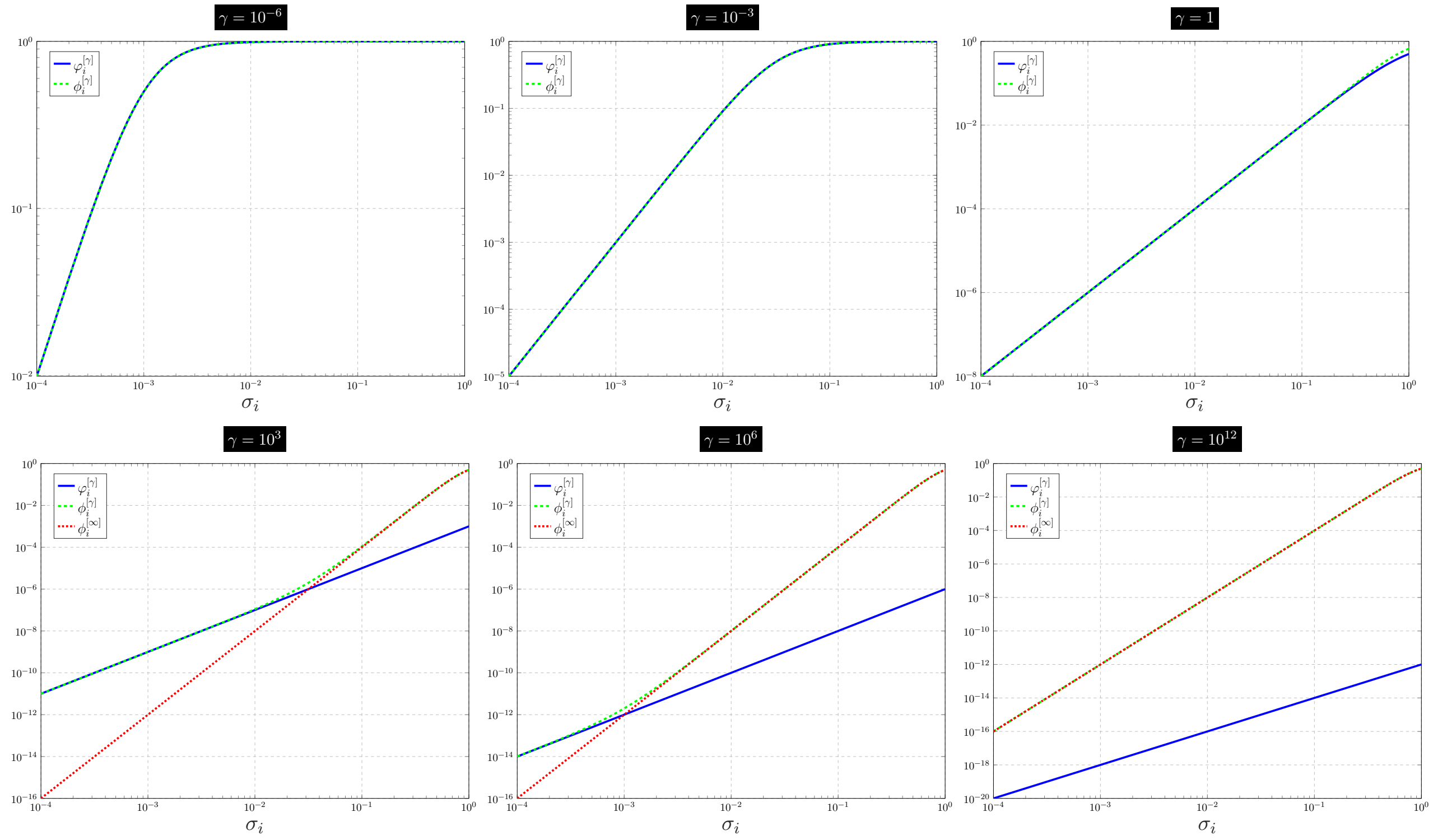

Figure 1: Filter factors from the Tikhonov regularization method $\varphi_i^{[\gamma]}$, $i = 1, \cdots, n$, (in blue line) and from the stabilized-regularized method $\phi_i^{[\gamma]}$, $i = 1, \cdots, n$, (in green dashed line), with respect to the singular values $\sigma_i$ taken as continuous variables in $[10^{-4}, 1]$ for different values of the parameter $\gamma = 10^{-6}, 10^{-3}, 1, 10^3, 10^6, 10^{12}$. The limit filter factors $\phi_i^{[\infty]}$, $i = 1, \cdots, n$, (in dotted red) are also represented in the last three cases ($\gamma = 10^3, 10^6$ and $10^{12}$) to illustrate how the filter factors $\phi_i^{[\gamma]}$ converges to $\phi_i^{[\infty]}$ (filter factors $\phi_i^{[\infty]}$ being supposed to provide the best approximated solution, when the mesh is refined). These are displayed on a log-log scale.

Hence elements of the covariance are small, specially when $\gamma$ is small and when $\gamma$ is large. This reduction of variance is achieved with a cost of bias on the stabilized-regularized solution $\widetilde{\boldsymbol{u}}_\gamma^\star$ computed from (25) (where $\boldsymbol{b}$ is repaced by $\widetilde{\boldsymbol{b}}$):

$$\begin{aligned}\mathbb{E}(\widetilde{\boldsymbol{u}}_\gamma^\star) &= \sum_{i=1}^{n} \phi_i^{[\gamma]} \langle \boldsymbol{u}^\star, \boldsymbol{v}_i \rangle \boldsymbol{v}_i \\ &= \boldsymbol{u}^\star - \sum_{i=1}^{n} \left(1 - \phi_i^{[\gamma]}\right) \langle \boldsymbol{u}^\star, \boldsymbol{v}_i \rangle \boldsymbol{v}_i. \end{aligned} \tag{32}$$

Since $\left(1 - \phi_i^{[\gamma]}\right) = \frac{\gamma}{\sigma_i^2 + \gamma\sigma_i^4 + \gamma}$ is smaller than $\left(1 - \varphi_i^{[\gamma]}\right) = \frac{\gamma}{\sigma_i^2 + \gamma}$, for any $i = 1, \cdots, n$ and $\gamma > 0$, then the bias from the the stabilized-regularized method proposed here is smaller than the one from the standard Tikhonov method. Moreover, if we assume that the discrete Picard condition is satisfied, we may expect that this bias is also small compared to $\|\boldsymbol{u}^\star\|_2$ (see [21] - Chap 4).

## 3. Convergence analysis of the stabilized-regularized problem

### *3.1.* Analysis of the stabilized-regularized problem under unperturbed data

We analyse here the accuracy and the stability of the stabilized-regularized formulation (18) under unperturbed data, by investigating the relationship between the minimal norm solution $\boldsymbol{u}^\star$ (resp. its associated data $\mathbf{A}\boldsymbol{u}^\star$) and the stabilized-regularized solution $\boldsymbol{u}_\gamma^\star$ (resp. its associated data $\mathbf{A}\boldsymbol{u}_\gamma^\star$). We assume that the problem (18) admits a unique solution for each $\gamma$.

**Theorem 2.** *The minimal norm solution $\boldsymbol{u}^\star$ defined in (6) and the stabilized-regularized solution $\boldsymbol{u}^\star_\gamma$ of the problem (18) verify the following equality*

$$\mathbf{A}\boldsymbol{u}^\star = \mathbf{A}\boldsymbol{u}^\star_\gamma,\ \forall\, \gamma > 0. \tag{33}$$

*Moreover, if the matrix $\mathbf{A}$ is of maximal rank, then these solutions coincide*

$$\boldsymbol{u}^\star = \boldsymbol{u}^\star_\gamma,\ \forall\, \gamma > 0. \tag{34}$$

*Proof.* *Let us introduce the functional $\mathscr{R}_\gamma(\boldsymbol{v})$ over $\mathbb{R}^n$ as follows (see [3] - Proof of Theorem 3.2)*

$$\mathscr{R}_\gamma(\boldsymbol{v}) := \frac{1}{2}\mathcal{A}(\boldsymbol{u}^\star - \boldsymbol{v}, \boldsymbol{u}^\star - \boldsymbol{v}) + \frac{\gamma}{2}\left\|\mathbf{A}^t\mathbf{A}\boldsymbol{v} - \mathbf{A}^t\boldsymbol{b}\right\|_n^2 + \frac{\gamma}{2}\left\|\mathbf{L}\boldsymbol{v} - \mathbf{g}\right\|_*^2,\ \forall\, \boldsymbol{v} \in \mathbb{R}^n.$$

*Developing and reordering terms from $\mathscr{R}_\gamma(\boldsymbol{v})$ by taking account (15), (19) and (12), we obtain*

$$\begin{aligned}
\mathscr{R}_\gamma(\boldsymbol{v}) &:= \frac{1}{2}\mathcal{A}(\boldsymbol{u}^\star, \boldsymbol{u}^\star) + \frac{1}{2}\mathcal{A}(\boldsymbol{v}, \boldsymbol{v}) - \mathcal{A}(\boldsymbol{u}^\star, \boldsymbol{v}) + \frac{\gamma}{2}\left\|\mathbf{A}^t\mathbf{A}\boldsymbol{v} - \mathbf{A}^t\boldsymbol{b}\right\|_n^2 + \frac{\gamma}{2}\|\mathbf{L}\boldsymbol{v} - \mathbf{g}\|_*^2 \\
&= \frac{1}{2}\mathcal{A}(\boldsymbol{v}, \boldsymbol{v}) - \mathcal{L}(\boldsymbol{v}) + \frac{\gamma}{2}\left\|\mathbf{A}^t\mathbf{A}\boldsymbol{v} - \mathbf{A}^t\boldsymbol{b}\right\|_n^2 + \frac{\gamma}{2}\|\mathbf{L}\boldsymbol{v} - \mathbf{g}\|_*^2 + \frac{1}{2}\mathcal{A}(\boldsymbol{u}^\star, \boldsymbol{u}^\star) \\
&= \mathcal{J}(\boldsymbol{v}) + \frac{\gamma}{2}\left\|\mathbf{A}^t\mathbf{A}\boldsymbol{v} - \mathbf{A}^t\boldsymbol{b}\right\|_n^2 + \frac{\gamma}{2}\|\mathbf{L}\boldsymbol{v} - \mathbf{g}\|_*^2 + \frac{1}{2}\mathcal{A}(\boldsymbol{u}^\star, \boldsymbol{u}^\star) \\
&= \mathcal{J}_\gamma(\boldsymbol{v}) + \frac{1}{2}\mathcal{A}(\boldsymbol{u}^\star, \boldsymbol{u}^\star),\ \forall\, \boldsymbol{v} \in \mathbb{R}^n
\end{aligned}$$

*Thus, functionals $\mathscr{R}_\gamma(\cdot)$ and $\mathcal{J}_\gamma(\cdot)$ are equal up to the additive constant $\frac{1}{2}\mathcal{A}(\boldsymbol{u}^\star, \boldsymbol{u}^\star)$. Hence minimizing $\mathscr{R}_\gamma(\cdot)$ or $\mathcal{J}_\gamma(\cdot)$ over $\mathbb{R}^n$ gives the same solution. Defining the sequence $\boldsymbol{v}_k = \boldsymbol{u}^\star + \frac{1}{\gamma^k}\boldsymbol{u}^\star$, for all $k \in \mathbb{Z}$, and taking into account that the minimal norm solution $\boldsymbol{u}^\star$ belongs to the convex subset $\boldsymbol{K}$, then we obtain the following*

$$\begin{aligned}
\mathscr{R}_\gamma(\boldsymbol{u}^\star_\gamma) &\leq \mathscr{R}_\gamma\left(\boldsymbol{u}^\star + \frac{1}{\gamma^k}\boldsymbol{u}^\star\right) \\
&= \frac{1}{2\gamma^{2k}}\mathcal{A}(\boldsymbol{u}^\star,\, \boldsymbol{u}^\star) + \frac{\gamma}{2\gamma^{2k}}\left(\left\|\mathbf{A}^t\boldsymbol{b}\right\|_n^2 + \|\mathbf{g}\|_*^2\right) \\
\mathcal{A}(\boldsymbol{u}^\star - \boldsymbol{u}^\star_\gamma, \boldsymbol{u}^\star - \boldsymbol{u}^\star_\gamma) &\leq \frac{1}{\gamma^{2k}}\left(\mathcal{A}(\boldsymbol{u}^\star,\, \boldsymbol{u}^\star) + \gamma\left(\left\|\mathbf{A}^t\boldsymbol{b}\right\|_2^2 + \|\mathbf{g}\|_*^2\right)\right) \\
\left\|\mathbf{A}\boldsymbol{u}^\star - \mathbf{A}\boldsymbol{u}^\star_\gamma\right\|_m^2 &\leq \frac{1}{\gamma^{2k}}\left(\|\mathbf{A}\boldsymbol{u}^\star\|_m^2 + \gamma\left(\left\|\mathbf{A}^t\boldsymbol{b}\right\|_n^2 + \|\mathbf{g}\|_*^2\right)\right),\ \forall\, k \in \mathbb{Z}.
\end{aligned} \tag{35}$$

*Letting thus $k$ goes to $+\infty$ in the case $\frac{1}{\gamma} \in ]0, 1[$, we thus obtain from (35) the result (33). Otherwise, if $\frac{1}{\gamma} > 1$, we let $k$ goes to $-\infty$ to also obtain the result (33).*

*If the matrix $\mathbf{A}$ is of maximal rank, then the symmetric matrix $\mathbf{A}^t\mathbf{A}$ is positive definite. Let $\lambda_1 \geq \ldots \geq \lambda_n \geq 0$ be its eigenvalues and let $\boldsymbol{a}_1, \boldsymbol{a}_2, \ldots, \boldsymbol{a}_n$ be the corresponding eigenvectors that form an orthonormal basis in $\mathbb{R}^n$ (see [34] - Chap. 5, [35] - Chap. 5.5). For any given vector $\boldsymbol{v} \in \mathbb{R}^n$, we thus have*

$$\boldsymbol{v} = \sum_{i=1}^{n} v_i\boldsymbol{a}_i, \quad \textit{with}\ \ \|\boldsymbol{v}\|_n^2 = \sum_{i=1}^{n} v_i^2.$$

*On the other hand, we obtain*

$$\begin{aligned}\|\mathbf{A}\boldsymbol{v}\|_m^2 &= \langle \mathbf{A}\boldsymbol{v}, \mathbf{A}\boldsymbol{v}\rangle_m \\ &= \left\langle \boldsymbol{v}, \mathbf{A}^t\mathbf{A}\boldsymbol{v}\right\rangle_n \\ &= \left\langle \sum_{i=1}^n v_i\boldsymbol{a}_i, \sum_{i=1}^n \lambda_i v_i\boldsymbol{a}_i \right\rangle_n \\ &= \sum_{i=1}^n \lambda_i v_i^2, \ \forall\, \boldsymbol{v} \in \mathbb{R}^n. \end{aligned} \tag{36}$$

*With respect to* (36)*, we finally obtain the following inequalities*

$$\lambda_n\|\boldsymbol{v}\|_n^2 \leq \|\mathbf{A}\boldsymbol{v}\|_m^2 \leq \lambda_1\|\boldsymbol{v}\|_n^2, \ \forall\, \boldsymbol{v} \in \mathbb{R}^n. \tag{37}$$

*By taking up* (35) *by using* (37)*, we finally arrive at the following estimate*

$$\begin{aligned}\lambda_n\|\boldsymbol{u}^\star - \boldsymbol{u}_\gamma^\star\|_n^2 &\leq \left\|\mathbf{A}\boldsymbol{u}^\star - \mathbf{A}\boldsymbol{u}_\gamma^\star\right\|_m^2 \\ &\leq \frac{1}{\gamma^k}\left(\lambda_1\left\|\boldsymbol{u}^\star\right\|_n^2 + \gamma\left(\left\|\mathbf{A}^t\boldsymbol{b}\right\|_n^2 + \|\mathbf{g}\|_*^2\right)\right) \\ \|\boldsymbol{u}^\star - \boldsymbol{u}_\gamma^\star\|_n^2 &\leq \frac{1}{\gamma^k}\left(\frac{\lambda_1}{\lambda_n}\left\|\boldsymbol{u}^\star\right\|_n^2 + \frac{\gamma}{\lambda_n}\left(\left\|\mathbf{A}^t\boldsymbol{b}\right\|_n^2 + \|\mathbf{g}\|_*^2\right)\right), \ \forall\, k \in \mathbb{Z}. \end{aligned} \tag{38}$$

*Letting thus $k$ goes to $+\infty$ in the case $\frac{1}{\gamma} \in ]0,1[$, we thus obtain from* (38) *the result* (34)*. Otherwise, if $\frac{1}{\gamma} > 1$, we let $k$ goes to $-\infty$ to also obtain the result* (34)*.* □

**Remark 7.**

- *Results established in Theorem 2 (namely equalities* (33) *and* (34)*) mean that the fair balance in incorporating information about the solution through* $\mathbf{L}\boldsymbol{v} = \mathbf{g}$ *and in keeping the normal equation* $\mathbf{A}^t\mathbf{A}\boldsymbol{v} = \mathbf{A}^t\boldsymbol{b}$ *satisfied, is reached at the minimal norm solution* $\boldsymbol{u}^\star$ *of the least squares problem* (1) *when the stabilized-regularized formulation* (18) *is used. This solution is geometrically interprated as the unique vector on* $(\mathbf{Ker}\,\mathbf{A})^\perp$ *whose its image under the matrix* $\mathbf{A}$ *is equal to the projection of* $\boldsymbol{b}$ *onto* $\mathbf{Im}\,\mathbf{A}$ *(see [1] - Chapter 7).*

- *Another important point achieved in the results* (33) *and* (34) *is, whatever the vector norm used in the regularization term of the functional* $\mathcal{J}_\gamma(\cdot)$ *in* (19)*, we obtain through the formulation* (18) *the same solution (the minimal norm solution). In a statistical viewpoint, it means that the stabilized Lasso regularization method (when* $\|\cdot\|_* = \|\cdot\|_{\ell_1}$ *in* (19)*, where* $\|\cdot\|_{\ell_1}$ *is the* $\ell_1-$*norm) and the stabilized Ridge regularization method (corresponding to the euclidean norm* $\|\cdot\|_* = \|\cdot\|_p$*) select the same solution.*

- *Another very important result from the Theorem 2 is that the formulation* (18) *provides a solution with a vanishing regularization error (i.e.* $\|\boldsymbol{u}^\star - \boldsymbol{u}_\gamma^\star\|_n = 0$*), when the matrix* $\mathbf{A}$ *is of maximum rank. Whereas, the error on the data is always zero (i.e.* $\|\mathbf{A}\boldsymbol{u}^\star - \mathbf{A}\boldsymbol{u}_\gamma^\star\|_m = 0$*), whatever the matrix* $\mathbf{A}$*.*

### *3.2.* Analysis of the stabilized-regularized problem under perturbed data

In the study of the sensitivity of the solution $\boldsymbol{u}^\star$ of the problem (1) to variations of the right hand data $\boldsymbol{b}$, when taking into account the results in Theorem 2, we may focus on analysing the sensitivity of the stabilized-regularized solution $\boldsymbol{u}_\gamma^\star$ from problem (21).
Sensitivity analysis is very important in the analysis of the stability of a given problem. In the solving of a given problem, even if we analytically know the exact right hand term $\boldsymbol{b}$, we numerically introduce errors due to rounding, discretization and inaccuracy of the measurement of this data. We thus work with a different right hand term $\tilde{\boldsymbol{b}}$ in the order of

$$\|\boldsymbol{b} - \tilde{\boldsymbol{b}}\|_m / \|\boldsymbol{b}\|_m \simeq \varepsilon,$$

where $\varepsilon$ is the relative precision of the machine used. If we denote by $\widetilde{\boldsymbol{u}}_{\gamma}^{\star}$ the solution of the minimization problem (21) associated to the right hand term $\tilde{\boldsymbol{b}}$, we thus seek to estimate the error $\|\boldsymbol{u}^{\star}-\widetilde{\boldsymbol{u}}_{\gamma}^{\star}\|_{n}$ with respect to the previous variation on the data.

**Theorem 3.** *Let $\boldsymbol{u}_{\gamma}^{\star}$ and $\widetilde{\boldsymbol{u}}_{\gamma}^{\star}$ be the solutions of problem (21) corresponding to the right hand terms $\mathbf{b}$ and $\tilde{\mathbf{b}}$, respectively. Then, we obtain the estimate*

$$\left\|\mathbf{A}\boldsymbol{u}_{\gamma}^{\star}-\mathbf{A}\widetilde{\boldsymbol{u}}_{\gamma}^{\star}\right\|_{m} \leq \left\|\left(\mathbf{I}+\gamma\mathbf{A}\mathbf{A}^{t}\right)\left(\boldsymbol{b}-\tilde{\boldsymbol{b}}\right)\right\|_{m}, \ \forall\,\gamma>0. \tag{39}$$

*Moreover if the matrix $\mathbf{A}$ is of maximal rank, we thus obtain the forthcoming estimate*

$$\left\|\boldsymbol{u}_{\gamma}^{\star}-\widetilde{\boldsymbol{u}}_{\gamma}^{\star}\right\|_{n} \leq \sqrt{\frac{1}{\lambda_{n}}}\left\|\left(\mathbf{I}+\gamma\mathbf{A}\mathbf{A}^{t}\right)\left(\boldsymbol{b}-\tilde{\boldsymbol{b}}\right)\right\|_{m}, \ \forall\,\gamma>0, \tag{40}$$

*where $\lambda_{n}>0$ is the smallest eigenvalue of the matrix $\mathbf{A}^{t}\mathbf{A}$.*

*Proof. Considering the problem (21) with respect to the right hand term $\boldsymbol{b}$ by taking test functions as $\boldsymbol{v}=\boldsymbol{u}_{\gamma}^{\star}-\widetilde{\boldsymbol{u}}_{\gamma}^{\star}$, we obtain the following*

$$\begin{aligned}\left\langle\mathbf{A}\boldsymbol{u}_{\gamma}^{\star},\,\mathbf{A}\left(\boldsymbol{u}_{\gamma}^{\star}-\widetilde{\boldsymbol{u}}_{\gamma}^{\star}\right)\right\rangle_{m}&+\gamma\left\langle\mathbf{A}^{t}\mathbf{A}\boldsymbol{u}_{\gamma}^{\star},\mathbf{A}^{t}\mathbf{A}\left(\boldsymbol{u}_{\gamma}^{\star}-\widetilde{\boldsymbol{u}}_{\gamma}^{\star}\right)\right\rangle_{n}+\gamma\left\langle\mathbf{L}\boldsymbol{u}_{\gamma}^{\star},\mathbf{L}\left(\boldsymbol{u}_{\gamma}^{\star}-\widetilde{\boldsymbol{u}}_{\gamma}^{\star}\right)\right\rangle_{p}\\&=\left\langle\mathbf{A}^{t}\boldsymbol{b},\,\left(\boldsymbol{u}_{\gamma}^{\star}-\widetilde{\boldsymbol{u}}_{\gamma}^{\star}\right)\right\rangle_{n}+\gamma\left\langle\mathbf{A}^{t}\mathbf{A}\mathbf{A}^{t}\boldsymbol{b},\,\left(\boldsymbol{u}_{\gamma}^{\star}-\widetilde{\boldsymbol{u}}_{\gamma}^{\star}\right)\right\rangle_{n}+\gamma\left\langle\mathbf{L}^{t}\mathbf{g},\,\left(\boldsymbol{u}_{\gamma}^{\star}-\widetilde{\boldsymbol{u}}_{\gamma}^{\star}\right)\right\rangle_{n}\end{aligned} \tag{41}$$

*From the problem (21) with respect to the right hand term $\tilde{\boldsymbol{b}}$, when making the same choice $\boldsymbol{v}=\boldsymbol{u}_{\gamma}^{\star}-\widetilde{\boldsymbol{u}}_{\gamma}^{\star}$ we also obtain the forthcoming equality*

$$\begin{aligned}\left\langle\mathbf{A}\widetilde{\boldsymbol{u}}_{\gamma}^{\star},\,\mathbf{A}\left(\boldsymbol{u}_{\gamma}^{\star}-\widetilde{\boldsymbol{u}}_{\gamma}^{\star}\right)\right\rangle_{m}&+\gamma\left\langle\mathbf{A}^{t}\mathbf{A}\widetilde{\boldsymbol{u}}_{\gamma}^{\star},\mathbf{A}^{t}\mathbf{A}\left(\boldsymbol{u}_{\gamma}^{\star}-\widetilde{\boldsymbol{u}}_{\gamma}^{\star}\right)\right\rangle_{n}+\gamma\left\langle\mathbf{L}\widetilde{\boldsymbol{u}}_{\gamma}^{\star},\mathbf{L}\left(\boldsymbol{u}_{\gamma}^{\star}-\widetilde{\boldsymbol{u}}_{\gamma}^{\star}\right)\right\rangle_{p}\\&=\left\langle\mathbf{A}^{t}\tilde{\boldsymbol{b}},\,\left(\boldsymbol{u}_{\gamma}^{\star}-\widetilde{\boldsymbol{u}}_{\gamma}^{\star}\right)\right\rangle_{n}+\gamma\left\langle\mathbf{A}^{t}\mathbf{A}\mathbf{A}^{t}\tilde{\boldsymbol{b}},\,\left(\boldsymbol{u}_{\gamma}^{\star}-\widetilde{\boldsymbol{u}}_{\gamma}^{\star}\right)\right\rangle_{n}+\gamma\left\langle\mathbf{L}^{t}\mathbf{g},\,\left(\boldsymbol{u}_{\gamma}^{\star}-\widetilde{\boldsymbol{u}}_{\gamma}^{\star}\right)\right\rangle_{n}\end{aligned} \tag{42}$$

*Differentiating equations (41) and (42) by using Cauchy-Schwarz inequality, we obtain*

$$\begin{aligned}\left\|\mathbf{A}\left(\boldsymbol{u}_{\gamma}^{\star}-\widetilde{\boldsymbol{u}}_{\gamma}^{\star}\right)\right\|_{m}^{2}&+\gamma\left\|\mathbf{A}^{t}\mathbf{A}\left(\boldsymbol{u}_{\gamma}^{\star}-\widetilde{\boldsymbol{u}}_{\gamma}^{\star}\right)\right\|_{n}^{2}+\gamma\left\|\mathbf{L}\left(\boldsymbol{u}_{\gamma}^{\star}-\widetilde{\boldsymbol{u}}_{\gamma}^{\star}\right)\right\|_{p}^{2}\\&=\left\langle\left(\boldsymbol{b}-\tilde{\boldsymbol{b}}\right),\,\mathbf{A}\left(\boldsymbol{u}_{\gamma}^{\star}-\widetilde{\boldsymbol{u}}_{\gamma}^{\star}\right)\right\rangle_{m}+\gamma\left\langle\mathbf{A}\mathbf{A}^{t}\left(\boldsymbol{b}-\tilde{\boldsymbol{b}}\right),\,\mathbf{A}\left(\boldsymbol{u}_{\gamma}^{\star}-\widetilde{\boldsymbol{u}}_{\gamma}^{\star}\right)\right\rangle_{m}\\&=\left\langle\left(\mathbf{I}+\gamma\mathbf{A}\mathbf{A}^{t}\right)\left(\boldsymbol{b}-\tilde{\boldsymbol{b}}\right),\,\mathbf{A}\left(\boldsymbol{u}_{\gamma}^{\star}-\widetilde{\boldsymbol{u}}_{\gamma}^{\star}\right)\right\rangle_{m}\\&\leq\left\|\left(\mathbf{I}+\gamma\mathbf{A}\mathbf{A}^{t}\right)\left(\boldsymbol{b}-\tilde{\boldsymbol{b}}\right)\right\|_{m}\left\|\mathbf{A}\left(\boldsymbol{u}_{\gamma}^{\star}-\widetilde{\boldsymbol{u}}_{\gamma}^{\star}\right)\right\|_{m}\end{aligned} \tag{43}$$

*From the estimate (43), we obtain the following estimate which proves (39)*

$$\left\|\mathbf{A}\left(\boldsymbol{u}_{\gamma}^{\star}-\widetilde{\boldsymbol{u}}_{\gamma}^{\star}\right)\right\|_{m} \leq \left\|\left(\mathbf{I}+\gamma\mathbf{A}\mathbf{A}^{t}\right)\left(\boldsymbol{b}-\tilde{\boldsymbol{b}}\right)\right\|_{m}. \tag{44}$$

*If the matrix $\mathbf{A}$ is of maximal rank, then taking into account (37) we establish from the estimate (44) the result (40).* □

**Corollary 1.** *Under assumptions of Theorem 2 and Theorem 3, we obtain the estimate*

$$\|\mathbf{A}\boldsymbol{u}^{\star}-\mathbf{A}\widetilde{\boldsymbol{u}}_{\gamma}^{\star}\|_{m} \leq \left\|\left(\mathbf{I}+\gamma\mathbf{A}\mathbf{A}^{t}\right)\left(\boldsymbol{b}-\tilde{\boldsymbol{b}}\right)\right\|_{m}, \ \forall\,\gamma>0. \tag{45}$$

*Moreover if the matrix $\mathbf{A}$ is of maximal rank, we hence obtain the forthcoming estimate*

$$\|\boldsymbol{u}^{\star}-\widetilde{\boldsymbol{u}}_{\gamma}^{\star}\|_{n} \leq \sqrt{\frac{1}{\lambda_{n}}}\left\|\left(\mathbf{I}+\gamma\mathbf{A}\mathbf{A}^{t}\right)\left(\boldsymbol{b}-\tilde{\boldsymbol{b}}\right)\right\|_{m}, \ \forall\,\gamma>0, \tag{46}$$

*where $\lambda_{n}>0$ is the smallest eigenvalue of the matrix $\mathbf{A}^{t}\mathbf{A}$.*

*Proof. We establish (45) and (46) by simultaneously considering results from theorems 2 and 3.* □

## 4. Numerical examples

In this section, we present three numerical well known examples in the litterature: the **shaw** 4.1, the **heat** 4.2 and the **phillips** 4.3 test problems [24]. These classical examples of *improperly posed* problems are such that any small errors (including computational errors) in the given data make the solution impossible or grossly inaccurate. The computations carried out here are through Matlab and illustrations of the proposed regularization technique are made using codes from Hansen's Matlab packages *Regularization Tools* [22]. We consider as regularization norm $\|\cdot\|_*$ on the functional $\mathcal{J}_\gamma(\cdot)$ in (19), the euclidean vector norm $\|\cdot\|_p$ of the space $\mathbb{R}^p$. These numerical examples are obtained from the Fredholm integral equation of the first kind, through two different discretization techniques: the quadrature method (for the **shaw** and **heat** test problems) and the Galerkin method with orthonormal basis functions (for the **phillips** test problem). However, these discrete problems often occur with a high conditioning and rank-deficient matrix, even under moderate dimensions. We consider for simplicity square matrices by taking $m = n$.

Our goal through this section is first to illustrate results established under unperturbed data (33) and (34) in Theorem 2. We use the **shaw** and **heat** numerical examples, discretized by a quadrature method and that leads to a linear system of the form $\mathbf{A}\boldsymbol{u} = \boldsymbol{b}$.
As for the validating of the results (45) and (46) in Corollary 1 obtained under perturbed data, we use all the mentioned three test problems. In fact, we first analyse the **shaw** discrete ill-posed problem used in the image reconstruction [36] (see also [23]). We also investigate the inverse **heat** problem [13], whose the ill-conditioning is controlled by a parameter. We finally explore the **phillips** test problem [33] that we discretize by the Galerkin method.
For each case, we not only analyse the numerical stability with respect to the regularization parameter when facing perturbations (or not) on the data, but we also investigate how our stabilized-regularized formulation is efficient, accurate and stable specially on large-scale problems.

Each of the examples we investigate here is constructed by first computing the matrix $\mathbf{A}$ and the vector $\boldsymbol{b}$ characterizing the least squares problem (1), along with the regularization matrix $\mathbf{L}$ whose the choice is based on the regularity of the solution that we would like to take into account. These data are extrated from Matlab codes *shaw(n)*, *heat(n, κ)*, *phillips(n)* and *get_l(n, d)* detailed in *Regularization Tools* [22], where $d$ is fixed with respect to the regularity of the solution. All the figures where we will have the regularization parameter $\gamma$ represented on the abscissa axis, we will use a logarithmic distribution of its values through the Matlab command "*logspace(⋯)*". Moreover, we displayed this abscissa axis on a log-scale.

### *4.1.* Example 1: One-dimensional image restoration model

An one-dimensional image reconstruction problem coming from the discretization of the Fredholm integral equation of the first kind, was investigated in [36] and defined by

$$\int_{-\frac{\pi}{2}}^{\frac{\pi}{2}} K(s,x)u^\star(x)dx = b(s), \quad -\frac{\pi}{2} \leq s \leq \frac{\pi}{2}, \tag{47}$$

where $K$ is the kernel, $b$ is the right hand side and $u^\star$ is the desired solution. This solution is the light intensity source defined by the following analytical expression

$$u^\star(x) = a_1 \exp\left(-c_1 (x - x_1)^2\right) + a_2 \exp\left(-c_2 (x - x_2)^2\right), \quad -\frac{\pi}{2} \leq x \leq \frac{\pi}{2}, \tag{48}$$

where parameters $a_1$, $a_2$, $c_1$, $c_2$, $x_1$ and $x_2$ are constants that determine the intensity, width, and position of the light sources. Throughout this numerical example, we consider two sets of these parameters whose the first one is given by

$$\begin{aligned} a_1 = 1,\ \ c_1 = 4,\ \ x_1 = 0.5 \\ a_2 = 1,\ \ c_2 = 4,\ \ x_2 = -0.5. \end{aligned} \tag{49}$$

Whereas the second set of parameters is defined as follows

$$\begin{aligned} a_1 = 2,\ \ c_1 = 6,\ \ x_1 = 0.8 \\ a_2 = 1,\ \ c_2 = 2,\ \ x_2 = -0.5. \end{aligned} \tag{50}$$

As for the exact kernel $K$, it is analytically defined by the following expression

$$K(s,x) = (\cos(s)+\cos(x))^2 \left(\frac{\sin(\psi(s,x))}{\psi(s,x)}\right)^2, \;\; -\frac{\pi}{2} \le s \le \frac{\pi}{2}, \tag{51}$$

with the given function $\psi(s,x) = \pi\,(\sin(s)+\cos(x))$.
To apply a quadrature formula on the integral equation of the first kind (47), we use as in [23, 24] the midpoint quadrature-collocation rule with $n$ equidistantly spaced points $s_i$ and $x_i$, defined as

$$s_i = -\frac{\pi}{2} + ih, \;\; x_i = -\frac{\pi}{2} + ih, \;\; h = \frac{\pi}{n}, \;\; i = 1, \cdots, n. \tag{52}$$

Rewriting the integral equation (47) through the discretizations in (52), we obtain the linear system

$$\sum_{j=1}^{n} w_j K(s_i, x_j^\star) u^\star(x_j^\star) = b(s_i), \;\; i = 1, \cdots, n, \tag{53}$$

where each midpoint is given by $x_j^\star = x_j - \frac{h}{2}$, and $w_j = h$ for $j = 1, \cdots, n$ are the weights of the quadrature formula. We thus derive the $n \times n$ matrix $\mathbf{A} = (a_{ij})_{1\le i,j\le n}$ and the $n$-vector exact solution $\boldsymbol{u}^\star = (u_j^\star)_{1\le j\le n}$ from the analytical formulas (51) and (48) as follows

$$a_{ij} = w_j K(s_i, x_j^\star), \;\; u_j^\star = u^\star(x_j^\star), \;\; i = 1, \cdots, n, \;\; j = 1, \cdots, n. \tag{54}$$

Since the exact solution in (48) has a continuous second derivative in $\left[-\frac{\pi}{2},\, \frac{\pi}{2}\right]$, we may thus consider as a reguarization matrix the following second derivative discrete operator of order $(n-2)\times n$

$$\mathbf{L} = \begin{bmatrix} 1 & -2 & 1 & 0 & \cdots & 0 \\ 0 & 1 & -2 & 1 & \ddots & \vdots \\ \vdots & \ddots & \ddots & \ddots & \ddots & 0 \\ 0 & \cdots & 0 & 1 & -2 & 1 \end{bmatrix} \in \mathbb{R}^{(n-2)\times n} \tag{55}$$

Here, $\mathbf{L}\boldsymbol{v}$ is a finite-difference approximation that is proportional to the second derivative of $\boldsymbol{v}$, and the penalty term $\|\mathbf{L}\boldsymbol{v} - \mathbf{g}\|_{n-2}$ when added in the functional $\mathcal{J}_\gamma(\cdot)$ in (19) means that the solution second derivative is controlled.

Throughout this numerical example, we use the Matlab codes *shaw(n)* and *get_l(n, 2)* implemented in *Regularization Tools* [22] to obtain the data $\mathbf{A}$, $\boldsymbol{b}$, $\mathbf{L}$ and $\boldsymbol{u}^\star$. We thus solve by a direct method the equation (22) by using the Matlab command \.

#### *4.1.1.* **The unperturbed data problem**

We assume here that the data (namely $\mathbf{A}$ and $\boldsymbol{b}$) from the discretization of the integral equation of the first kind (47) are not perturbed, despite the computational errors that may come from the approximation of this integral by a quadrature rule in the computing of the matrix $\mathbf{A}$. This integration error should be very small, when we consider very large dimensions $m$ and $n$. However, the matrix $\mathbf{A}$ should be very ill-conditioned and of rank-deficient. Moreover, we use the fact that the solution we're trying to approximate (the exact solution in (48)) is two times differentiable with a known derivative approximated by $\boldsymbol{g} = \mathbf{L}\boldsymbol{u}^\star$ (in real applications $\boldsymbol{g}$ may not be known).

For this first validating part, we fix the dimensions to $m = n = 64$ whose any change will be mentioned accordingly, and we use the parameters in (49). We thus face a huge conditioning of the matrix $\mathbf{cond}(\mathbf{A}) = 9.3099e^{+18}$ with a rank-deficient $\mathbf{rank}(\mathbf{A}) = 20$. We therefore present in Figure 2 - C and D, the absolute and relative errors between the exact data $\mathbf{A}\boldsymbol{u}^\star$ and the reconstructed data $\mathbf{A}\boldsymbol{u}_\gamma^\star$, with respect to the regularization parameter $\gamma \in [10^{-5}, 1]$ and $\gamma \in [1, 10^5]$. We also present in Figure 2 - A and B the absolute and relative errors between the exact analytical solution $\boldsymbol{u}^\star$ from (48) and the stabilized-regularized solution $\boldsymbol{u}_\gamma^\star$ from the solving of equation (22). These errors show how the stabilized-regularized approach we propose is accurate, stable (with respect to the regularization parameter $\gamma$) and efficient. Taking into account errors that come from the numerical integration scheme, we may accept that these results illustre very well Theorem 2 despite the rank-deficient faced. These numerical results are as good as those obtained through the multi-parameter regularization technique proposed in [8] (see table 3 and table 4).
Results in Figure 2 confirm the null regularization error (i.e. $\|\boldsymbol{u}^\star - \boldsymbol{u}_\gamma^\star\|_n = 0$) we gain from the formulation (22), under an unperturbed data problem. Moreover, the central problem of choosing the best regularization parameter seem also to be overcomed. We observe that from a certain threshold $\gamma_0$, the stabilized-regularized solution $\boldsymbol{u}_\gamma^\star$ becomes a stable and accurate approximation for any $\gamma > \gamma_0$.

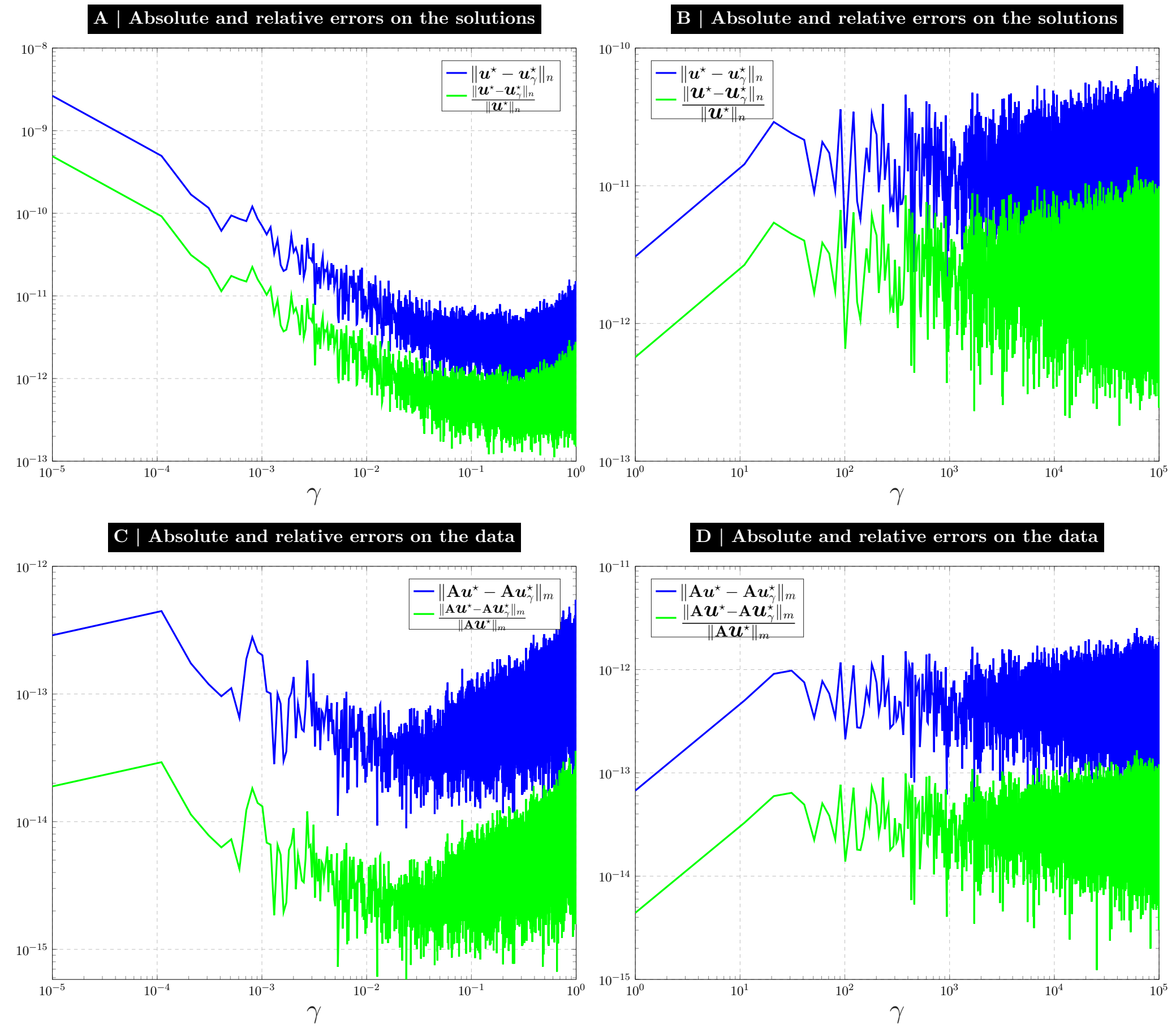

Figure 2: Absolute and relative errors between the exact solution $\boldsymbol{u}^\star$ computed from (48) and the stabilized-regularized solution $\boldsymbol{u}_\gamma^\star$ obtained from the solving of (22) by the Matlab command \, with respect to the regularization parameter $\gamma \in [10^{-5}, 1]$ (in A and C) and $\gamma \in [1, 10^5]$ (in B and D). The data $\mathbf{A}$, $\mathbf{L}$ and $\boldsymbol{b}$ are computed from the Matlab codes *shaw(n)* and *get_l(n, 2)* provided in *Regularization Tools* [22] under the dimensions $m = n = 64$ ($\mathbf{cond}(\mathbf{A}) = 9.3099e^{+18}$, $\mathbf{rank}(\mathbf{A}) = 20$), and the derivative is approximated as $\boldsymbol{g} = \mathbf{L}\boldsymbol{u}^\star$. Parameters in (49) are used, and the errors are displayed on a log-log scale.

In order to deepen these points, we test in the following subsection the formulation (22) under the context of perturbed data as in [23, 25] where the only *type of prior information* we have on the *unknown* exact solution $\boldsymbol{u}^\star$ is that, it is smooth in the sense defined by the operator in (55) (but we do not know the result of the smoothness, i.e. the function $\boldsymbol{g} = \mathbf{0}$). This means that in the formulation (22), the right hand term $\gamma \mathbf{L}^t \boldsymbol{g}$ should not been taken into account since being considered as a non available information.

### *4.1.2.* Perturbations by white noises

We first investigate the case where the right hand side $\boldsymbol{b} = (b_i)_{1 \le i \le n}$ is perturbed by a white noise, that is an uncorrelated normally distributed error $\eta\varepsilon_i$, $i = 1, \cdots, n$, with a zero mean and the deviation $\eta = 10^{-3}$ (any change of this deviation will be mentioned accordingly)

$$\tilde{b}_i = b_i + \eta\varepsilon_i, \;\text{ where }\; \varepsilon_i \sim \mathcal{N}(0,\, 1). \tag{56}$$

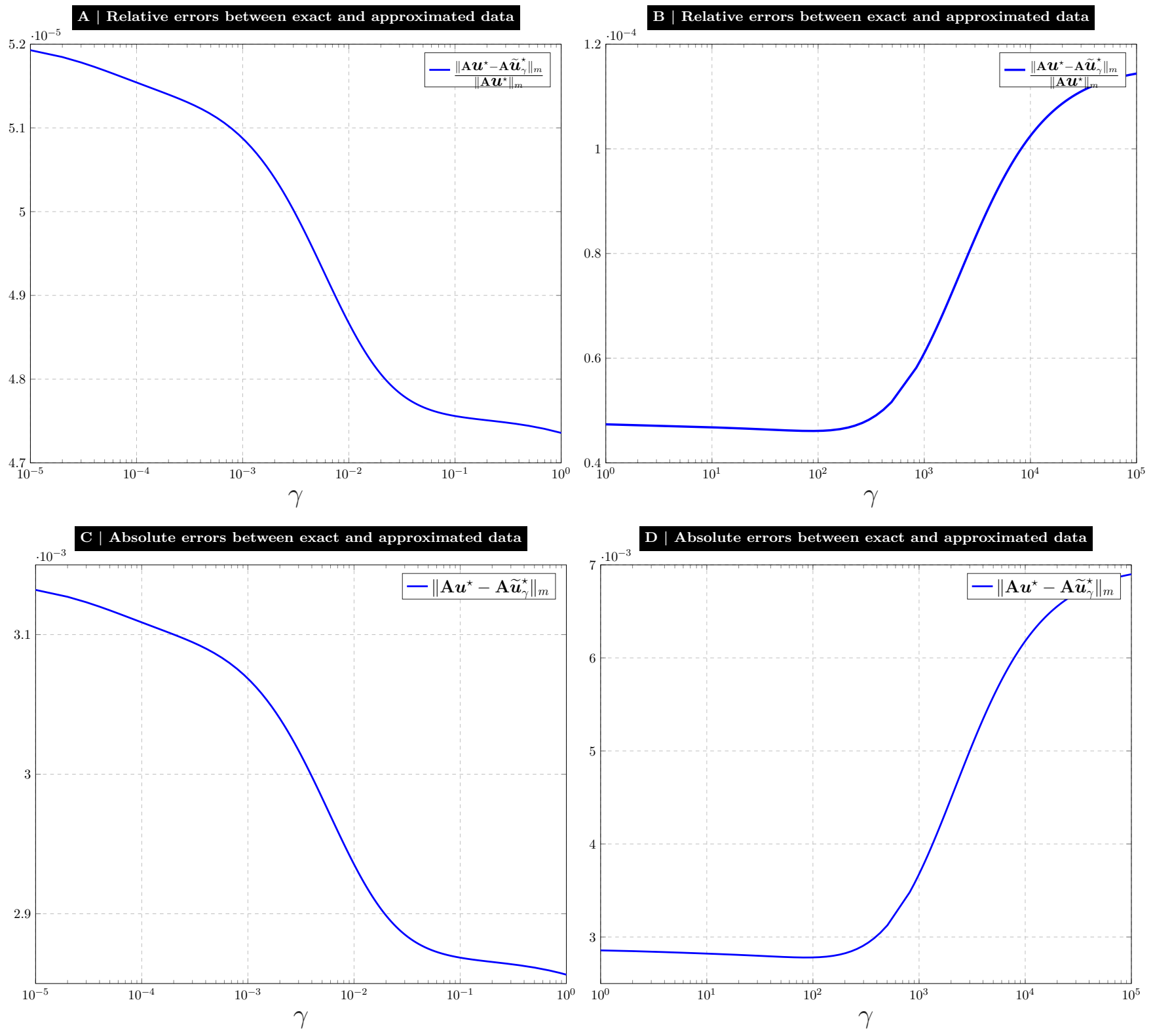

Figure 3: Absolute and relative errors between the exact data $\mathbf{A}\boldsymbol{u}^\star$ computed from (48), and the stabilized-regularized data $\mathbf{A}\widetilde{\boldsymbol{u}}_\gamma^\star$ obtained from the solving of (22) by the Matlab command \, with respect to the parameter $\gamma \in [10^{-5}, 1]$ (in A and C) and $\gamma \in [1, 10^5]$ (in B and D). The data $\mathbf{A}$, $\mathbf{L}$ and $\boldsymbol{b}$ are computed from the Matlab codes *shaw(n)* and *get_l*$(n, 2)$ provided in *Regularization Tools* [22] under the dimensions $m = n = 1\,000$ ($\mathbf{cond}(\mathbf{A}) = 7.697e^{+20}$, $\mathbf{rank}(\mathbf{A}) = 20$), and the right hand term $\gamma \mathbf{L}^t \boldsymbol{g}$ in (22) is not taken into account since $\boldsymbol{g} = \mathbf{0}$. The computations are made under the perturbed data in (56) and the parameters in (49) are used. The horizontal axis is displayed on a log-scale.

In this part, we fix $\boldsymbol{g} = \mathbf{0}$ and in order to reduce the quadrature error induced in the approximation of the integral in (47), we set a fairly high number of integration points to $m = n = 1\,000$. We hence face a poorly conditioned matrix of $\mathbf{cond}(\mathbf{A}) = 7.697e^{+20}$ whose the rank is also very small $\mathbf{rank}(\mathbf{A}) = 20$.

We present in Figure 3 relative and absolute errors between the exact data $\mathbf{A}\boldsymbol{u}^\star$ (where $\boldsymbol{u}^\star$ is computed from (48)) and the stabilized-regularized data $\mathbf{A}\widetilde{\boldsymbol{u}}_\gamma^\star$ (where $\widetilde{\boldsymbol{u}}_\gamma^\star$ is obtained from the solving of the system (22)), under the perturbation in (56) and the parameters in (49). These errors being computed with respect to the regularization parameter $\gamma \in [10^{-5}, 1]$ (for A and C) and $\gamma \in [1, 10^5]$ (for B and D), show how much the stabilized-regularized approach in (22) provides accurate and stable data that do not change much with respect to the parameter $\gamma$ (hence does not depend on a predefined optimal parameter). The results in C and D are coherent with the result (45) from Corollary 1.

We also present in Figure 4 the corresponding relative and absolute errors between the exact solution $\boldsymbol{u}^\star$ and the stabilized-regularized solution $\widetilde{\boldsymbol{u}}_\gamma^\star$, under the perturbed right hand side provided in (56) and the parameters in (49). Despite the fact that we face poor conditioning and rank-deficient matrix, we observe the result in D that is not only consistent with the one in (46) of Corollary 1, and is also stable especially when the regularization parameter $\gamma > 1$ (that is to say when there is over-regularization).

We display in Figure 5 the exact data $\mathbf{A}\boldsymbol{u}^\star$, the data reconstructed from the stabilized-regularized approach $\mathbf{A}\widetilde{\boldsymbol{u}}_\gamma^\star$, and the ones from Tikhonov regularization $\mathbf{A}\widetilde{\boldsymbol{u}}_{\text{tik},\gamma}^\star$ given in (5) (with $\boldsymbol{g} = \mathbf{0}$), for different values of the regularization parameter $\gamma = 10^{-5}, 1, 10^5, 10^{10}$. We observe that the curves of these three data are overlaid for the values of $\gamma = 10^{-5}$ and $\gamma = 1$. For larger values of the regularization parameter $\gamma$ however (for instance $\gamma = 10^5, 10^{10}$), we note that only the data reconstructed from the stabilized-regularized approach remain equal to the exact data. The data from the Tikhonov regularization behave differently, when facing large regularization parameters (see also Figure 2 in [23]).

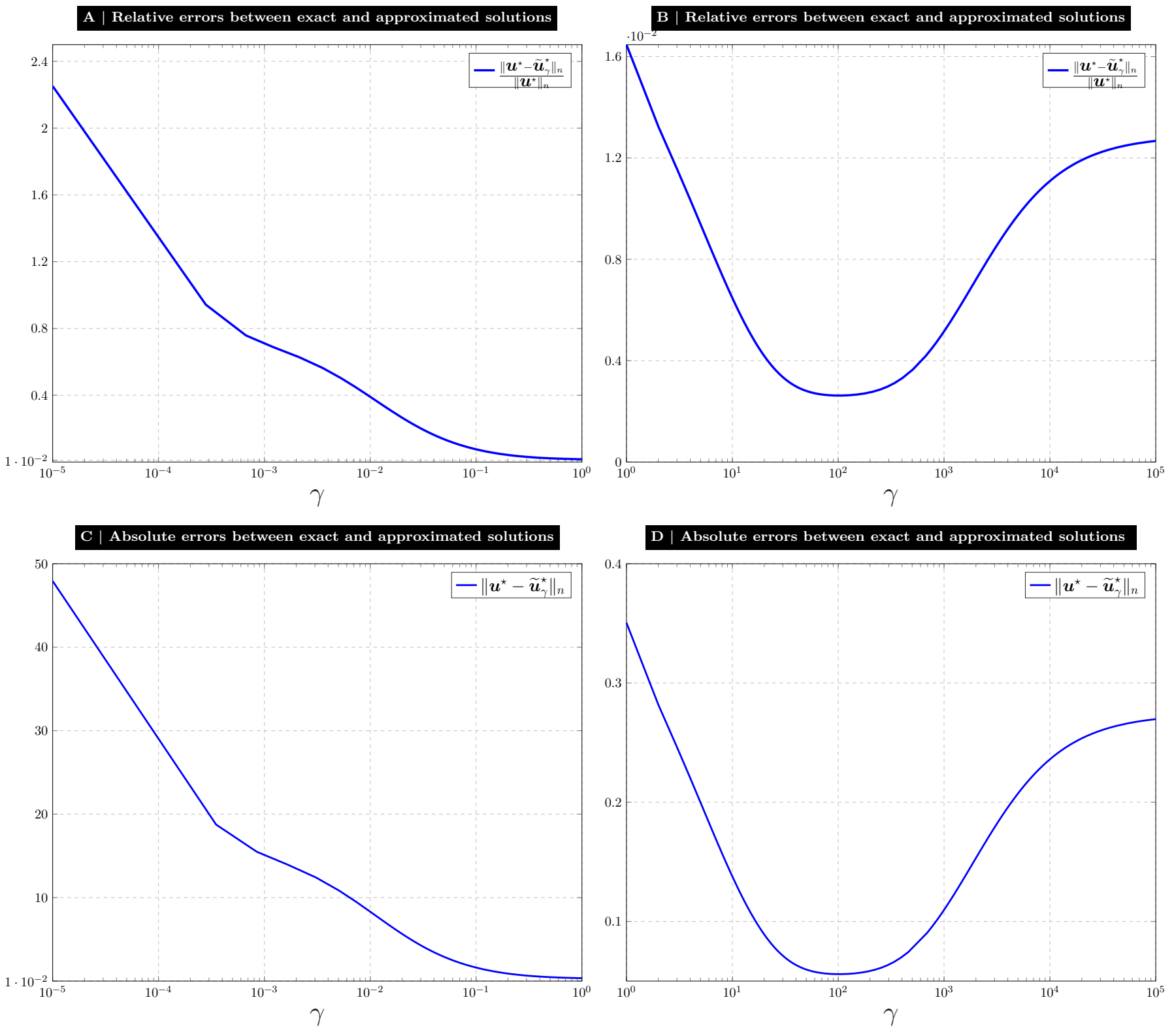

Figure 4: Absolute and relative errors between the exact solution $\boldsymbol{u}^\star$ computed from (48), and the stabilized-regularized solution $\widetilde{\boldsymbol{u}}^\star_\gamma$ obtained from the solving of (22) by the Matlab command \, with respect to the parameter $\gamma \in [10^{-5}, 1]$ (in A and C) and $\gamma \in [1, 10^5]$ (in B and D). The data $\mathbf{A}$, $\mathbf{L}$ and $\boldsymbol{b}$ are computed from the Matlab codes *shaw(n)* and *get_l*$(n, 2)$ provided in *Regularization Tools* [22] under the dimensions $m = n = 1\,000$ ($\mathbf{cond}(\mathbf{A}) = 7.697e^{+20}$, $\mathbf{rank}(\mathbf{A}) = 20$), and the right hand term $\gamma \mathbf{L}^t \boldsymbol{g}$ in (22) is not taken into account since $\boldsymbol{g} = \mathbf{0}$. The computations are made under the perturbed data in (56) and the parameters in (49). The horizontal axis is displayed on a log-scale.

We also display in Figure 6, the corresponding solutions $\boldsymbol{u}^\star$, $\widetilde{\boldsymbol{u}}^\star_\gamma$ and $\widetilde{\boldsymbol{u}}^\star_{\text{tik},\gamma}$ in order to graphically illustrate the stability of stabilized-regularized approach (22), despite that the matrix $\mathbf{A}$ is rank-deficient. We thus present for different values of the regularization parameter $\gamma = 10^{-5}, 1, 10^5, 10^{10}$, these solutions with respect to the index $i = 1, \cdots, 1\,000$ of the domain's mesh. The Tikhonov regularization solution was already investigated in [23] (see also sections 4.8.1 and 4.8.2 in [24]), where the appropriated value of the regularization parameter was chosen through the L-curve procedure (we can also see Figures $1, 2, 3, 4, 5$ in [36] for others approximations). Here for $\gamma = 10^{-5}$ in Figure 6, both solutions obtained from the stabilized-regularized formulation (22) and from the Tikhonov formulation (5) are oscillating. Whereas for $\gamma = 1$, we observe both solutions well fitting to the exact solution. However, for large values of the regularization parameter $\gamma$ (for $\gamma = 10^5$, $10^{10}$), only the stabilized-regularized solution remains fitted to the exact solution even if we observe small shift at the boundary of the domain. We thus obtain precise and stable approximation of the non-linear exact solution (48), under an over-regularization (i.e. for any $\gamma \geq 1$).

We present in Figure 7, the associated L-curves of our stabilized-regularized method, that is to say a representation of the size of the stabilized-regularized solution (through its semi-norm $\|\mathbf{L}\widetilde{\boldsymbol{u}}^\star_\gamma\|_{n-2}$ and norm $\|\widetilde{\boldsymbol{u}}^\star_\gamma\|_n$) versus the corresponding residual $\|\mathbf{A}\widetilde{\boldsymbol{u}}^\star_\gamma - \widetilde{\boldsymbol{b}}\|_m$. We thus note very small variations of these semi-norm and norm with respect to very small variations of the residual (compared to Figure 1 in [23] and Figure 1 in [25]), despite the wide range of values of the regularization parameter chosen $\gamma \in \left[1, 10^5\right]$. This is an additional proof of the stability of the formulation proposed here, and which no longer requires the graphical analysis of the L-curve method as it is proposed in [23, 25] (or any other strategy for the parameter's choice).

For a more investigation of the stability of the stabilized-regularized approach proposed in (22),

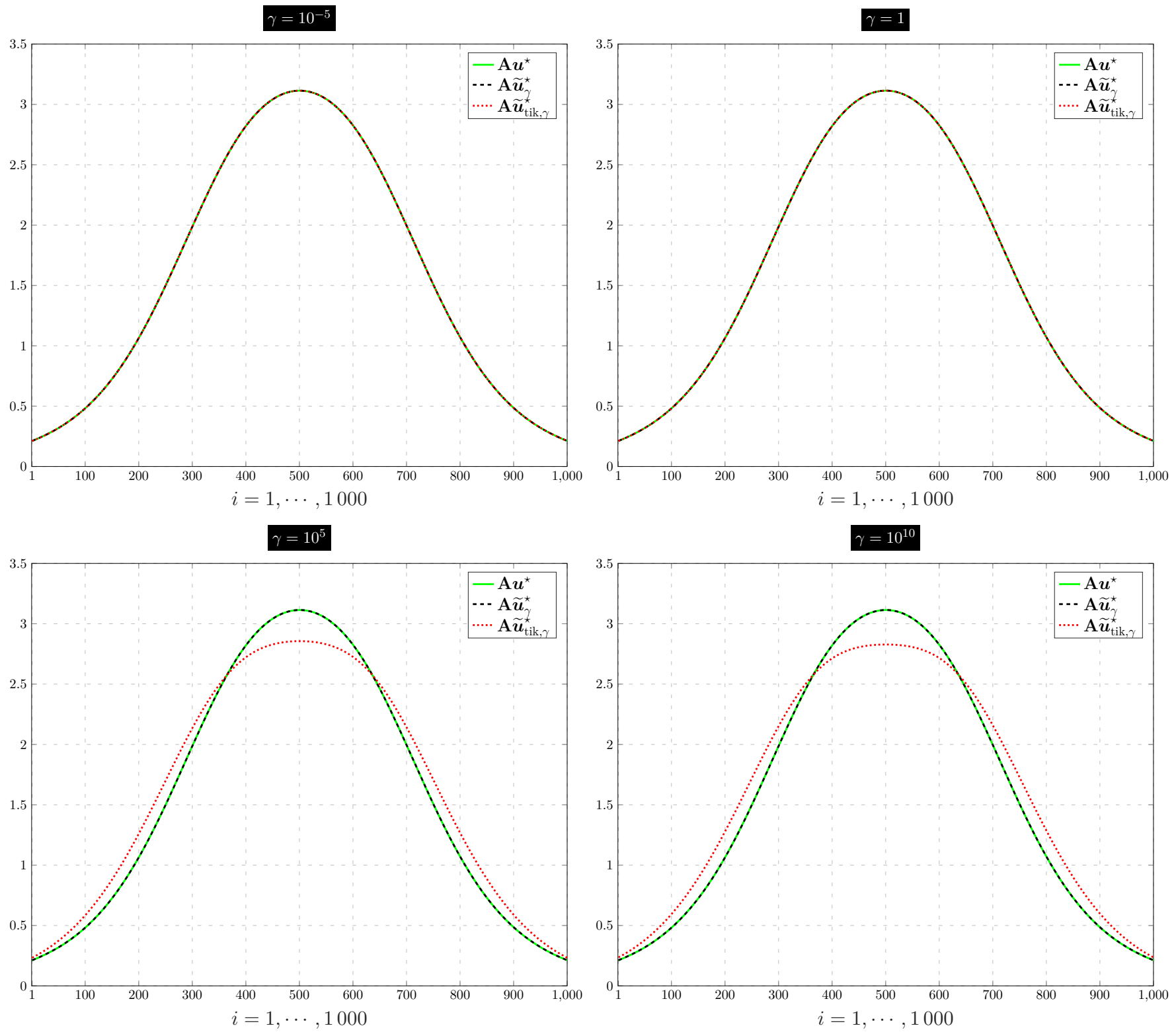

Figure 5: The exact data $\mathbf{A}\boldsymbol{u}^\star$ (in green solid line), the data reconstructed from the stabilized-regularized approach $\mathbf{A}\widetilde{\boldsymbol{u}}^\star_\gamma$ (in black dashed line), and the data made from the Tikhonov regularization $\mathbf{A}\widetilde{\boldsymbol{u}}^\star_{\mathrm{tik},\gamma}$ (in red dotted line), are presented with respect to the discretisation index $i = 1, \cdots, 1\,000$, for the different values of the regularization parameter $\gamma = 10^{-5}, 1, 10^5, 10^{10}$ (given in the top of the plots). The stabilized-regularized solution $\widetilde{\boldsymbol{u}}^\star_\gamma$ is obtained from the solving of (22) through the Matlab command \, and the Tikhonov regularization solution $\widetilde{\boldsymbol{u}}^\star_{\mathrm{tik},\gamma}$ is computed in the solving of (5) by the Matlab code *tikhonov( )* implemented in *Regularization Tools* [22]. The data $\mathbf{A}$, $\mathbf{L}$ and $\boldsymbol{b}$ are computed from the Matlab codes *shaw(n)* and *get_l*$(n, 2)$ in *Regularization Tools* [22] and the approximate derivative is fixed to $\boldsymbol{g} = \mathbf{0}$. The computations are made under the perturbed data in (56), with respect to the dimensions $m = n = 1\,000$ ($\mathbf{cond}(\mathbf{A}) = 7.697e^{+20}$, $\mathbf{rank}(\mathbf{A}) = 20$) and the parameters in (49).

we have also analysed the modified shaw problem by replacing the parameters in (49) by those in (50). We have used here two different noise levels, namely $\eta = 10^{-3}$ and $\eta = 10^{-2}$.
For the first noise level $\eta = 10^{-3}$, we present in Figure 8 the exact solution $\boldsymbol{u}^\star$ computed from (48), the stabilized-regularized solution $\widetilde{\boldsymbol{u}}^\star_\gamma$ obtained from the solving of (22) through the Matlab command \, and the Tikhonov regularization solution $\widetilde{\boldsymbol{u}}^\star_{\mathrm{tik},\gamma}$ obtained in the solving of (5) by the Matlab code *tikhonov( )*, with respect to the discretization index $i = 1, \cdots, 1\,000$, for the different values of the regularization parameter $\gamma = 1, 10^3, 10^6, 10^{12}$ (these values are among the ones used in the computing of the filter factors in Figure 1). For small regularization parameters $\gamma = 1, 10^3$, solutions computed from the Tikhonov method and the stabilized-regularized method are very similar, which is coherent with the behaviour of the filter factors in Figure 1. However for larger regularization parameters, namely $\gamma = 10^6, 10^{12}$, the solution from the stabilized-regularized method keeps still trying to capte the two bumps of the exact solution but with big deviations at the end of the interval. These deviations are simply due to the relatively small dimensions. Indeed in Figure 9, by successively taking the dimensions $m = n = 10\,000$, $m = n = 20\,000$, $m = n = 25\,000$ and $m = n = 30\,000$, we observe an increasingly exact approximation solution $\widetilde{\boldsymbol{u}}^\star_\gamma$.
As for the noise level $\eta = 10^{-2}$, we take back the same computations in Figure 10 and Figure 11, and we observe the same behavior. These tests bring the same accuracy and stability.

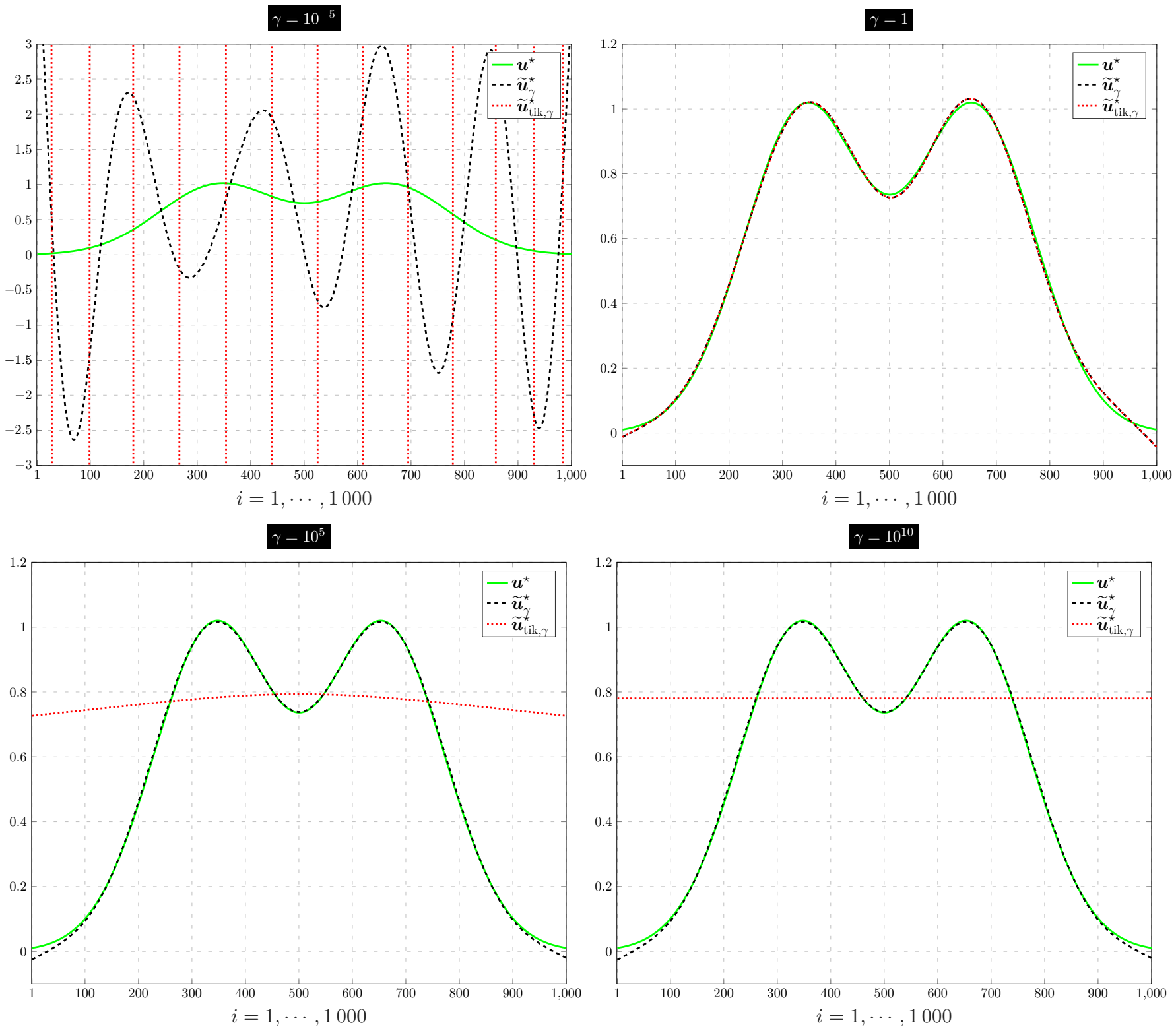

Figure 6: The exact solution $\boldsymbol{u}^\star$ (in green solid line) computed from (48), the stabilized-regularized solution $\widetilde{\boldsymbol{u}}^\star_\gamma$ (in black dashed line) obtained from the solving of (22) by the Matlab command \, and the Tikhonov regularization solution $\widetilde{\boldsymbol{u}}^\star_{\mathrm{tik},\gamma}$ (in red dotted line) obtained in the solving of (5) by the Matlab code *tikhonov( )* implemented in *Regularization Tools* [22], are presented with respect to the discretization index $i = 1, \cdots, 1\,000$, for different values of the parameter $\gamma = 10^{-5}, 1, 10^5, 10^{10}$. The data $\mathbf{A}$, $\mathbf{L}$ and $\boldsymbol{b}$ are computed from the Matlab codes *shaw(n)* and *get_l*$(n, 2)$ in *Regularization Tools* [22] and the term $\gamma \mathbf{L}^t \boldsymbol{g}$ in (5) and (22) is not taken into account ($\boldsymbol{g} = \mathbf{0}$). The computations are made under the perturbed data in (56), with the dimensions $m = n = 1\,000$ ($\mathbf{cond}(\mathbf{A}) = 7.697e^{+20}$, $\mathbf{rank}(\mathbf{A}) = 20$) and the parameters in (49).

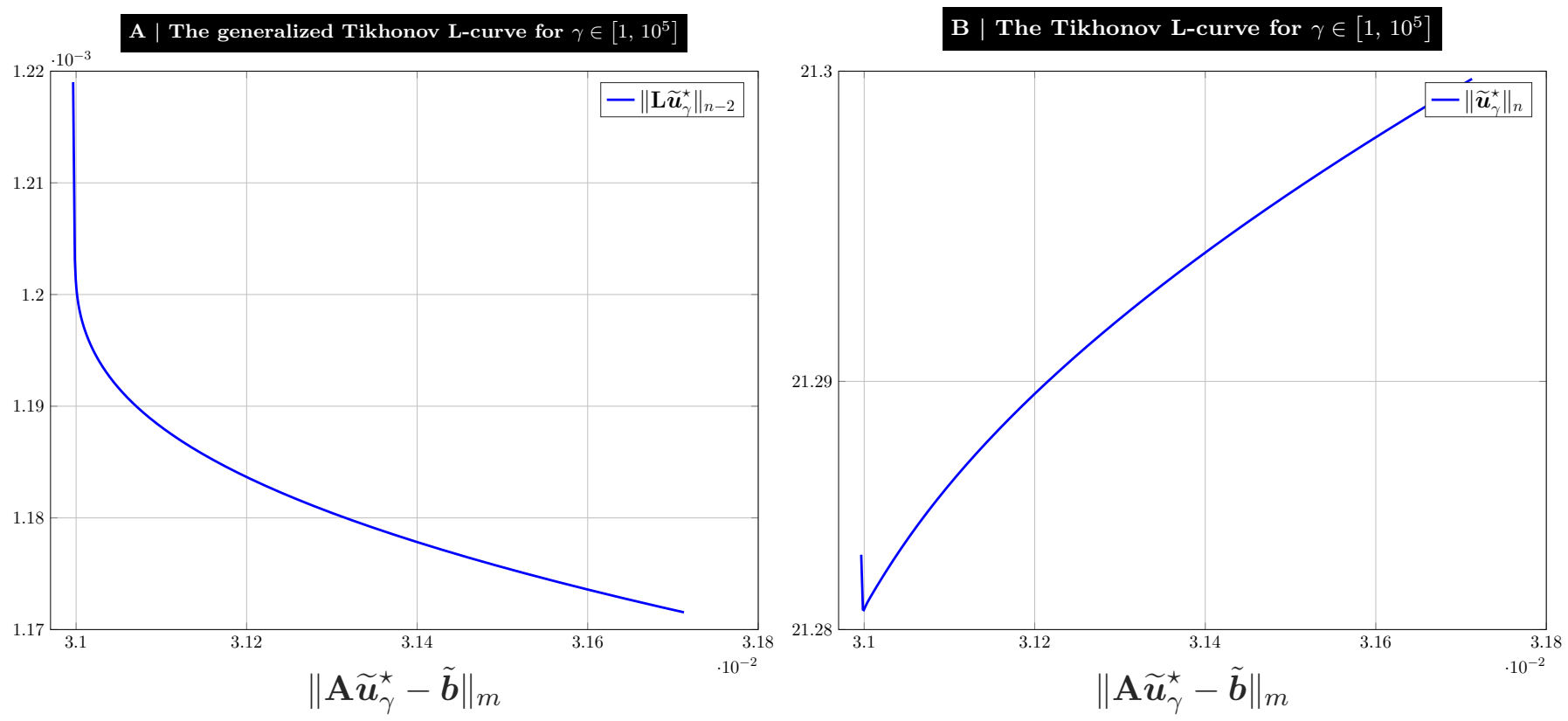

Figure 7: The residual norm $\|\mathbf{A}\widetilde{\boldsymbol{u}}^\star_\gamma - \widetilde{\boldsymbol{b}}\|_m$ versus the semi-norm $\|\mathbf{L}\widetilde{\boldsymbol{u}}^\star_\gamma\|_{n-2}$ and the norm $\|\widetilde{\boldsymbol{u}}^\star_\gamma\|_n$ represented for the range of values of the parameter $\gamma \in \left[1, 10^5\right]$. The data $\mathbf{A}$, $\mathbf{L}$ and $\boldsymbol{b}$ are computed from the codes *shaw(n)* and *get_l*$(n, 2)$ provided in *Regularization Tools* [22]. The stabilized-regularized solution $\boldsymbol{u}^\star_\gamma$ obtained from the solving of (22) by the Matlab command \, without taken into account the right hand term $\gamma \mathbf{L}^t \boldsymbol{g}$ in (22) ($\boldsymbol{g} = \mathbf{0}$). The computations are made under the perturbed data in (56), with the parameters in (49) and the dimensions $m = n = 1\,000$ ($\mathbf{cond}(\mathbf{A}) = 7.697e^{+20}$, $\mathbf{rank}(\mathbf{A}) = 20$).

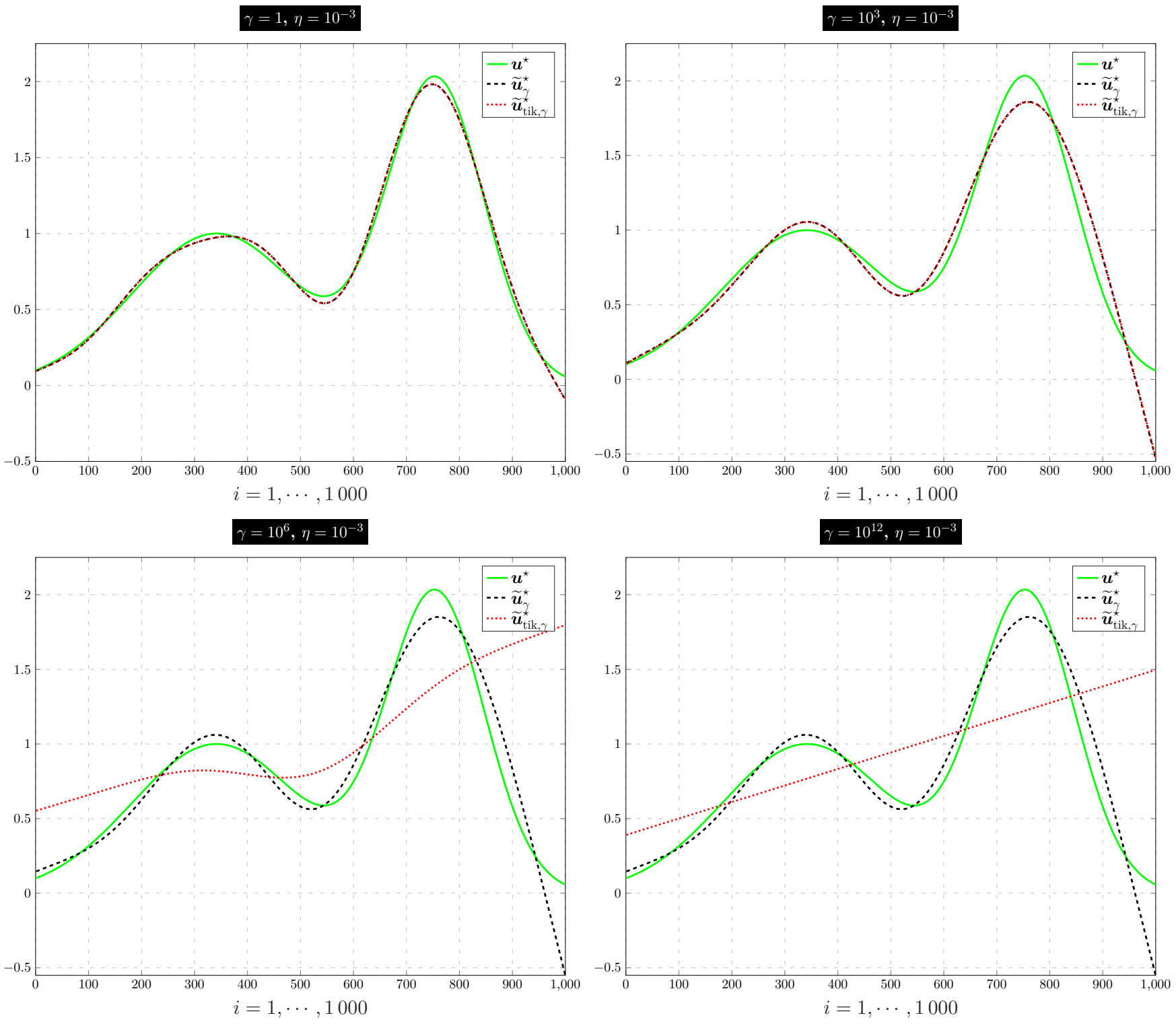

Figure 8: The exact solution $\boldsymbol{u}^\star$ (in green solid line) computed from (48), the stabilized-regularized solution $\widetilde{\boldsymbol{u}}_\gamma^\star$ (in black dashed line) obtained from the solving of (22) through the Matlab command \, and the Tikhonov regularization solution $\widetilde{\boldsymbol{u}}_{\text{tik},\gamma}^\star$ (in red dotted line) obtained in the solving of (5) by the Matlab code *tikhonov( )* in *Regularization Tools* [22], presented with respect to the discretization index $i = 1, \cdots, 1\,000$, for different values of the parameter $\gamma = 1, 10^3, 10^6, 10^{12}$. The data $\mathbf{A}$, $\mathbf{L}$ and $\boldsymbol{b}$ are computed from the Matlab codes *shaw(n)* and *get_l*$(n, 2)$ in *Regularization Tools* [22] and the term $\gamma \mathbf{L}^t \boldsymbol{g}$ in (5) and (22) is not taken into account ($\boldsymbol{g} = \mathbf{0}$). Computations are made under perturbation in (56) with deviation $\eta = 10^{-3}$. The dimensions are fixed to $m = n = 1\,000$ ($\mathbf{cond}(\mathbf{A}) = 7.697e^{+20}$, $\mathbf{rank}(\mathbf{A}) = 20$), and the parameters in (50) are used.

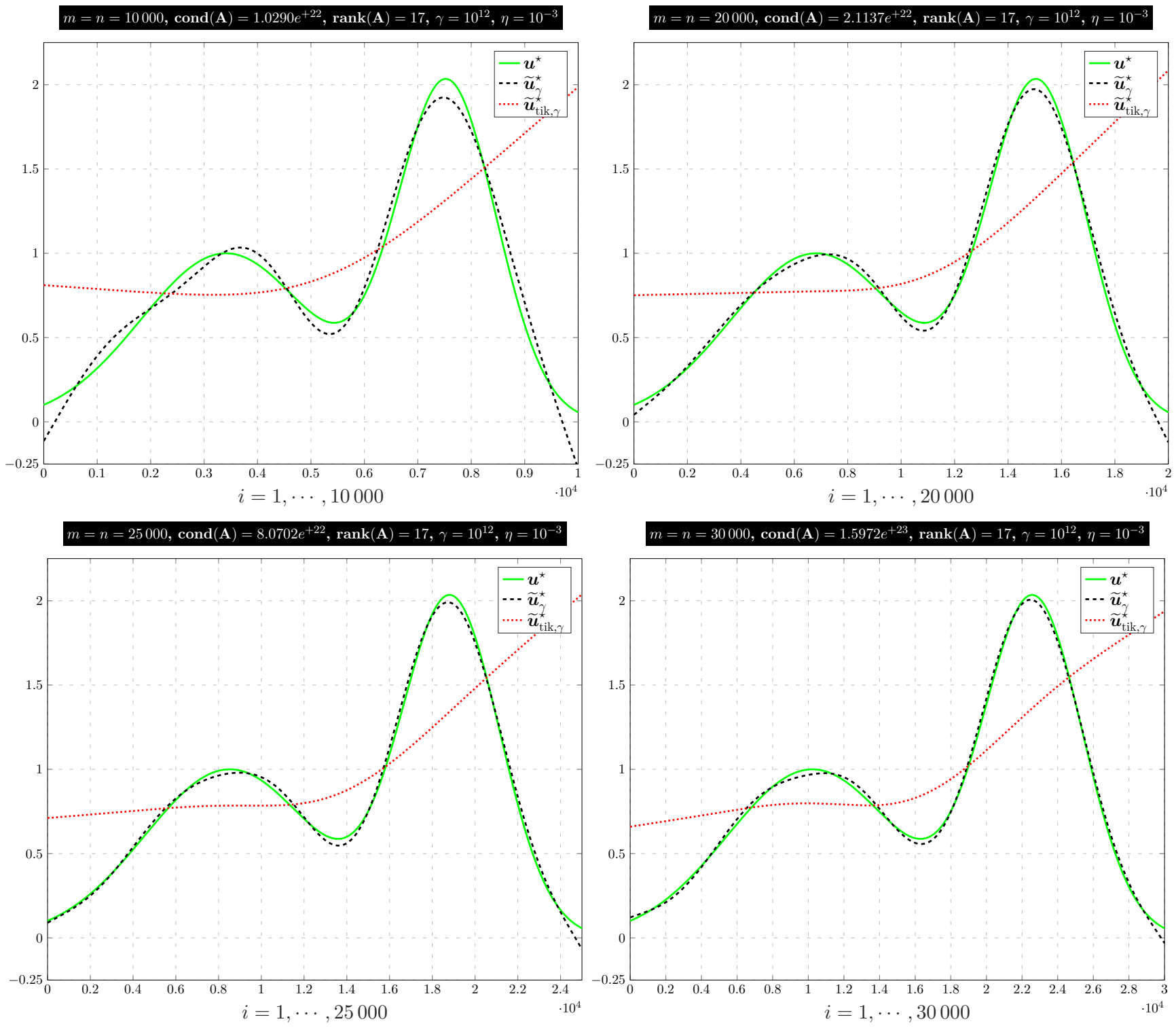

Figure 9: The exact solution $\boldsymbol{u}^\star$ (in green solid line) computed from (48), the stabilized-regularized solution $\widetilde{\boldsymbol{u}}_\gamma^\star$ (in black dashed line) obtained from the solving of (22) through the Matlab command \, and the Tikhonov regularization solution $\widetilde{\boldsymbol{u}}_{\text{tik},\gamma}^\star$ (in red dotted line) obtained in the solving of (5) by the Matlab code *tikhonov( )* implemented in *Regularization Tools* [22], presented with respect to the discretization index $i = 1, \cdots, m$, for different values of the dimensions $m = n = 10\,000, 20\,000, 25\,000, 30\,000$. The data $\mathbf{A}$, $\mathbf{L}$ and $\boldsymbol{b}$ are computed from the Matlab codes *shaw(n)* and *get_l*$(n, 2)$ in *Regularization Tools* [22] and the term $\gamma \mathbf{L}^t \boldsymbol{g}$ in (5) and (22) is not taken into account ($\boldsymbol{g} = \mathbf{0}$). Computations are made under perturbed data in (56) with the deviation $\eta = 10^{-3}$. The regularization parameter is fixed to $\gamma = 10^{12}$, and the parameters in (50) are used.

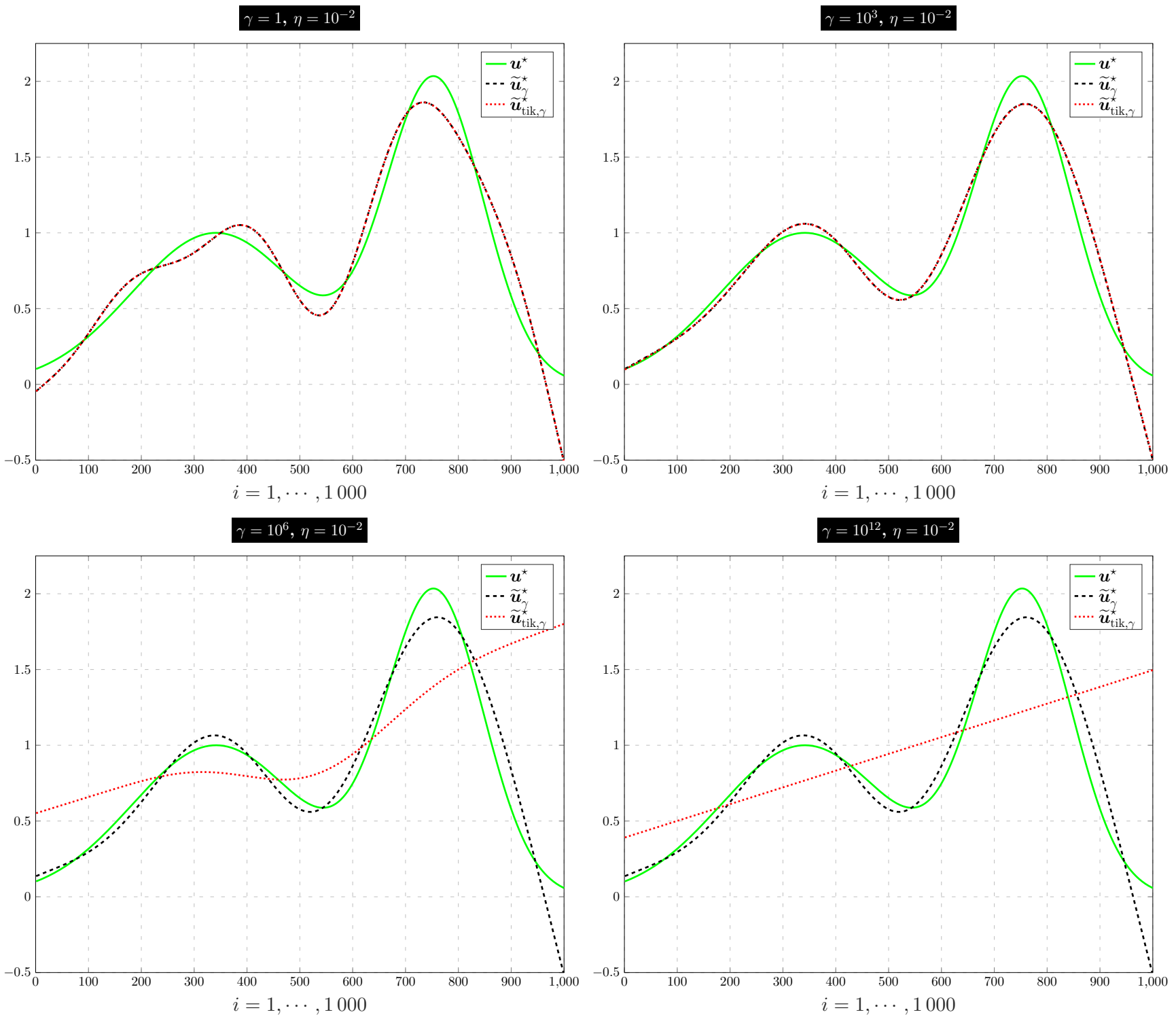

Figure 10: The exact solution $\boldsymbol{u}^\star$ (in green solid line) computed from (48), the stabilized-regularized solution $\widetilde{\boldsymbol{u}}^\star_\gamma$ (in black dashed line) obtained from the solving of (22) through the Matlab command \, and the Tikhonov regularization solution $\widetilde{\boldsymbol{u}}^\star_{\mathrm{tik},\gamma}$ (in red dotted line) obtained in the solving of (5) by the Matlab code *tikhonov( )* in *Regularization Tools* [22], are presented with respect to the discretization index $i = 1, \cdots, 1\,000$, for different values of the regularization parameter $\gamma = 1, 10^3, 10^6, 10^{12}$. The data $\mathbf{A}$, $\mathbf{L}$ and $\boldsymbol{b}$ are computed from the Matlab codes *shaw(n)* and *get_l*$(n, 2)$ in *Regularization Tools* [22] and the term $\gamma \mathbf{L}^t \boldsymbol{g}$ in (5) and (22) is not taken into account ($\boldsymbol{g} = \mathbf{0}$). The computations are made under the perturbed data in (56) with the more deviation $\eta = 10^{-2}$. The dimensions are fixed to $m = n = 1\,000$ ($\mathbf{cond}(\mathbf{A}) = 7.697e^{+20}$, $\mathbf{rank}(\mathbf{A}) = 20$), and the parameters in (50) are used.

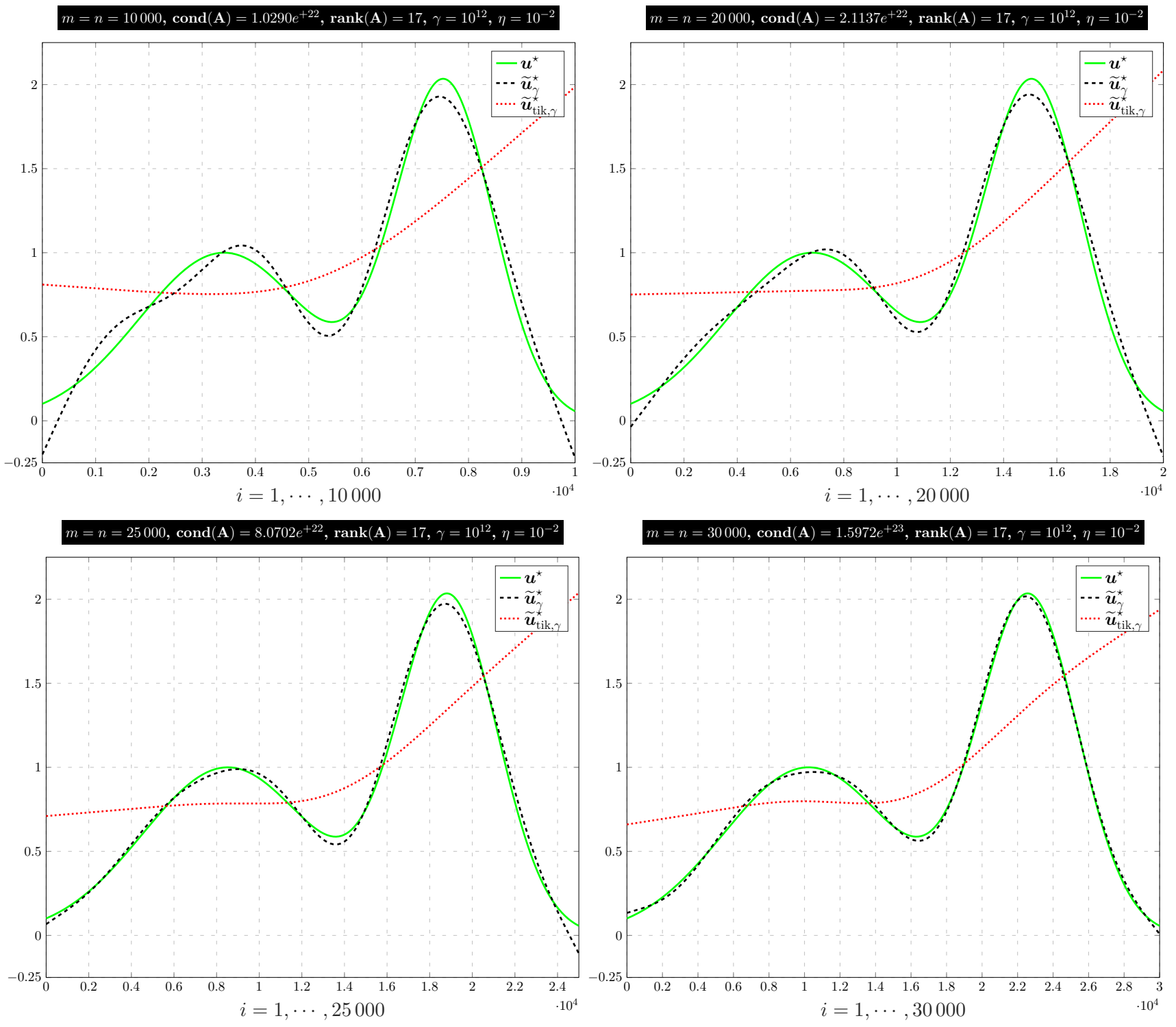

Figure 11: The exact solution $\boldsymbol{u}^\star$ (in green solid line) computed from (48), the stabilized-regularized solution $\widetilde{\boldsymbol{u}}^\star_\gamma$ (in black dashed line) obtained from the solving of (22) through the Matlab command \, and the Tikhonov regularization solution $\widetilde{\boldsymbol{u}}^\star_{\mathrm{tik},\gamma}$ (in red dotted line) obtained in the solving of (5) by the Matlab code *tikhonov( )* in *Regularization Tools* [22], are represented with respect to the discretization index $i = 1, \cdots, m$, for different values of the dimensions $m = n = 10\,000, 20\,000, 25\,000, 30\,000$. The data $\mathbf{A}$, $\mathbf{L}$ and $\boldsymbol{b}$ are computed from the Matlab codes *shaw(n)* and *get_l*$(n, 2)$ in *Regularization Tools* [22] and the term $\gamma \mathbf{L}^t \boldsymbol{g}$ in (5) and (22) is not taken into account ($\boldsymbol{g} = \mathbf{0}$). The computations are made under the perturbed data in (56) with the more deviation $\eta = 10^{-2}$. The regularization parameter is fixed to $\gamma = 10^{12}$, and the parameters in (50) are used.

#### *4.1.3.* Perturbations by correlated noises

The stabilized-regularized approach we propose here also deals very well with correlated noises (whether it is white noises or not). We consider here the following two highly correlated perturbations investigated in [23, 25]:

- **Filtered white noises**

This perturbation is defined from the filtered white noises $e_i$, generated by the following formula

$$e_i = ae_{i-1} + \eta\varepsilon_i, \text{ where } \varepsilon_i \sim \mathcal{N}(0,1), \ \eta = 10^{-3}, \ 0 \le a \le 1, \quad i = 2, \cdots, n.$$

Here, we fix the parameter $a = 0.5$ and both the matrix $\mathbf{A}$ and the right hand term $\boldsymbol{b}$ of the least squares problem are perturbed and become $\widetilde{\mathbf{A}} = (\tilde{a}_{ij})_{1\le i,j\le n}$ and $\widetilde{\boldsymbol{b}} = (\tilde{b}_i)_{1\le i\le n}$, where

$$\tilde{a}_{ij} = a_{ij} + e_i \text{ and } \tilde{b}_i = b_i + e_i, \quad i = 1, \cdots, n, \ j = 1, \cdots, n. \tag{57}$$

Under the perturbation (57), we have also made similar computations than the previous case and have presented the results in the Annexe 6.1.1 (see Figure 21, Figure 22, Figure 23 and Figure 24). These numerical results also align well with the theoretical results in (45) - (46) and show the same robustness as the previous white noise perturbation. In addition of being effective on the approximation of the non-linear solution $\boldsymbol{u}^\star$, we also notice here that our approach is also very precise and stable under the *filtered white noises* correlated data (57).

- **Perturbations from smoothing of the matrix A and the right hand term *b***

We also investigate perturbations that we denote $e_i^b$ and $e_{ij}^a$ and that are defined from regular smoothing of the matrix $\mathbf{A}$ and the vector $\boldsymbol{b}$. Indeed, we investigate here the following perturbations introduced in [23, 25] by

$$\begin{aligned} e_i^b &= \mu\left(b_{i-1} + b_{i+1}\right), \ i = 2, \cdots, n-1, \\ e_{ij}^a &= \mu\left(a_{i-1,j} + a_{i+1,j} + a_{i,j-1} + a_{i,j+1}\right), \ i = 2, \cdots, n-1, \ j = 2, \cdots, n-1, \end{aligned}$$

whete $\mu$ is a smoothing parameter that has to be fixed (that we set here to $\mu = 0.01$). Thus, the perturbed matrix $\widetilde{\mathbf{A}} = (\tilde{a}_{ij})_{1\le i,j\le n}$ and right hand side $\widetilde{\boldsymbol{b}} = (\tilde{b}_i)_{1\le i\le n}$ are given as follows

$$\tilde{a}_{ij} = a_{ij} + e_{ij}^a \text{ and } \tilde{b}_i = b_i + e_i^b, \ i = 1, \cdots, n, \ j = 1, \cdots, n. \tag{58}$$

In a practical point of view, the type of errors in (58) may represent sampling errors or even errors coming from the finite element approximation of partial differential equations models.

We present in the Annexe 6.1.2 at the Figure 25, the relative errors between the exact solution $\boldsymbol{u}^\star$ and the stabilized-regularized solution $\widetilde{\boldsymbol{u}}_\gamma^\star$, under the perturbations in (58). These errors also show how the stabilized-regularized solution is a precise and stable approximation of the non-linear exact solution. We also represent in the same Figure 25 the associated L-curves, under the range of the regularization parameter $\gamma \in [1, 10^5]$. Moreover, we display at Figure 26 the solutions $\boldsymbol{u}^\star$, $\widetilde{\boldsymbol{u}}_\gamma^\star$ and $\widetilde{\boldsymbol{u}}_{\mathrm{tik},\gamma}^\star$ with respect to the index $i = 1, \cdots, 1\,000$ of the domain's mesh, for different values of the regularization parameter $\gamma = 10^{-4}, 1, 10^5, 10^{10}$. These curves graphically illustrate the stability of our approach when $\gamma > 1$, despite that the matrix $\mathbf{A}$ is rank-deficient and poorly conditioned. Even by taking a huge regularization parameter (over-reguralizing), namely $\gamma = 10^5, 10^{10}$, we still have very accurate and stable solution.

### 4.2. **Example 2**: *Inverse heat equation*

The inverse heat equation that we investigate here is the underlying ill-posed problem of the Volterra integral equation of the first kind, formulated as follows

$$\int_0^s K(s,t)u^\star(t)dt = b(s), \quad 0 \le s \le 1, \tag{59}$$

with a convolution type kernel $K(s,t) = \mathcal{K}(s-t)$ defined by

$$\mathcal{K}(s-t) = \frac{1}{2\kappa\sqrt{\pi}(s-t)^{3/2}} \exp\left(-\frac{1}{4\kappa^2(s-t)}\right).$$

The parameter $\kappa$ controls the ill-conditioning of the linear system built from the discretization of the problem (59), by means of a quadrature rule (the midpoint quadrature rule described previously is used here). Throughout this numerical example, we fix this parameter to the worst ill-conditioning value, namely $\kappa = 1$. In order to validate theoretical results from Theorem 2 and Corollary 1, we use the Matlab codes *heat(n, κ)* and *get_l*$(n,2)$ implemented in *Regularization Tools* [22] for the computing of the data $\mathbf{A}$, $\mathbf{L}$, $\boldsymbol{b}$ and $\boldsymbol{u}^\star$ (the matrix $\mathbf{A}$ being a Toeplitz matrix because of the convolution type kernel $\mathcal{K}$ (see [13])). The exact solution $\boldsymbol{u}^\star = (u_j^\star)_{1\le j\le n}$, where $u_j^\star = u^\star(t_j^\star)$, is defined as follows (see [14] and [20] - Section 6.7)

$$u^\star(t) = \begin{cases} \dfrac{300}{4}t^2, & \text{if } \; 0 \le \; t \le \dfrac{1}{10}, \\ \dfrac{3}{4} + (20t-2)(3-20t), & \text{if } \; \dfrac{1}{10} < t \le \dfrac{3}{20}, \\ \dfrac{3}{4}e^{2(3-20t)}, & \text{if } \; \dfrac{3}{20} < t \le \dfrac{1}{2}, \\ 0, & \text{if } \; \dfrac{1}{2} < t \le 1. \end{cases} \tag{60}$$

The exact solution in (60) is piecewise continuous and has piecewise continuous derivatives in $[0,1]$. So, any discrete derivation operator introduced as a regularization matrix could generated more numerical errors than expected. However, in order to better capture the curvature of this solution, we nevertheless need to choose the second derivative operator provided in (55) unlike the first order operator used in [13].

#### 4.2.1. **The unperturbed data problem**

Here again, we assume that the matrix $\mathbf{A}$ and right hand term $\boldsymbol{b}$ from the discretization of the integral equation (59) are not perturbed, despite the computational errors that come from the approximation of the integral by a quadrature rule. However, we face here a very ill-conditioned (due to the fixed parameter to $\kappa = 1$) and rank-deficient matrix $\mathbf{A}$. Moreover, we use the fact that the solution we try to approximate (the exact solution in (60)) is two times piecewise differentiable with a known approximative derivative computed as follows $\boldsymbol{g} = \mathbf{L}\boldsymbol{u}^\star$.

For this validating purpose of the results in Theorem 2, we fix the dimensions to $m = n = 50$ as in [13] (any change of these dimensions will be mentioned accordingly). The built matrix has the conditioning of $\mathbf{cond}(\mathbf{A}) = 2.9474e^{+26}$ with the rank of $\mathbf{rank}(\mathbf{A}) = 48$. We thus present in Figure 12, the absolute and relative errors between the exact solution in (60) and the reconstructed stabilized-regularized solution from the solving of (22) by the Matlab command \. These errors show how the stabilized-regularized method we propose, provides accurate and stable solution regardless the choice of the regularization parameter $\gamma \in [10^{-5}, 1]$ and $\gamma \in [1, 10^5]$, and confirm theoretical results in (33) and (34).

In order to place ourselves in a more realistic context, we also investigate in the forthcoming subsection the case of perturbed data and the only *type of prior information* which is known on the *unknown* exact solution in (60) is its smoothness in the sense defined by the operator (55). We thus fix the approximation of the derivative $\boldsymbol{g}$ to zero ($\boldsymbol{g} = \mathbf{0}$).

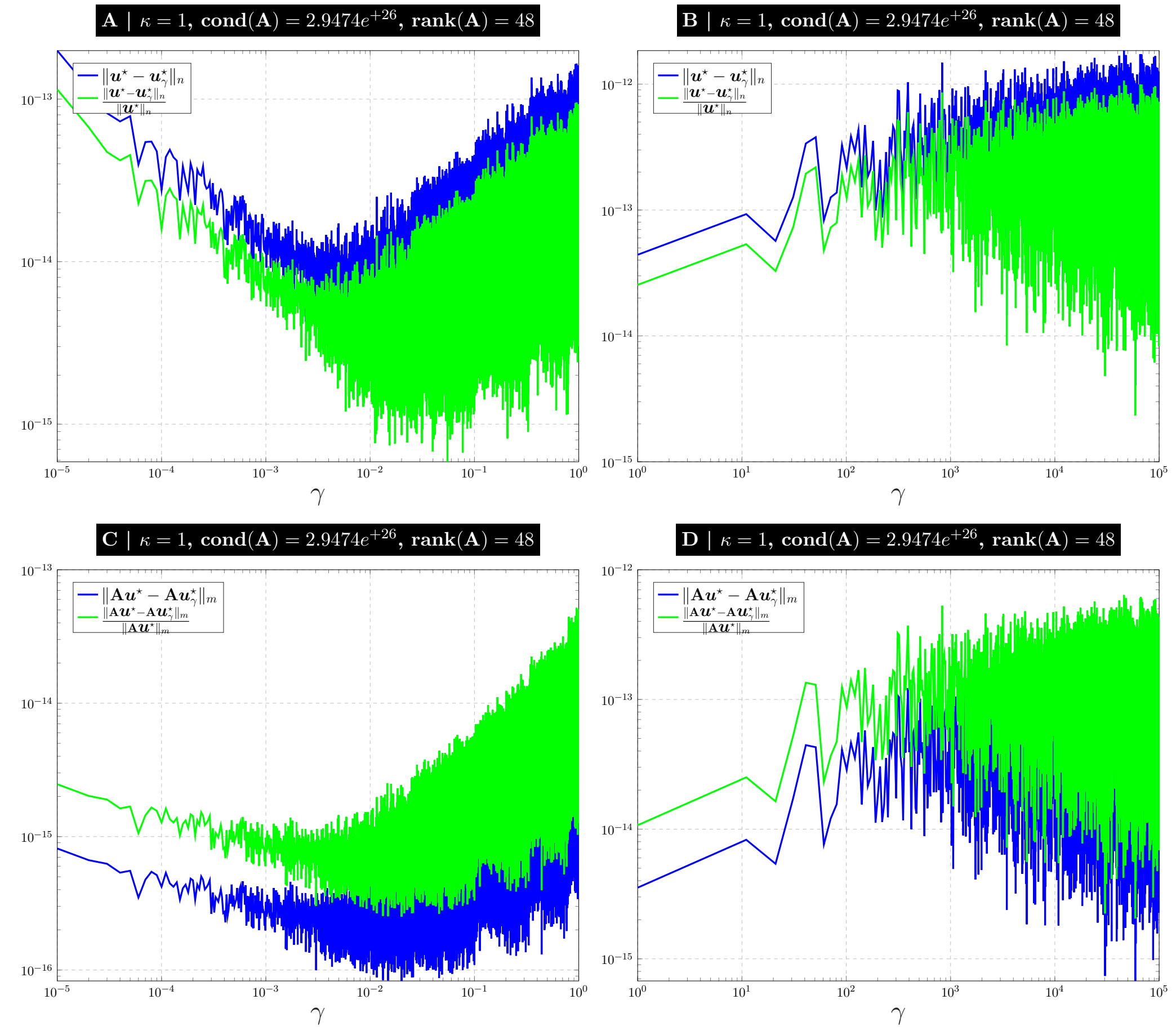

Figure 12: Absolute and relative errors between the exact solution $\boldsymbol{u}^\star$ computed from (60) and the stabilized-regularized solution $\boldsymbol{u}^\star_\gamma$ obtained from the solving of (22) by the Matlab command \, with respect to the regularization parameter $\gamma \in [10^{-5}, 1]$ (A and C) and $\gamma \in [1, 10^5]$ (B and D). The data $\mathbf{A}$, $\mathbf{L}$ and $\boldsymbol{b}$ are computed from the Matlab codes *heat(n, κ)* and *get_l(n, 2)* provided in *Regularization Tools* [22] under the dimensions $m = n = 50$, and the derivative is approximated as $\boldsymbol{g} = \mathbf{L}\boldsymbol{u}^\star$. These errors are displayed on a log-log scale

### *4.2.2.* White noise perturbed data problem

We again consider here the noise defined in (56) where only the right hand term is perturbed, and we use the matrix in (55) as the regularization matrix. As we pointed out previously, the term $\gamma\mathbf{L}^t\boldsymbol{g}$ is not taken into account in the solving of the stabilized-regularized formulation (22). We fix the dimensions $m = n = 1\,000$, and thus face both rank-deficient and ill-conditioned problem (even at lower dimensions, see [22, 24]) given by $\mathbf{rank}(\mathbf{A}) = 588$ and $\mathbf{cond}(\mathbf{A}) = 2.0332e^{+232}$, respectively. We have performed an approximation of the exact solution (60) through the stabilized-regularized formulation (22) and presented the absolute and relative errors in Figure 13 and Figure 15. The absolute errors show how the stabilized-regularized approach we propose is accurate and stable in the data (but less accurate in the solution), and confirm the result (45) in Corollary 1.

We show in Figure 16 the exact solution in (60) and the stabilized-regularized solution from (22) with respect to the index $i = 1, \cdots, n$ of the domain's mesh, for different values of the dimensions $m = n = 500, 1\,000, 5\,000$ and $10\,000$ (given in the top of the plots). We also provide in the top of these plots, the conditioning and the rank of the associated matrix. While in Figure 14, we present at the same time the corresponding exact and regularized data to the solutions represented in Figure 16. We can first notice how the conditioning of the matrix is so high just for $m = n = 500$ (and becomes $+\infty$ when $m = n = 5\,000$ and $m = n = 10\,000$), with a rank which does not evolve. These figures show how the integration error that we commit when computing the matrix $\mathbf{A}$ (in the approximation of the integral in (59)) is also harmful, when we face such ill-conditioned and rank-deficient system. Figure 16 and Figure 14 also show how the approach we propose is stable by continuing to generate an increasingly accurate approximation of the solution in (60), as long as the mesh is refined and despite a conditioning which goes to $+\infty$.

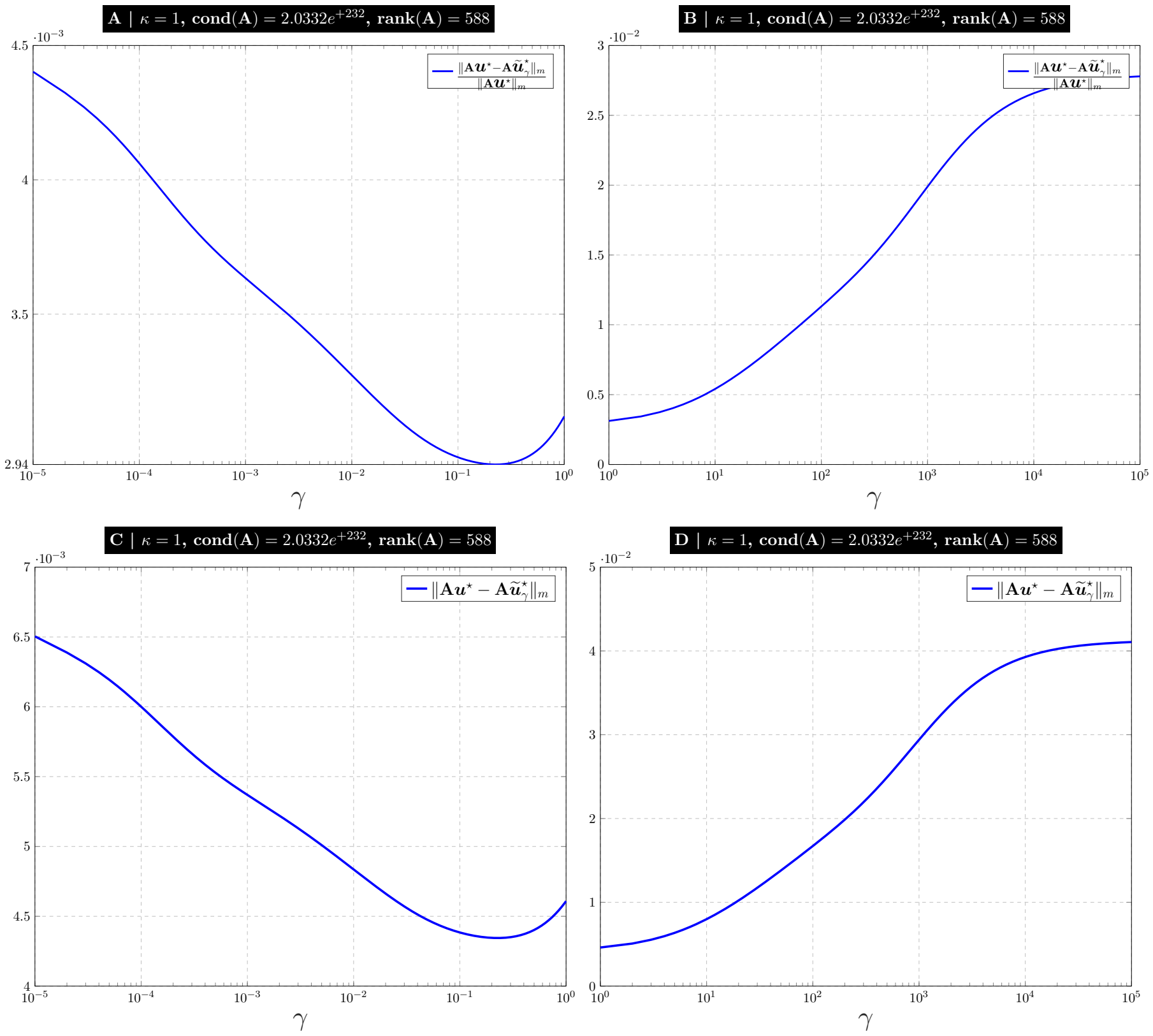

Figure 13: Absolute and relative errors between the exact data $\mathbf{A}\boldsymbol{u}^\star$ computed from (60) and the stabilized-regularized data $\mathbf{A}\widetilde{\boldsymbol{u}}^\star_\gamma$ obtained from the solving of (22) by the Matlab command \, are represented with respect to the parameter $\gamma \in [10^{-5}, 1]$ (in A and C) and $\gamma \in [1, 10^5]$ (in B and D), with a fixed $\kappa = 1$. The data $\mathbf{A}$, $\mathbf{L}$ and $\boldsymbol{b}$ are computed from the codes *heat(n, κ)* and *get_l(n, 2)* given in *Regularization Tools* [22], and the right hand term $\gamma \mathbf{L}^t \boldsymbol{g}$ in (22) is not taken account since considered as non available informations ($\boldsymbol{g} = \mathbf{0}$). The computations are made under the perturbed data in (56) and the dimensions are fixed to $m = n = 1\,000$. The horizontal axis is on a log-scale.

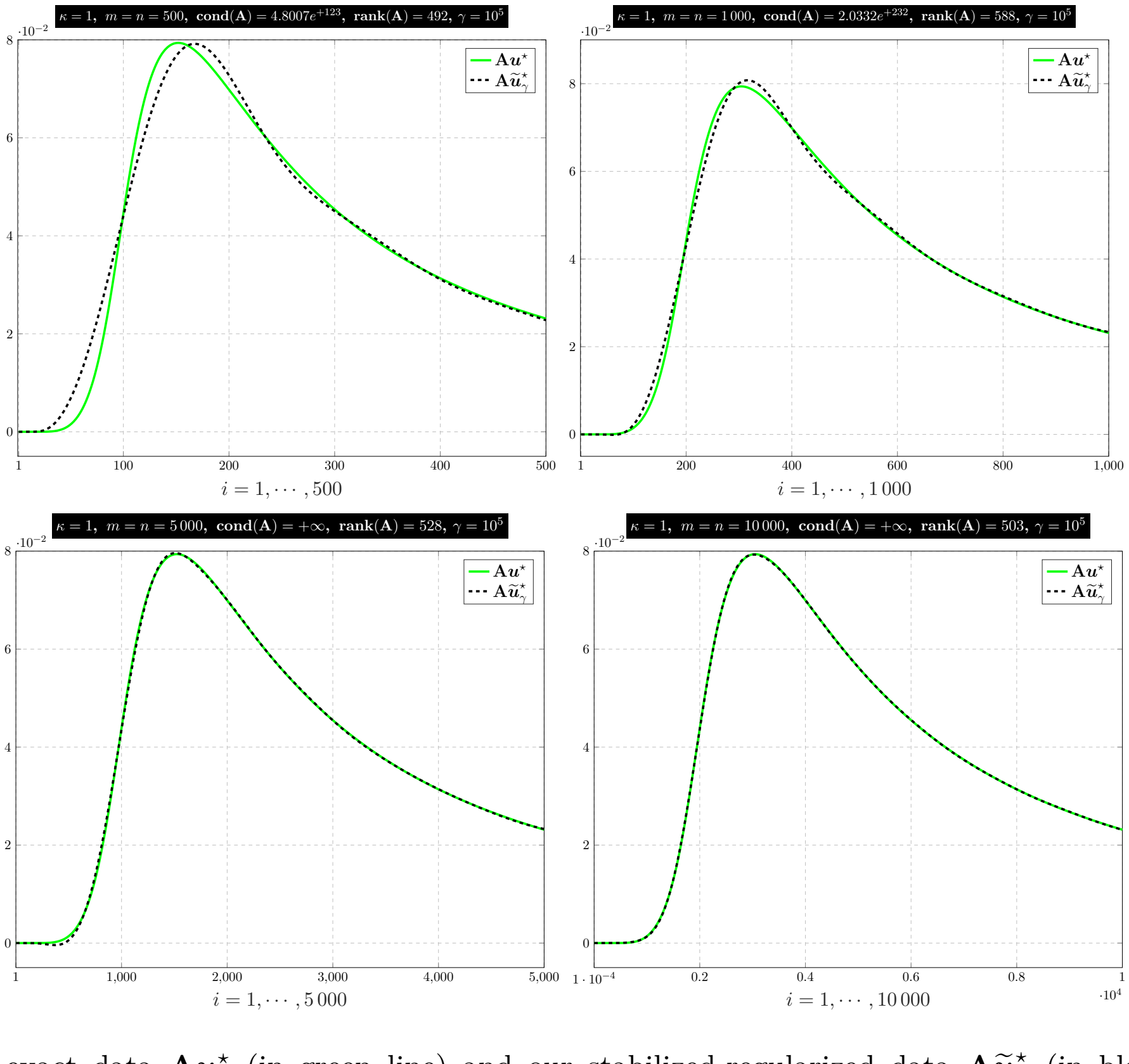

Figure 14: The exact data $\mathbf{A}\boldsymbol{u}^\star$ (in green line) and our stabilized-regularized data $\mathbf{A}\widetilde{\boldsymbol{u}}^\star_\gamma$ (in black dashed line) represented with respect to the discretization index $i = 1, \cdots, n$, for different dimensions $m = n$ (given in the top of the plots) with a fixed regularization parameter $\gamma = 10^5$ and with $\kappa = 1$. The regularized solution $\widetilde{\boldsymbol{u}}^\star_\gamma$ obtained from the solving of (22) by the Matlab command \, is computed under the perturbed data in (56) where the right hand term $\gamma \mathbf{L}^t \boldsymbol{g}$ in (22) is not taken account since being considered as non available informations ($\boldsymbol{g} = \mathbf{0}$). The data $\mathbf{A}$, $\mathbf{L}$ and $\boldsymbol{b}$ are computed from the codes *heat(n, κ)* and *get_l(n, 2)* given in *Regularization Tools* [22].

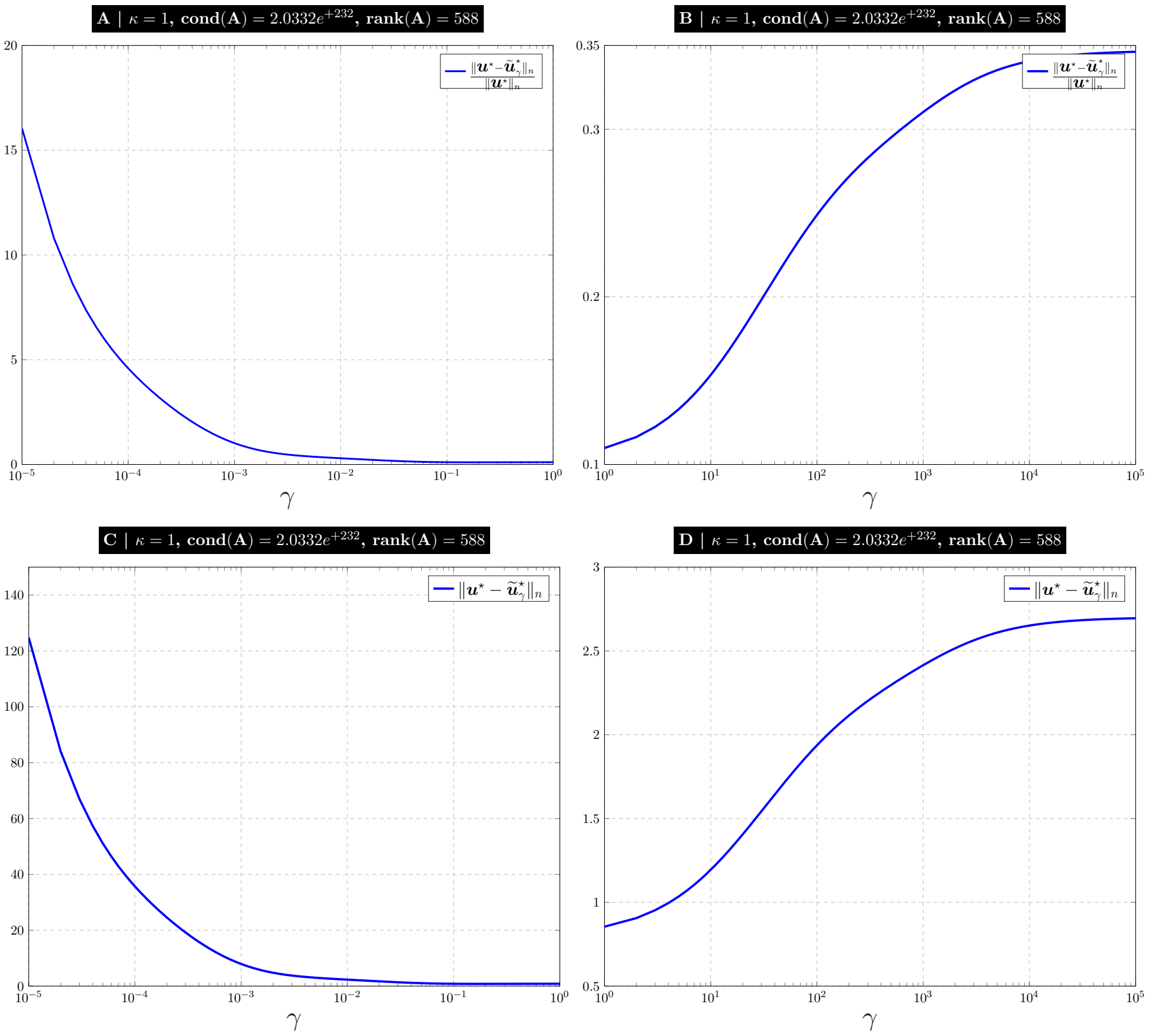

Figure 15: Absolute and relative errors between the exact solution $\boldsymbol{u}^\star$ computed from (60) and the stabilized-regularized solution $\widetilde{\boldsymbol{u}}^\star_\gamma$ obtained from the solving of (22) by the Matlab command \, represented with respect to the parameter $\gamma \in [10^{-5}, 1]$ (in A and C) and $\gamma \in [1, 10^5]$ (in B and D), with a fixed $\kappa = 1$. The data $\mathbf{A}$, $\mathbf{L}$ and $\boldsymbol{b}$ are computed from the codes *heat(n, κ)* and *get_l(n, 2)* given in *Regularization Tools* [22], where the right hand term $\gamma \mathbf{L}^t \boldsymbol{g}$ in (22) is not taken account since considered as non available informations ($\boldsymbol{g} = \mathbf{0}$). The computations are made under the perturbed data in (56) and the dimensions are fixed to $m = n = 1\,000$. The horizontal axis is on a log-scale.

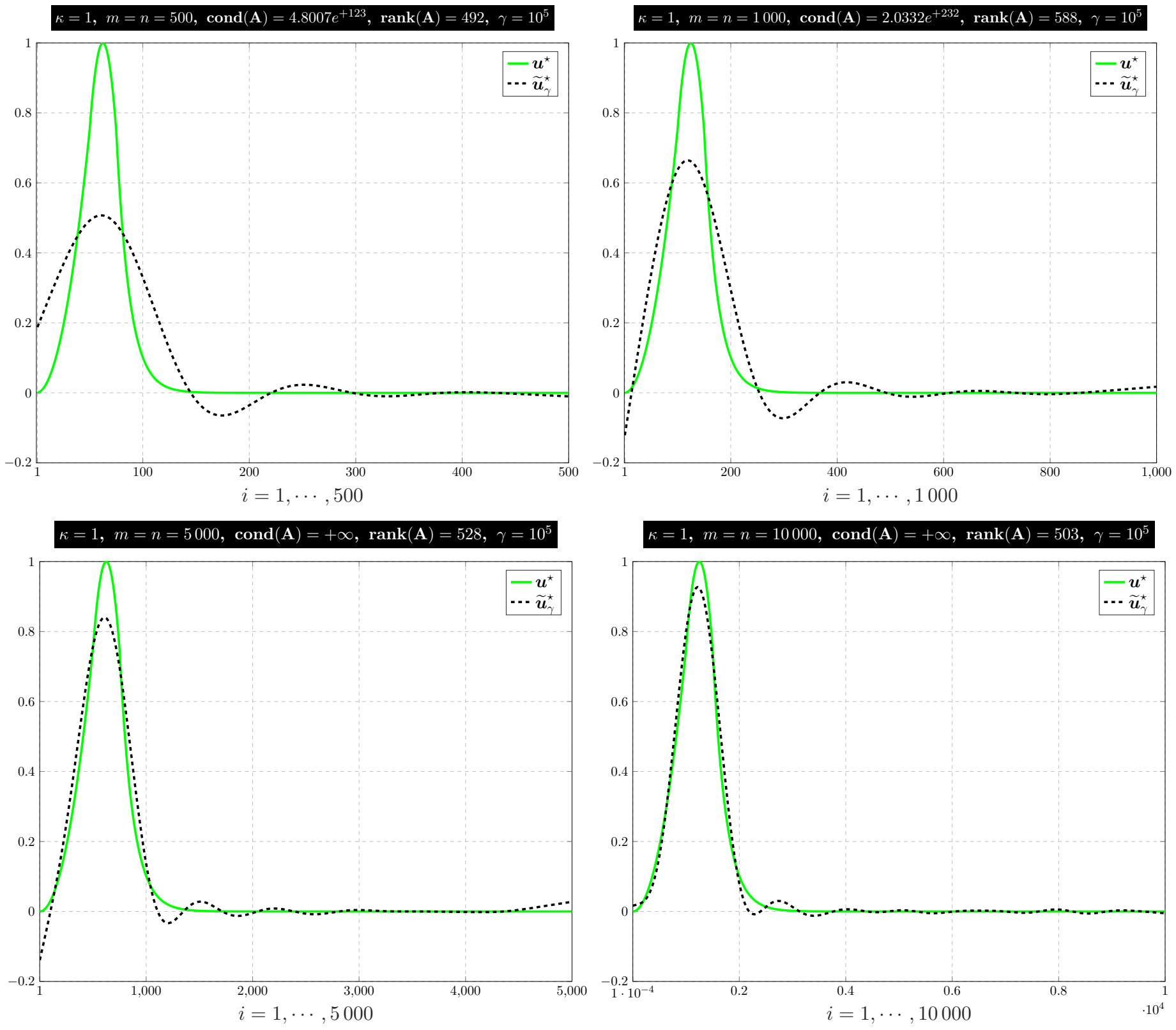

Figure 16: The exact solution $\boldsymbol{u}^\star$ (in green line) computed from (60) and our stabilized-regularized solution $\widetilde{\boldsymbol{u}}^\star_\gamma$ (in black line) obtained from the solving of (22) by the Matlab command \, represented with respect to the discretization index $i = 1, \cdots, n$, for different dimensions $m = n$ (given in the top of the plots) with a fixed regularization parameter $\gamma = 10^5$ and with $\kappa = 1$. The computations are made under the perturbed data in (56) and the right hand term $\gamma \mathbf{L}^t \boldsymbol{g}$ in (22) is not taken account since being considered as non available informations ($\boldsymbol{g} = \mathbf{0}$).

*4.3.* **Example 3***: Phillips' famous test problem*

We illustrate here the Phillips problem also provided in Regularization Tools [22], which is one of the classical test problems from the literature [33]. The Galerkin-type method with orthonormal basis functions (see [2] and [12] - Chap. 7) are used to discretize the Fredholm integral equation

$$\int_a^b K(s,x)u^\star(x)dx = b(s), \;\; c \le s \le d, \tag{61}$$

where the solution $u^\star$, the kernel $K$ and the right-hand side $b$ are defined by

$$u^\star(x) = \begin{cases} 1+\cos\left(\frac{\pi}{3}x\right), & \text{if } \; |x|<3 \\ 0, & \text{if } \; |x|\ge 3 \end{cases} \tag{62}$$

$$K(s,x) = u^\star(s-x), \tag{63}$$

$$b(s) = (6-|s|)\left(1+\frac{1}{2}\cos\left(\frac{\pi}{3}s\right)\right) + \left(\frac{9}{2\pi}\right)\sin\left(\frac{\pi}{3}|s|\right), \tag{64}$$

and where $a \le x \le b$ and $c \le s \le d$. For the building of the linear system, the components $a_{ij}$ and $b_i$, $1 \le i,j \le n$, of the matrix $\mathbf{A}$ and the right-hand term $\boldsymbol{b}$ are given by

$$a_{ij} = \int_a^b \int_c^d K(s,x)\phi_i(s)\psi_j(x)dsdx, \quad b_i = \int_c^d \phi_i(s)g(s)ds \tag{65}$$

where $\phi_i$ and $\psi_j$ are basis functions used in the Galerkin method and defined as follows [22]

$$\phi_i(s) = \begin{cases} h_s^{-\frac{1}{2}}, & s\in[s_{i-1},s_i] \\ 0, & \text{elsewhere} \end{cases}, \quad \psi_i(x) = \begin{cases} h_x^{-\frac{1}{2}}, & x\in[x_{i-1},x_i] \\ 0, & \text{elsewhere} \end{cases}$$

in which $h_s = (d-c)/n,\; h_x = (b-a)/n$, and $x_i = a + ih_x,\; s_i = c + ih_s,\; i = 1,\dots,n$, where the domain is defined by $a = c = -6$ and $b = d = 6$. The reference solution that we use here is denoted by the vector $\boldsymbol{u}^\star = (u_j^\star)_{1\le j\le n}$, where its components are obtained from

$$u_j^\star = \int_a^b \psi_j(x)u^\star(x)dx, \quad j = 1,\dots,n. \tag{66}$$

The analytical solution we face in (62) has a continuous first derivative in $[-6, 6]$, but presents a piecewise second derivative. In order to avoid possible huge numerical derivation errors by considering the second derivative operator in (55), we then use the following first derivative operator as a regularization matrix

$$\mathbf{L} = \begin{bmatrix} -1 & 1 & 0 & \cdots & 0 \\ 0 & -1 & 1 & \ddots & \vdots \\ \vdots & \ddots & \ddots & \ddots & 0 \\ 0 & \cdots & 0 & -1 & 1 \end{bmatrix} \in \mathbb{R}^{(n-1)\times n} \tag{67}$$

We fixe here the number of integration points to $m = n = 1\,000$ in order to minimize errors from the quadrature in the computing of the matrix and the right hand term in (65), but also the reference solution in (66). We use the code *phillips(n)* provided in Regularization Tools [22] for the computing of the data $\mathbf{A}$, $\mathbf{L}$, $\boldsymbol{b}$ and the solution components $u_j^\star$ in (66).

*4.3.1.* **White noise perturbed data problem**

The white noise perturbation provided in (56) is first considered here and the regularization matrix in (67) is used in the solving of the stabilized-regularized formulation (22), where the term $\gamma\mathbf{L}^t\boldsymbol{g}$ will not be taken into account. The obtained matrix has a maximal rank of $\mathbf{rank}(\mathbf{A}) = 1\,000$, but with a huge conditioning of $\mathbf{cond}(\mathbf{A}) = 2.6415e^{+10}$.

We computed relative and absolute errors on the solution and the data presented in Figure 17 and Figure 18. These results are represented with respect to the regularization parameter $\gamma \in [10^{-5}, 1]$ and $\gamma \in [1, 10^5]$, and are consistent compared to the theoretical results in (45) for any $\gamma > 0$, and in (46) only for $\gamma > 1$. We also present in Figure 20, the associated exact solution and stabilized-regularized solution as functions of the index $i = 1, \cdots, n$, for different values of the parameter $\gamma = 10^{-5}, 1, 10^5, 10^{10}$. This example also proves how our approach is clearly a remedy to ill-conditioning linear systems, and does not need the finding of a suitable value for the regularization parameter. We only need a regularization parameter as large as we want (which means an over-regularization). Sometimes, we also need to refine the mesh in order to obtain a more accurate approximation (as we saw in the previous numerical example).

Increasing the regularization parameter $\gamma$ produces greater smoothing and is what we see in Figure 20 by taking $\gamma = 1, 10^5, 10^{10}$, showing that our approach provides accurate and stable estimations of the exact solution. This stability is also seen in the Figure 19 - B where we observe very small variations of the semi-norm $\|\mathbf{L}\widetilde{\boldsymbol{u}}_\gamma^\star\|_{n-1}$ of the stabilized-regularized solution $\widetilde{\boldsymbol{u}}_\gamma^\star$ with respect to very small variations of the corresponding residual $\|\mathbf{A}\widetilde{\boldsymbol{u}}_\gamma^\star - \tilde{\boldsymbol{b}}\|_m$, despite the wide range of values of the regularization parameter $\gamma \in \left[1, 10^5\right]$. However for $\gamma \in \left[10^{-5}, 1\right]$ in A, we observe a less stable solution.

The results obtained in Figure 18 and Figure 20 are similar to those provided in [9] (see Example 4.3 and Example 4.5) with the same noise level but under lower dimensions. Under higher dimensions and noise level, results in Figure 20 are not only accurate but also stable with respect to the regularization parameter $\gamma$.

#### *4.3.2.* **Perturbations by correlated noise - filtered white noises**

Under the perturbation in (57) and with the regularization matrix in (67), we compute the approximated solution $\widetilde{\boldsymbol{u}}_\gamma^\star$ from solving of the stabilized-regularized formulation (22) by the Matlab command \, without taking account the term $\gamma \mathbf{L}^t \boldsymbol{g}$. We also face here an ill-conditioned system with the conditioning of $\mathbf{cond}(\mathbf{A}) = 2.8990e^{+10}$, even if the matrix is of maximal rank ($\mathbf{rank}(\mathbf{A}) =$ 1 000).

We present in the Annexe 6.2.1 at Figure 27 and Figure 28 absolute and relative errors on the solution and its associated data, with respect to the regularization parameter $\gamma$. Whereas in Figure 29, we represent the corresponding exact solution and stabilized-regularized solution as functions of the index $i = 1, \cdots, n$, for different values of the parameter $\gamma = 10^{-5}, 1, 10^5, 10^{10}$. Despite we face an ill-conditioned matrix but with a maximal rank, we obtain in Figure 27 (for $\gamma > 0$) and in Figure 28 (for $\gamma$ big enough, namely $\gamma > 10$) estimates illustrating theoretical results (46) and (45) in Corollary 1.

#### *4.3.3.* **Perturbations by correlated noise - perturbations from smoothing of the data**

Finally under the correlated noise given in (58) and by the use of the regularization matrix in (67), we also determine an approximation $\widetilde{\boldsymbol{u}}_\gamma^\star$ from the stabilized-regularized formulation (22) by the Matlab command \, without taking account the term $\gamma \mathbf{L}^t \boldsymbol{g}$. The computed matrix from this perturbation is of maximal rank and with the conditioning of $\mathbf{cond}(\mathbf{A}) = 2.7583e^{+11}$.

We present in the Annexe 6.2.2 at Figure 30 and Figure 31 absolute and relative errors on the data and on their associated solutions with respect to the regularization parameter $\gamma$, whereas in Figure 32 are displayed the corresponding solutions. We obtain accurate and stable results which are coherent to the theoretical results (46) and (45) in Corollary 1.

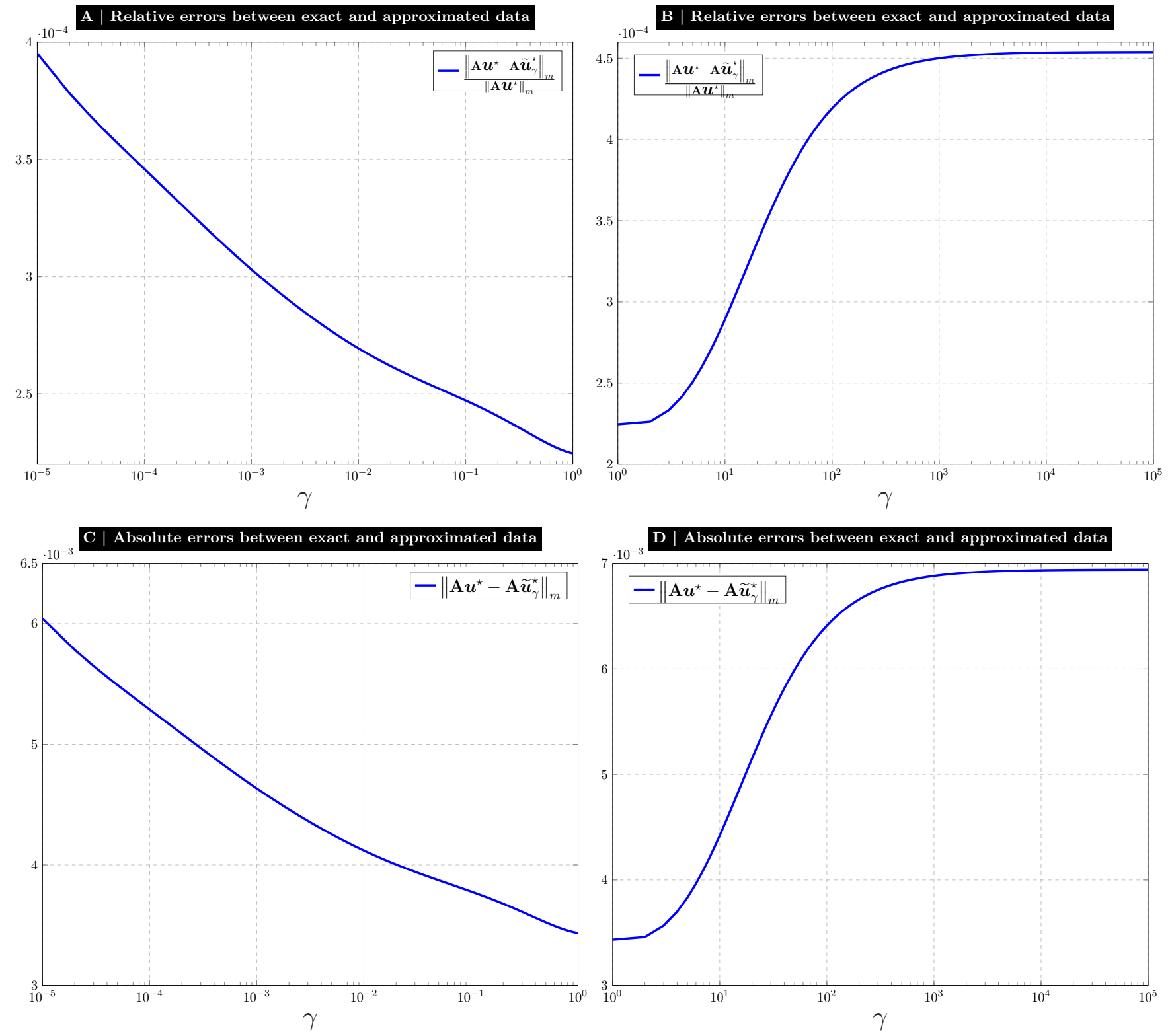

Figure 17: Absolute and relative errors between the exacte data $\mathbf{A}\boldsymbol{u}^\star$ obtained from (66) and our stabilized-regularized data $\mathbf{A}\widetilde{\boldsymbol{u}}_\gamma^\star$ obtained from the solving of (22) by the Matlab command \, with respect to the parameter $\gamma \in [10^{-5}, 1]$ (in A and C) and $\gamma \in [1, 10^5]$ (in B and D). The data $\mathbf{A}$, $\mathbf{L}$ and $\boldsymbol{b}$ are obtained from the Matlab codes *phillips(n)* and *get_l*$(n, 1)$ provided in *Regularization Tools* [22], the computations are made under the perturbed data in (56), the right hand term $\gamma \mathbf{L}^t \boldsymbol{g}$ in (22) is not taken account since being considered as non available informations,($\boldsymbol{g} = \mathbf{0}$), and the dimensions are fixed to $m = n = 1\,000$ (where $\mathbf{cond}(\mathbf{A}) = 2.6415e^{+10}$ and $\mathbf{rank}(\mathbf{A}) = 1\,000$).

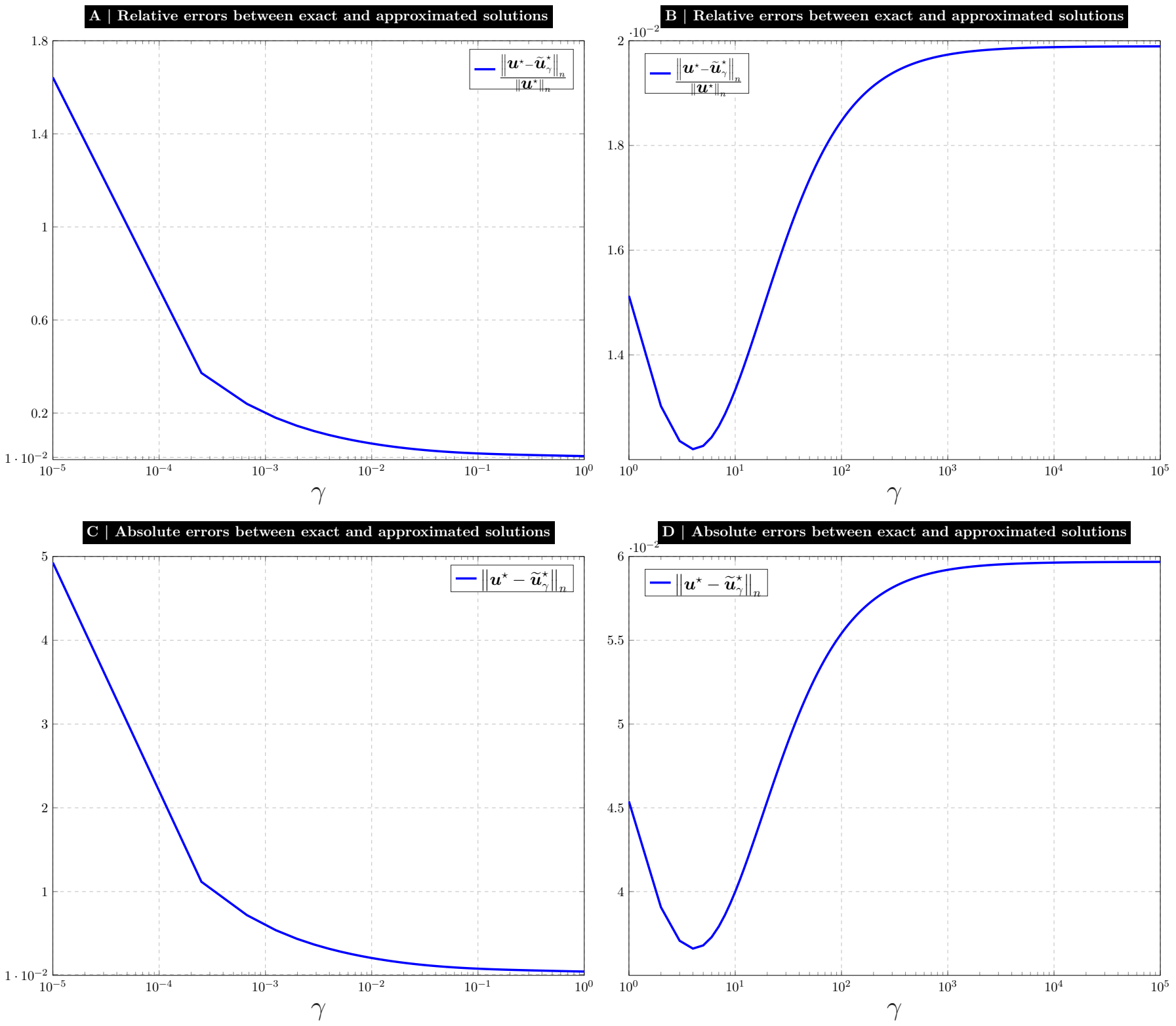

Figure 18: Relative and absolute errors between the exacte solution $\boldsymbol{u}^\star$ computed from (66) and our stabilized-regularized solution $\widetilde{\boldsymbol{u}}_\gamma^\star$ obtained from the solving of (22) by the Matlab command \, with respect to the parameter $\gamma \in [10^{-5}, 1]$ (in A and C) and $\gamma \in [1, 10^5]$ (in B and D). The data $\mathbf{A}$, $\mathbf{L}$ and $\boldsymbol{b}$ are obtained from the Matlab codes *phillips(n)* and *get_l*$(n, 1)$ provided in *Regularization Tools* [22], the computations are made under perturbed data in (56), the right hand term $\gamma \mathbf{L}^t \boldsymbol{g}$ in (22) is not taken account since being considered as non available informations ($\boldsymbol{g} = \mathbf{0}$), and the dimensions are fixed to $m = n = 1\,000$ (where $\mathbf{cond}(\mathbf{A}) = 2.6415e^{+10}$ and $\mathbf{rank}(\mathbf{A}) = 1\,000$).

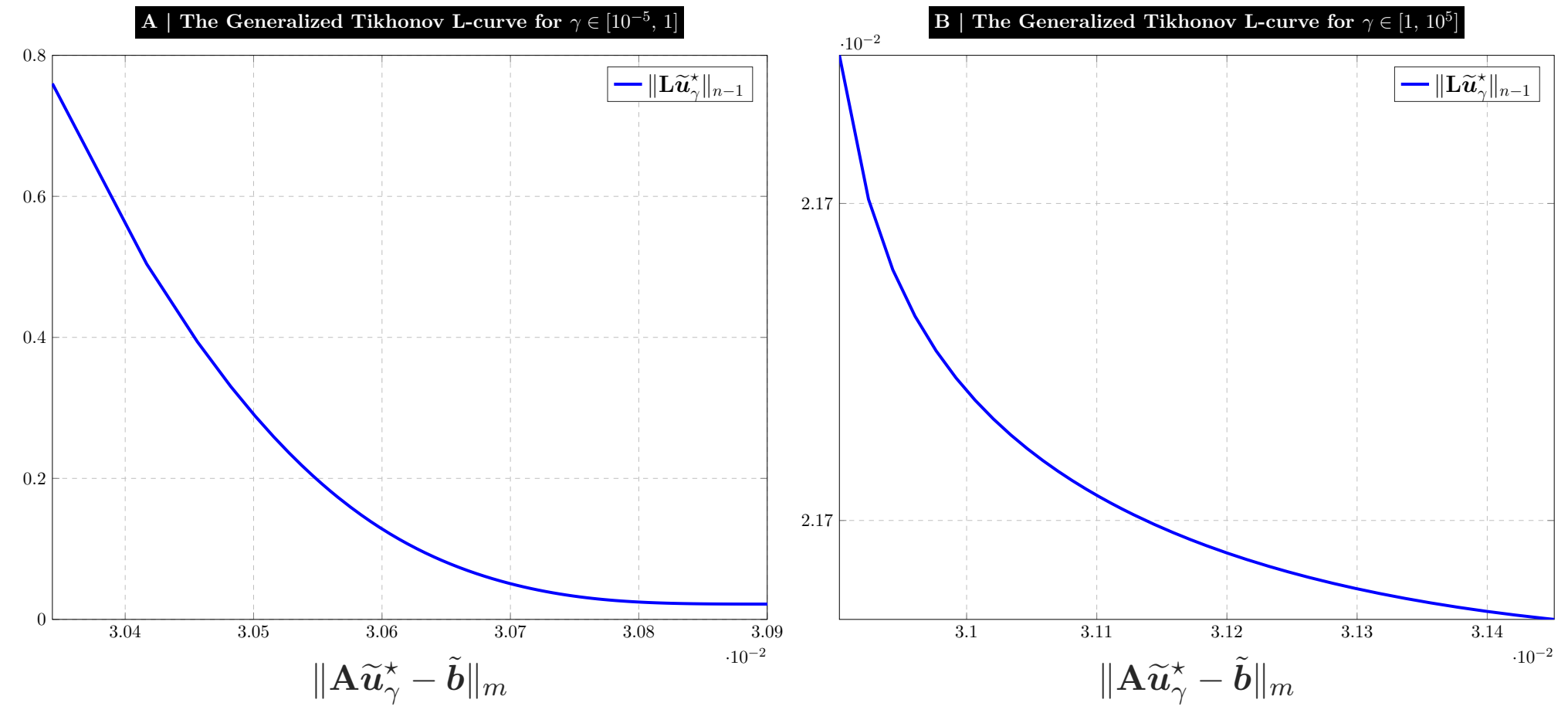

Figure 19: The residual norm $\|\mathbf{A}\widetilde{\boldsymbol{u}}^\star_\gamma - \tilde{\boldsymbol{b}}\|_m$ versus the semi-norm $\|\mathbf{L}\widetilde{\boldsymbol{u}}^\star_\gamma\|_{n-1}$ represented for the range of values of the regularization parameter $\gamma \in [10^{-5}, 1]$ (in A) and $\gamma \in [1, 10^5]$ (in B). The stabilized-regularized solution $\boldsymbol{u}^\star_\gamma$ is obtained from the solving of (22) by the Matlab command \ without taken into account the right hand term $\gamma \mathbf{L}^t \boldsymbol{g}$ in (22) (since being considered as non available informations). The data $\mathbf{A}$, $\mathbf{L}$ and $\boldsymbol{b}$ are obtained from the Matlab codes *phillips(n)* and *get_l*$(n, 1)$ provided in *Regularization Tools* [22] and the computations are made under the perturbed data in (56), with respect to the fixed dimensions $m = n = 1\,000$ and where $\mathbf{cond}(\mathbf{A}) = 2.6415e^{+10}$ and $\mathbf{rank}(\mathbf{A}) = 1\,000$.

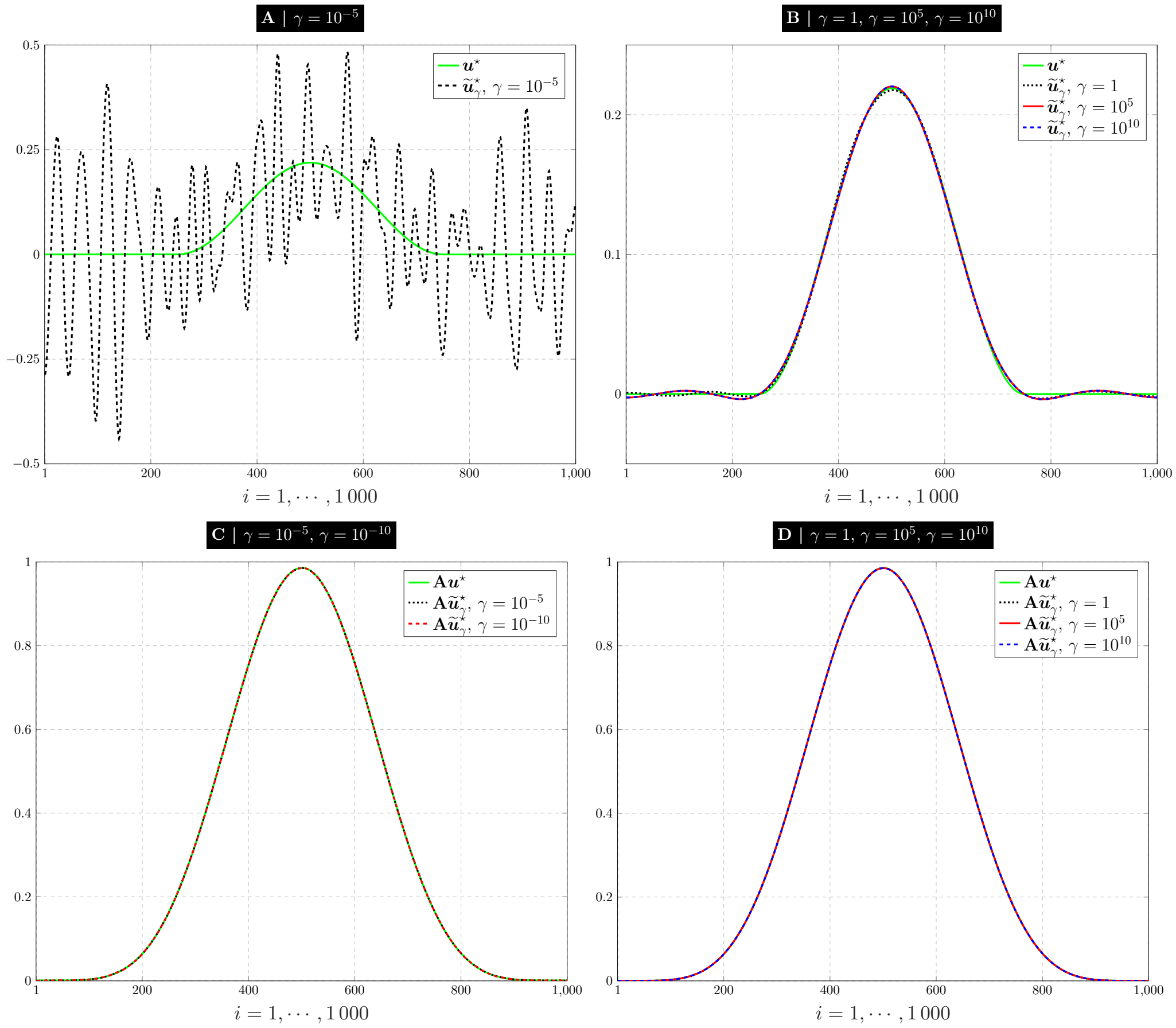

Figure 20: The exact solution $\boldsymbol{u}^\star$ (in green solid line) computed from (66) and our regularized solution $\widetilde{\boldsymbol{u}}^\star_\gamma$ (in black dashed line) obtained from the solving of (22) by the Matlab command \, presented with respect to the discretization index $i = 1, \cdots, 1\,000$, for different values of the regularization parameter $\gamma = 10^{-5}, 1, 10^5, 10^{10}$. The data $\mathbf{A}$, $\mathbf{L}$ and $\boldsymbol{b}$ are obtained from the Matlab codes *phillips(n)* and *get_l*$(n, 1)$ provided in *Regularization Tools* [22] and the computations are made under the perturbed data in (56) with respect to the fixed dimensions $m = n = 1\,000$, where $\mathbf{cond}(\mathbf{A}) = 2.6415e^{+10}$ and $\mathbf{rank}(\mathbf{A}) = 1000$.

## 5. General conclusion

Through this work, we have analyzed the stabilization of the Tikhonov regularization method for the solving of the least squares problem with or without perturbed data and under an ill-conditioned matrix. The stabilized-regularized approach developed here is based on the idea that, to solve the *least squares problem* (1) by the regularization method, we need to control at the same time the residual (through the *normal equation* (2)) and the solution (or one of its derivative, through some type of informations $\mathbf{L}\boldsymbol{u} = \mathbf{g}$).

This approach has spared us an intrinsic difficulty of the regularization of discrete ill-conditioned problems: to resort to a method for determining the regularization parameter to obtain an accurate solution. The stabilized-regularized method proposed here provides an accurate solution when the regularization parameter is large (over-regularizing), and somtimes with the need to refine the domain's mesh (see Figures 9, 11, 14 and 16).

Under unperturbed data and a maximal rank matrix $\mathbf{A}$, the stabilized-regularized method provides the minimal norm solution of the least squares problem (Theorem 2). This approach doesn't thus present regularization errors, despite the ill-conditioning of the matrix we may face. Figures 2 and 12 confirm these theoretical results.

Otherwise, under perturbed data and a maximal rank matrix, the stabilized-regularized method provides numerical results that are coherent to the estimate established at Corollary 1, which is what we hoped for since regularization errors are zero. Moreover, insights from Singular Value Decomposition show bias from the stabilized-regularized method smaller than the one from the standard Tikhonov method.

Even if the analysis we carried out through this work (Theorem 2 and Corollary 1) is independent of the used regularization norm on the functional in (19), we applied the Euclidean norm (Ridge regularization) to obtain the problem (21) from which we have made the numerical tests. A very important point to investigate is to analyse this stabilized-regularized approach, when Lasso regularization is used in place of Ridge regularization.

It will also be interesting to deeply investigate the Singular Value Decomposition and the spectral analysis of the stabilized-regularized method proposed here. This could be among other future avenues to explore.

## Statements and Declarations:

- I acknowledge the support of the Natural Sciences and Engineering Research Council of Canada (NSERC, https://www.nserc-crsng.gc.ca/).
- I confirm that this work is original and has not been published elsewhere, nor is it currently under consideration for publication elsewhere.
- I have no conflicts of interest to disclose.

## References

[1] Allaire G. and Kaber S. M., *Numerical Linear Algebra*, Springer (2008).

[2] Anderssen R.S. and Prenter P.M., *A formal comparison of methods proposed for the numerical solution of first kind integral equations*, J. Austral. Math. Soc. (Series B) 22, 488-500 (1981).

[3] Babuska, I., *The Finite Element Method with Penalty.* Mathematics of Computation, 27 (122), pp. 221-228 (1973).

[4] Baranger J., *Analyse numérique*, Collection enseignement des scineces, 38. Hermann (1991).

[5] Björck Å., *Numerical Methods in Matrix Computations*, Springer (2014).

[6] Boyd S. and Vandenberghe L., *Introduction to Applied Linear Algebra - Vectors, Matrices, and Least Squares*, Cambridge University Press (2018).

[7] Brezinski C., *Algorithme numérique*, Ellipses (1988).

[8] Brezinski C., Redivo-Zaglia M., Rodriguez G., Seatzu S., *Multi-parameter regularization techniques for ill-conditioned linear systems*, Numerische Mathematik, 94, pp. 203-228 (2003).

[9] Calvetti D. and Reichel L., *Tikhonov Regularization of Large Linear Problems*, BIT Numerical Mathematics, 43, pp. 263-283, (2003).

[10] Calvetti D. and Reichel L., *Tikhonov Regularization with a Solution Constraint*, SIAM Journal on Scientific Computing, 26 (1), pp. 224-239 (2004).

[11] Ciarlet P. G., *Introduction à l'analyse numérique matricielle et optimisation*, Masson (1982).

[12] Delves L. M. and Delves J., *Numerical Solution of Integral Equations*, Clarendon Press, Oxford, (1974).

[13] Eldén L., *The numerical solution of a non-characteristic Cauchy problem for a parabolic equation*, in P. Deuflhart and E. Hairer (Eds.), *Numerical Treatment of Inverse Problems in Differential and Integral Equations*, BirkhAauser, Boston (1983).

[14] Eldén L., *Numerical solution of the sideways heat equation*, Inverse Problem 11, 913-923 (1995).

[15] Engl H. W., Hanke M. and Neubauer A., *Regularization of Inverse Problems*, Springer (1996).

[16] Fierro R. D. and Bunch J. R., *Collinearity and Total Least Squares*, SIAM Journal on Matrix Analysis and Applications, 4, 15 (1994).

[17] Gallier J. and Quaintance J. *Linear Algebra and Optimization with Applications to Machine Learning. Volume II: Fundamentals of Optimization Theory with Applications to Machine Learning*, World Scientific (2020).

[18] Golub G. H., Heath M. and Wahba G., *Generalized Cross-Validation as a Method for Choosing a Good Ridge Parameter*, Technometrics, 21 (1979).

[19] Hackbusch W., *Integral Equations: Theory and Numerical Treatment*, Birkhäuser (2011).

[20] Hanke M., *Conjugate gradient type methods for ill-posed problems*, Longman Scientific and Technical, Wiley (1995).

[21] Hansen P. C., *Discrete inverse problems: insight and algorithms*, Society for Industrial and Applied Mathematics (2010).

[22] Hansen P. C., *REGULARIZATION TOOLS: A Matlab package for analysis and solution of discrete ill-posed problems*, Numerical Algorithms, 46, 189-194 (2007).

[23] Hansen P. C., *Analysis of Discrete Ill-Posed Problems by Means of the L-Curve*, SIAM Review, 4, 34, pp. 561-580 (1992).

[24] Hansen P. C., *Rank-deficient and discrete ill-posed problems: numerical aspects of linear inversion*, SIAM monographs on mathematical modeling and computation (1998).

[25] Hansen P. C. and O'Leary D. P., *The Use of the L-Curve in the Regularization of Discrete Ill-Posed Problems*, SIAM Journal on Scientific Computing, 6, 14 (1993).

[26] Hansen P. C., Nagy J. G., and O'Leary D. P., *Deburring Images: Matrices, Spectra, and Filtering*, SIAM, Philadelphia, (2006).

[27] Kern M., *Numerical Methods for Inverse Problems*, John Wiley & Sons (2016).

[28] Kilmer M. E. and O'Leary D. P., *Choosing Regularization Parameters in Iterative Methods for Ill-Posed Problems*, SIAM Journal on Matrix Analysis and Applications, 4, 22 (2001).

[29] Levin E. and Meltzer A. Y., *Estimation of the Regularization Parameter in Linear Discrete Ill-Posed Problems Using the Picard Parameter*, SIAM Journal on Scientific Computing, 6, 39 (2017).

[30] Morozov V. A., *Methods for Solving Incorrectly Posed Problems*, Springer-Verlag (1984).

[31] Morozov V. A., *Regularization Methods For Ill-Posed Problems*, CRC Press (1993).

[32] O'leary D. P., *Near-Optimal Parameters for Tikhonov and Other Regularization Methods*, SIAM Journal on Scientific Computing, 4, 23 (2001).

[33] Phillips D. L., *A Technique for the Numerical Solution of Certain Integral Equations of the First Kind*, Journal of the ACM, 1, 9, pp. 84-97 (1962).

[34] Quarteroni A., Sacco R. and Saleri F., *Numerical Mathematics*, Second Edition, Springer Berlin Heidelberg (2007).

[35] Quarteroni A., Saleri F. and Gervasio P., *Calcul Scientifique. Cours, exercices corrigés et illustrations en MATLAB et Octave*, Deuxième édition, Springer-Verlag Italia (2010).

[36] Shaw C. B., *Improvement ofthe resolution ofan instrument by numerical solution ofan integral equation*, Journal of Mathematical Analysis and Applications, 37, pp. 83-112 (1972).

[37] Sibony M., Mardon J. C., *Analyse numérique I. Systèmes linéaires et non linéaires*, Hermann (1988).

[38] Smola A. J. and Scholkopf B., *A tutorial on support vector regression*, Statistics and Computing, 14, pp. 199-222 (2004).

[39] Vapnik V. N., *The Nature of Statistical Learning Theory*, Springer-Verlag, New York (1995).

[40] Wahba G., *Spline Models for Observational Data*, SIAM (1990).

# 6. ANNEXE

## *6.1.* Example 1: One-dimensional image restoration model

### *6.1.1.* Filtered and corrolated white noises in (57)

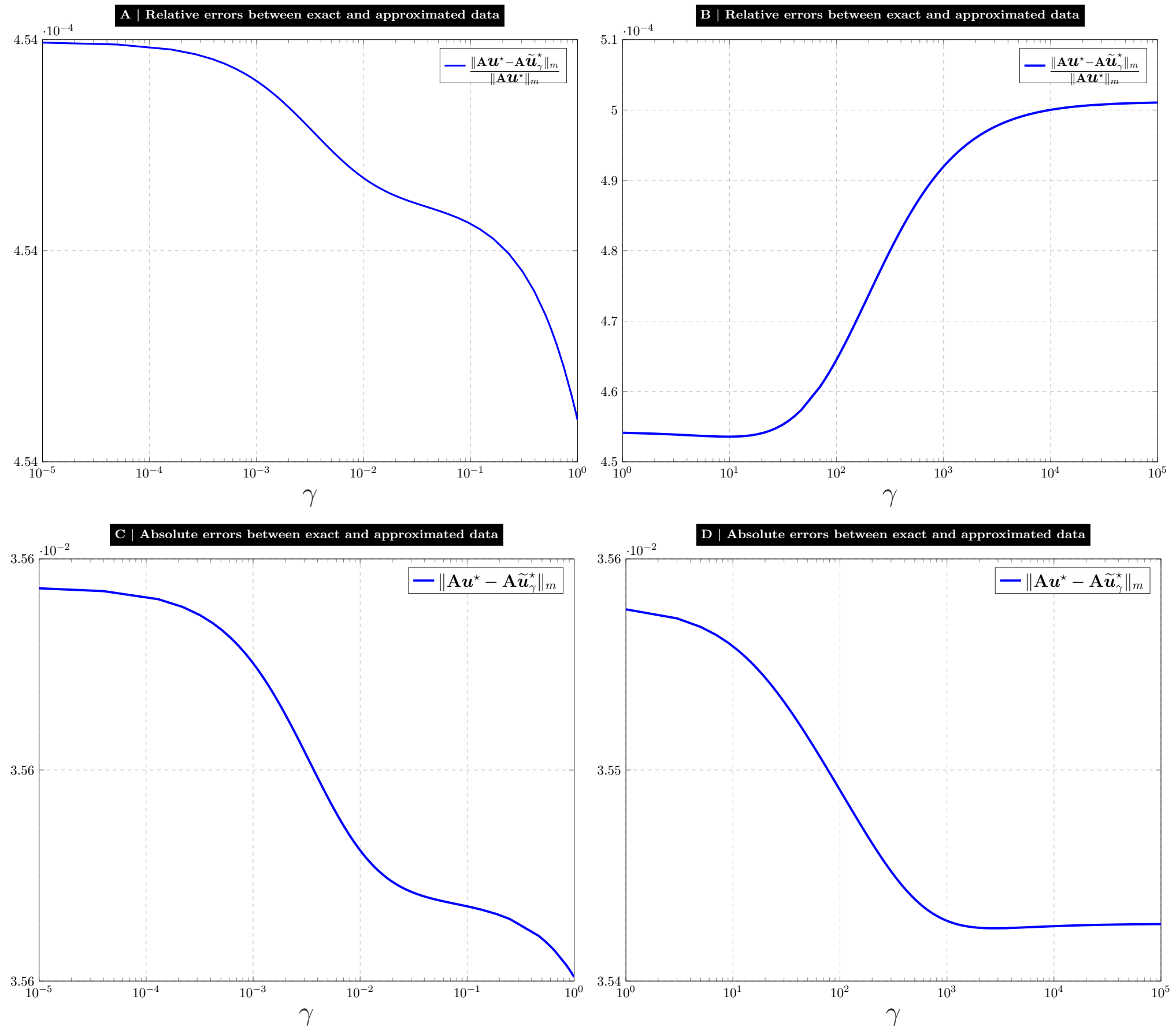

Figure 21: Absolute and relative errors between the exact data $\mathbf{A}\boldsymbol{u}^{\star}$ computed from (48), and the stabilized-regularized data $\mathbf{A}\widetilde{\boldsymbol{u}}^{\star}_{\gamma}$ obtained from the solving of (22) by the Matlab command \, with respect to the parameter $\gamma \in [10^{-5}, 1]$ (in A and C) and $\gamma \in [1, 10^{5}]$ (in B and D). The data $\mathbf{A}$, $\mathbf{L}$ and $\boldsymbol{b}$ are computed from the Matlab codes *shaw(n)* and *get_l*$(n, 2)$ provided in *Regularization Tools* [22] under the dimensions $m = n = 1\,000$ (with $\mathbf{cond}(\mathbf{A}) = 5.4935e^{+20}$ and $\mathbf{rank}(\mathbf{A}) = 20$), and the right hand term $\gamma \mathbf{L}^t \boldsymbol{g}$ in (22) is not taken into account since being considered as non available informations ($\boldsymbol{g} = \mathbf{0}$). The computations are made under parameters in (49) with the perturbed data (57), and the horizontal axis is displayed on a log-scale.

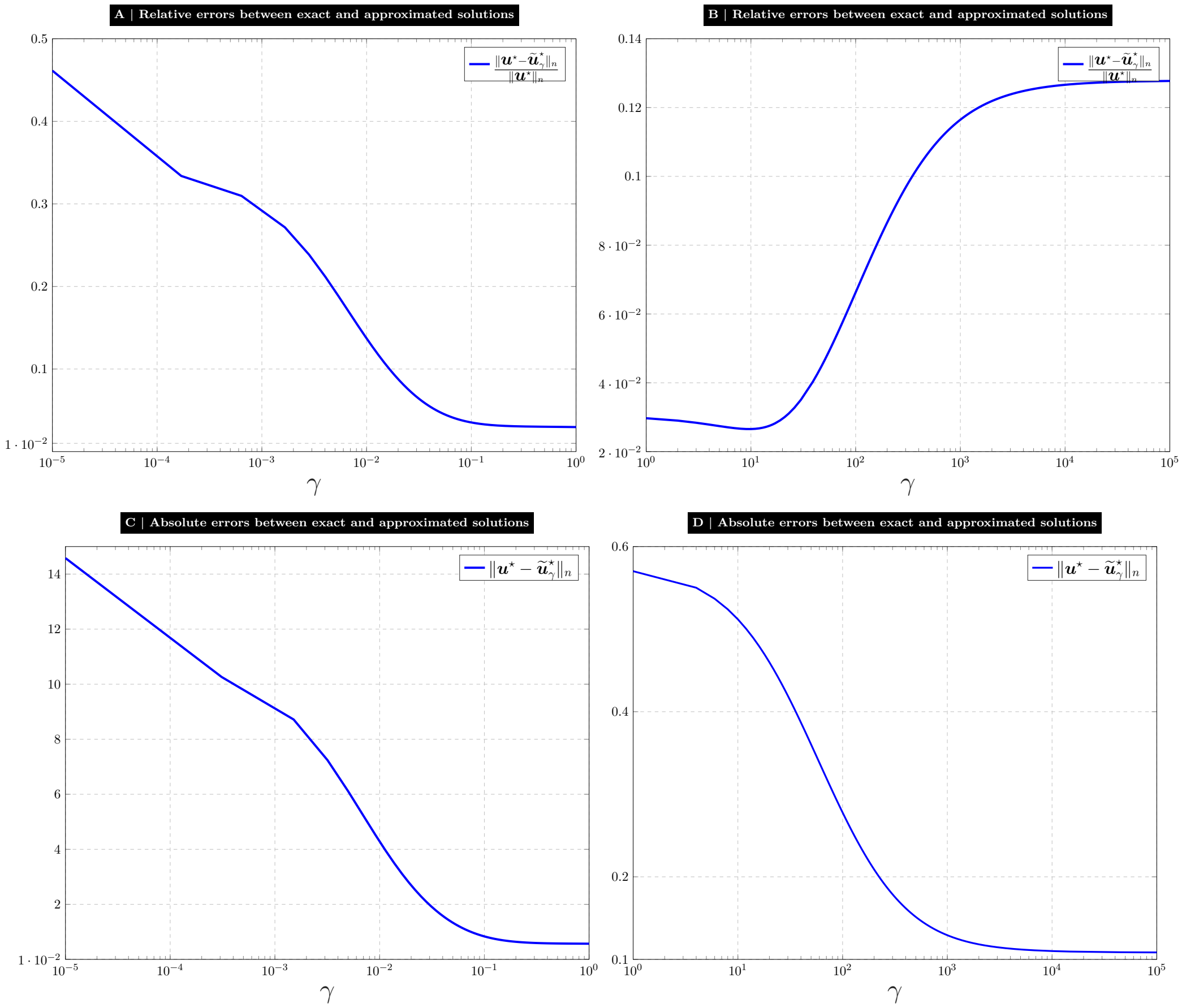

Figure 22: Absolute and relative errors between the exact solution $\boldsymbol{u}^\star$ computed from (48), and the stabilized-regularized solution $\widetilde{\boldsymbol{u}}^\star_\gamma$ obtained from the solving of (22) by the Matlab command \, with respect to the parameter $\gamma \in [10^{-5},\, 1]$ (in A and C) and $\gamma \in [1,\, 10^5]$ (in B and D). The data $\mathbf{A}$, $\mathbf{L}$ and $\boldsymbol{b}$ are computed from the Matlab codes *shaw(n)* and *get_l*$(n, 2)$ provided in *Regularization Tools* [22] under the dimensions $m = n = 1\,000$ (with $\mathbf{cond}(\mathbf{A}) = 5.4935e^{+20}$ and $\mathbf{rank}(\mathbf{A}) = 20$), and the right hand term $\gamma \mathbf{L}^t \boldsymbol{g}$ in (22) is not taken into account since being considered as non available informations ($\boldsymbol{g} = \mathbf{0}$). The computations are made under parameters in (49) with the perturbed data (57), and the horizontal axis is displayed on a log-scale.

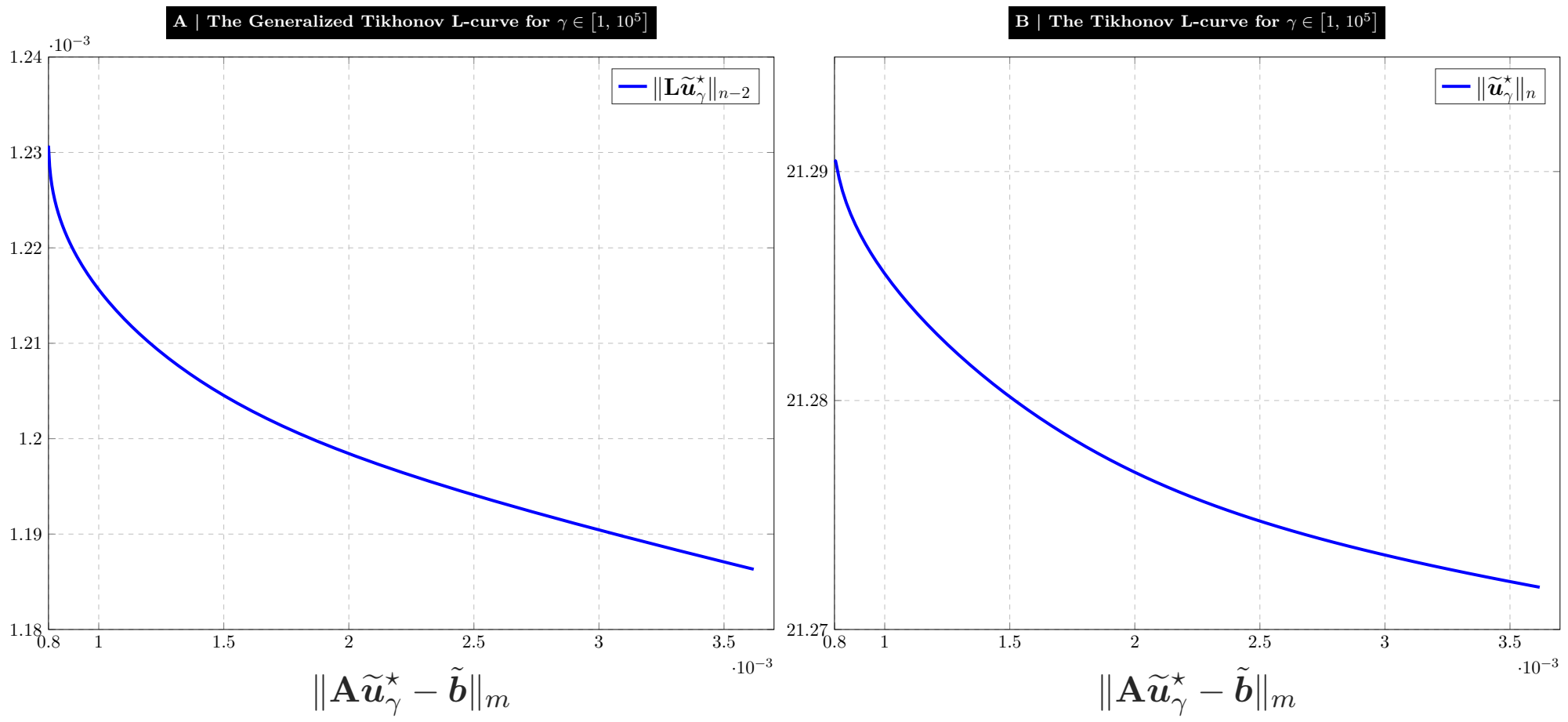

Figure 23: The residual norm $\|\mathbf{A}\widetilde{\boldsymbol{u}}^\star_\gamma - \tilde{\boldsymbol{b}}\|_m$ versus the semi-norm $\|\mathbf{L}\widetilde{\boldsymbol{u}}^\star_\gamma\|_{n-2}$ and the norm $\|\widetilde{\boldsymbol{u}}^\star_\gamma\|_n$ for $\gamma \in [1, 10^5]$ are represented in A and B, where the stabilized-regularized solution $\widetilde{\boldsymbol{u}}^\star_\gamma$ is obtained from the solving of (22) by the Matlab command \. The data $\mathbf{A}$, $\mathbf{L}$ and $\boldsymbol{b}$ are computed from the Matlab codes *shaw(n)* and *get_l*$(n, 2)$ provided in *Regularization Tools* [22] under the dimensions $m = n = 1\,000$ (with $\mathbf{cond}(\mathbf{A}) = 5.4935e^{+20}$ and $\mathbf{rank}(\mathbf{A}) = 20$). The computations are made under parameters in (49) with the perturbed data in (57), and the right hand term $\gamma \mathbf{L}^t \boldsymbol{g}$ in (22) is not taken into account since being considered as non available informations ($\boldsymbol{g} = \mathbf{0}$).

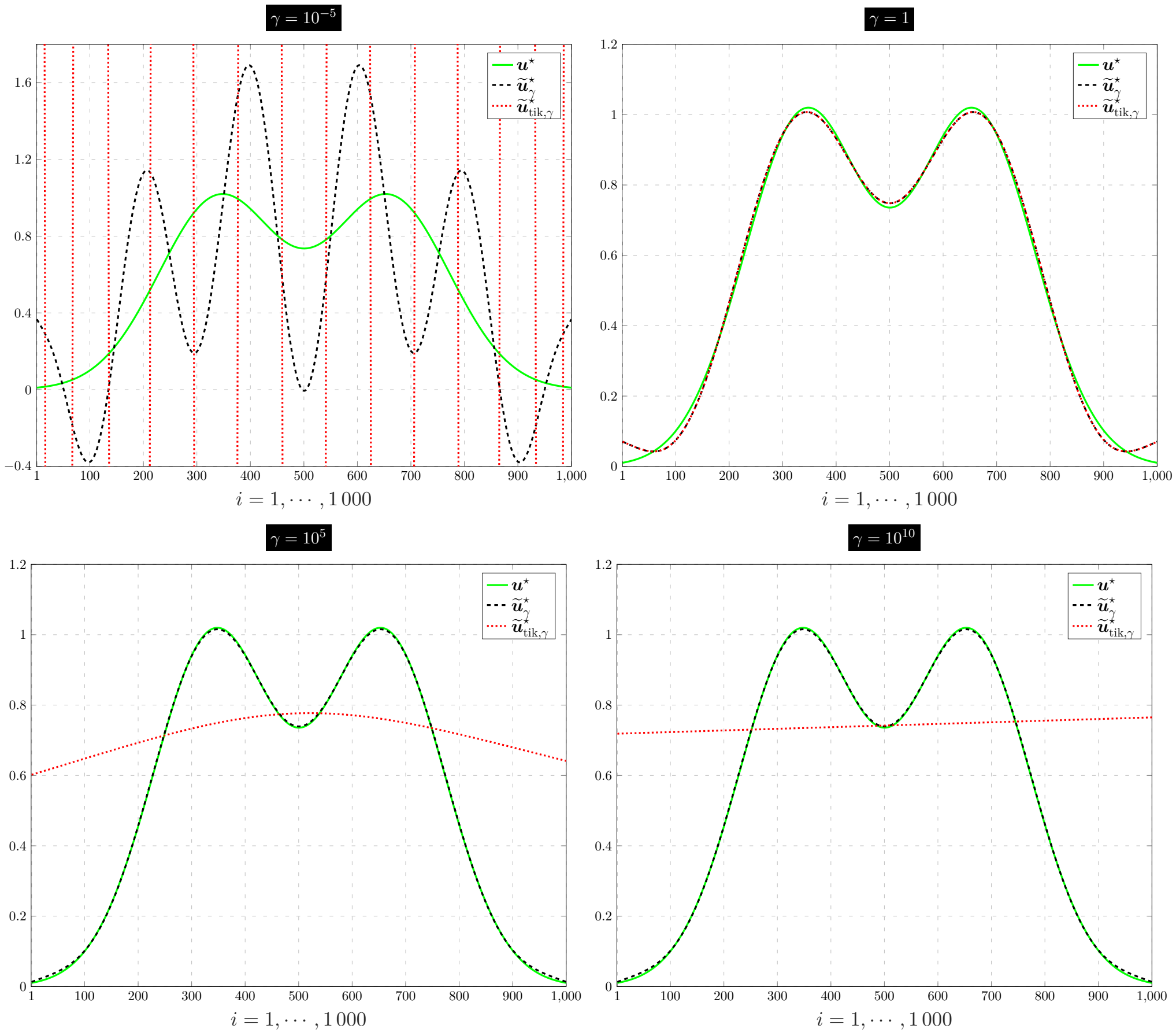

Figure 24: The exact solution $\boldsymbol{u}^\star$ (in green solid line) computed from (48), the stabilized-regularized solution $\widetilde{\boldsymbol{u}}^\star_\gamma$ (in black dashed line) obtained from the solving of (22) through the Matlab command \, and the Tikhonov regularization solution $\widetilde{\boldsymbol{u}}^\star_{\text{tik},\gamma}$ (in red dotted line) obtained in the solving of (5) by the code *tikhonov( )* implemented in *Regularization Tools* [22], are presented with respect to the discretization index $i = 1, \cdots, 1\,000$, for the different values of the regularization parameter $\gamma = 10^{-5}, 1, 10^5, 10^{10}$ (given in the top of the plots). The data $\mathbf{A}$, $\mathbf{L}$ and $\boldsymbol{b}$ are computed from the codes *shaw(n)* and *get_l*$(n, 2)$ provided in *Regularization Tools* [22] under the dimensions $m = n = 1\,000$ (with $\mathbf{cond}(\mathbf{A}) = 5.4935e^{+20}$ and $\mathbf{rank}(\mathbf{A}) = 20$). The computations are made under parameters in (49) with the perturbed data in (57), and the right hand term $\gamma \mathbf{L}^t \boldsymbol{g}$ in (22) is not taken into account since being considered as non available informations ($\boldsymbol{g} = \mathbf{0}$).

### *6.1.2.* Perturbations from smoothing of A and *b* in (58)

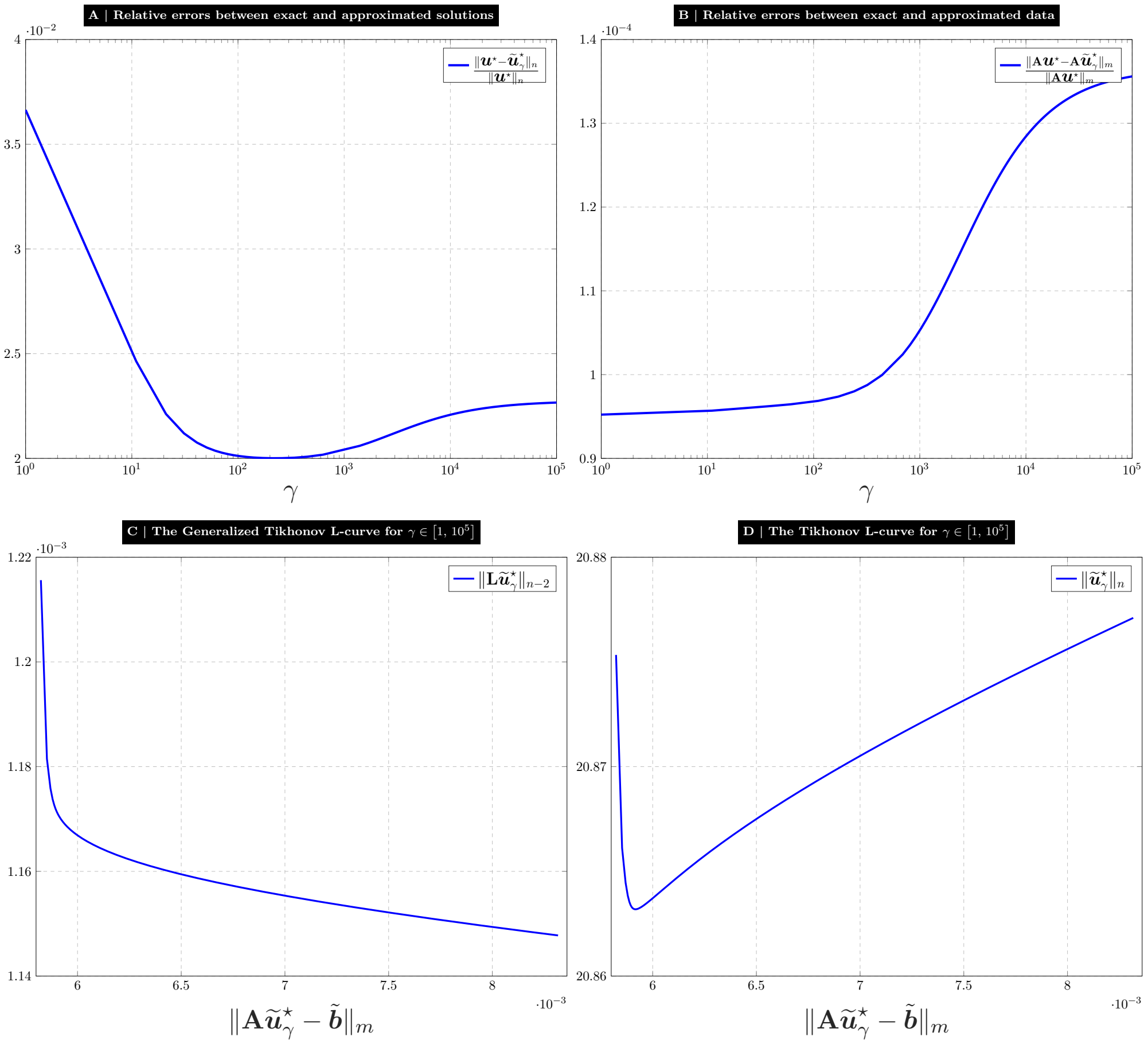

Figure 25: Relative errors between the exact solution $\boldsymbol{u}^{\star}$ computed from (48) and the stabilized-regularized solution $\widetilde{\boldsymbol{u}}_{\gamma}^{\star}$ obtained from the solving of (22) by the Matlab command \ are provided in A and B with respect to the parameter $\gamma \in [1, 10^5]$ (on a log-scale). The residual norm $\|\mathbf{A}\widetilde{\boldsymbol{u}}_{\gamma}^{\star} - \widetilde{\boldsymbol{b}}\|_m$ versus the semi-norm $\|\mathbf{L}\widetilde{\boldsymbol{u}}_{\gamma}^{\star}\|_{n-2}$ and the norm $\|\widetilde{\boldsymbol{u}}_{\gamma}^{\star}\|_n$ are presented in C and D for the range of values $\gamma \in [1, 10^5]$. In the solving of (22), the right hand term $\gamma \mathbf{L}^t \boldsymbol{g}$ is not taken account since being considered as non available informations ($\boldsymbol{g} = \mathbf{0}$). The computations are made under parameters in (49) with the perturbed data in (58), and with respect to the dimensions fixed to $m = n = 1\,000$, where $\mathbf{cond}(\mathbf{A}) = 2.1681e^{+21}$ and $\mathbf{rank}(\mathbf{A}) = 22$.

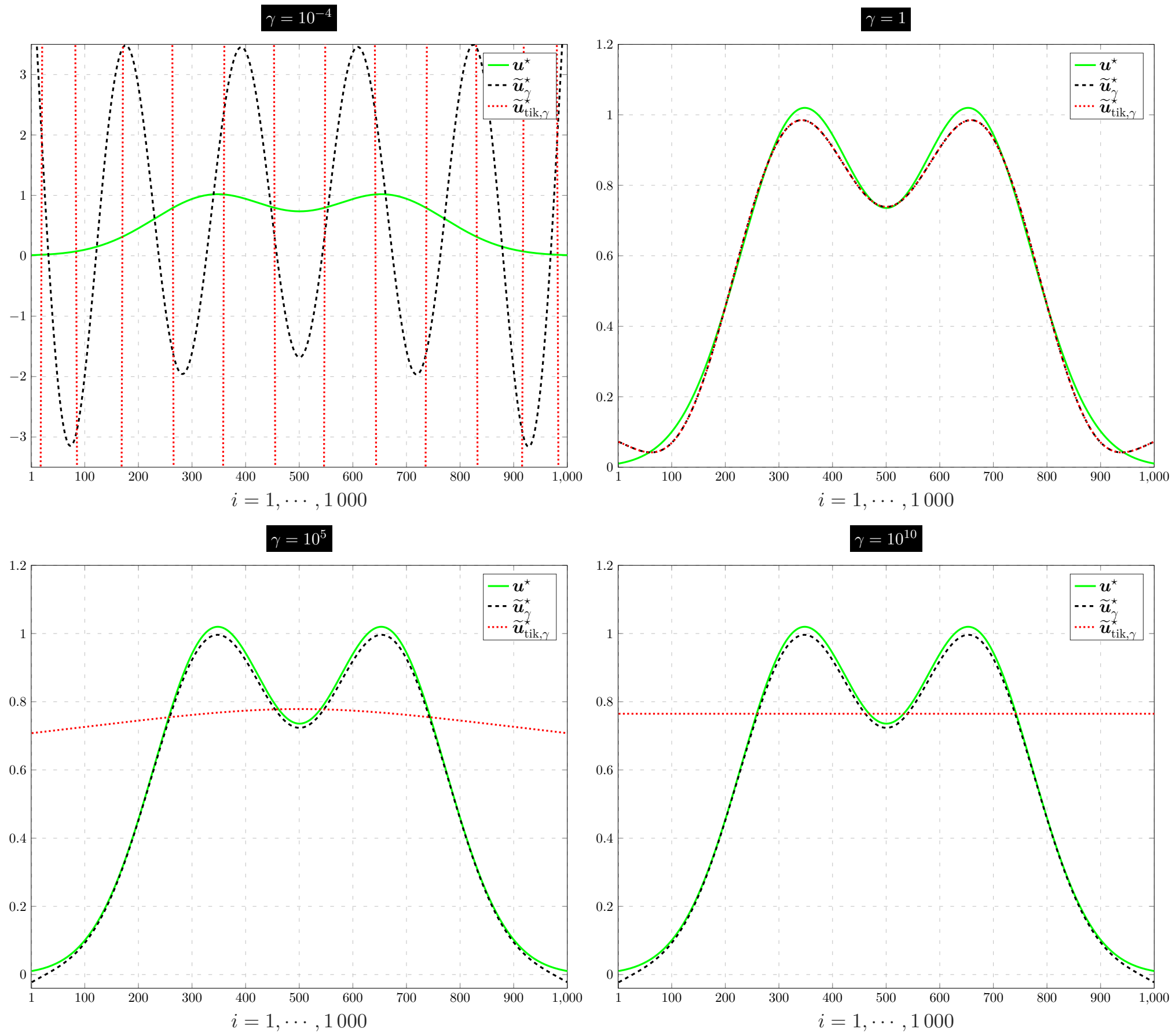

Figure 26: The exact solution $\boldsymbol{u}^\star$ (in green solid line) computed from (48), the stabilized-regularized solution $\widetilde{\boldsymbol{u}}^\star_\gamma$ (in black dashed line) obtained from the solving of (22) through the Matlab command \ and the Tikhonov regularization solution $\widetilde{\boldsymbol{u}}^\star_{\mathrm{tik},\gamma}$ (in red dotted line) obtained in the solving of (5) by the code *tikhonov( )* implemented in *Regularization Tools* [22] are presented with respect to the discretization index $i = 1, \cdots, 1\,000$, for the different parameter $\gamma = 10^{-4}, 1, 10^5, 10^{10}$ (given in the top of the plots). In (5) and (22), the right hand term $\gamma \mathbf{L}^t \boldsymbol{g}$ is not taken account since being considered as non available informations ($\boldsymbol{g} = \mathbf{0}$). The computations are made under parameters in (49) with the perturbed data in (58), and with respect to the dimensions fixed to $m = n = 1\,000$ where $\mathbf{cond}(\mathbf{A}) = 2.1681e^{+21}$ and $\mathbf{rank}(\mathbf{A}) = 22$.

### *6.2.* **Example 3**: *Phillips' famous test problem*

#### *6.2.1.* **Filtered and corrolated white noises in** (57)

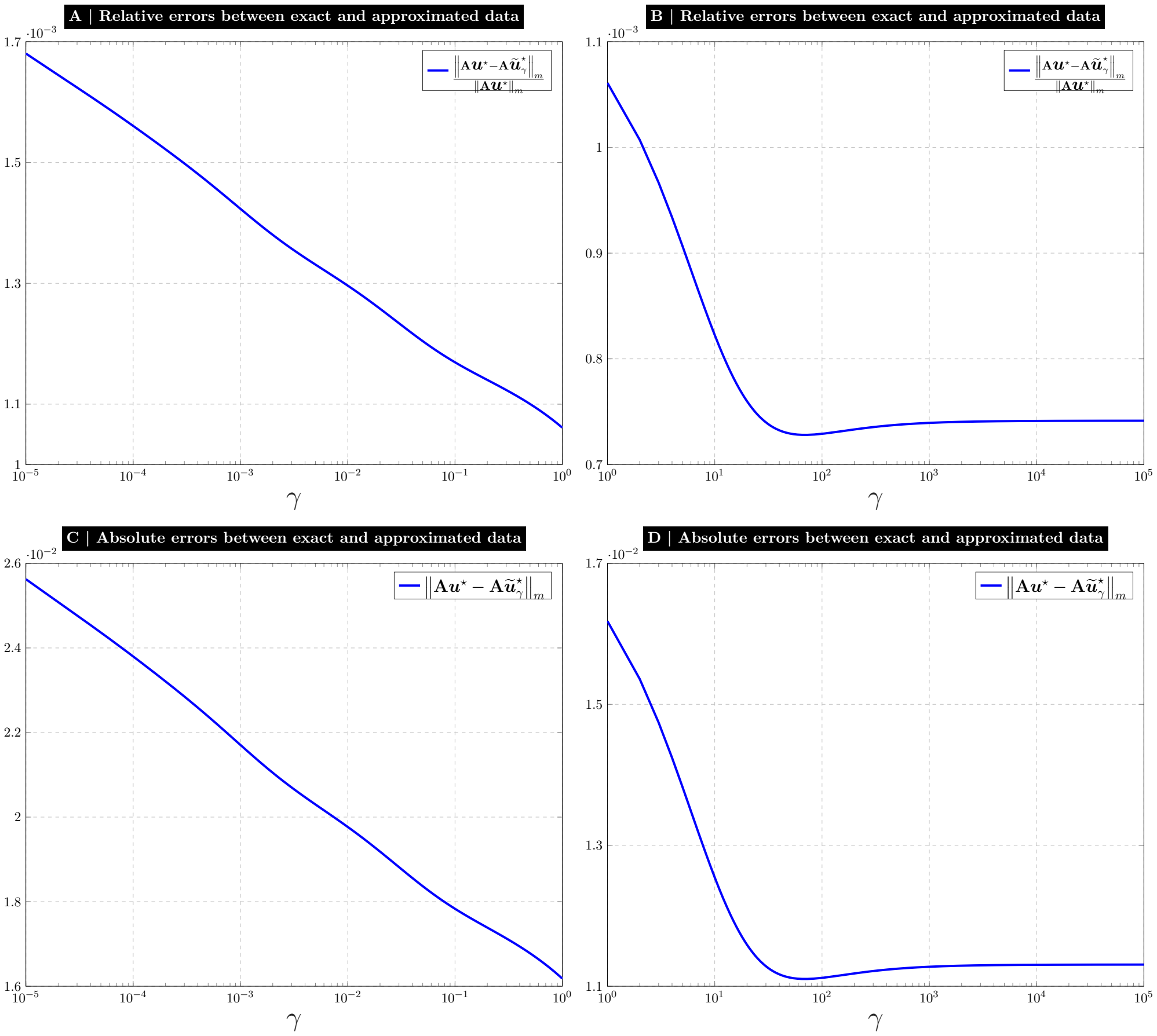

Figure 27: Absolute and relative errors between the exacte data $\mathbf{A}\boldsymbol{u}^\star$ obtained from (66) and our stabilized-regularized data $\mathbf{A}\widetilde{\boldsymbol{u}}^\star_\gamma$ obtained from the solving of (22) by the Matlab command \, with respect to the parameter $\gamma \in [10^{-5}, 1]$ (in A and C) and $\gamma \in [1, 10^5]$ (in B and D). The data $\mathbf{A}$, $\mathbf{L}$ and $\boldsymbol{b}$ are obtained from the Matlab codes *phillips(n)* and *get_l*$(n, 1)$ provided in *Regularization Tools* [22], the computations are made under the perturbed data in (57), the right hand term $\gamma \mathbf{L}^t \boldsymbol{g}$ in (22) is not taken account since being considered as non available informations ($\boldsymbol{g} = \mathbf{0}$), and the dimensions are fixed to $m = n = 1\,000$ where $\mathbf{cond}(\mathbf{A}) = 2.8990e^{+10}$ and $\mathbf{rank}(\mathbf{A}) = 1\,000$.

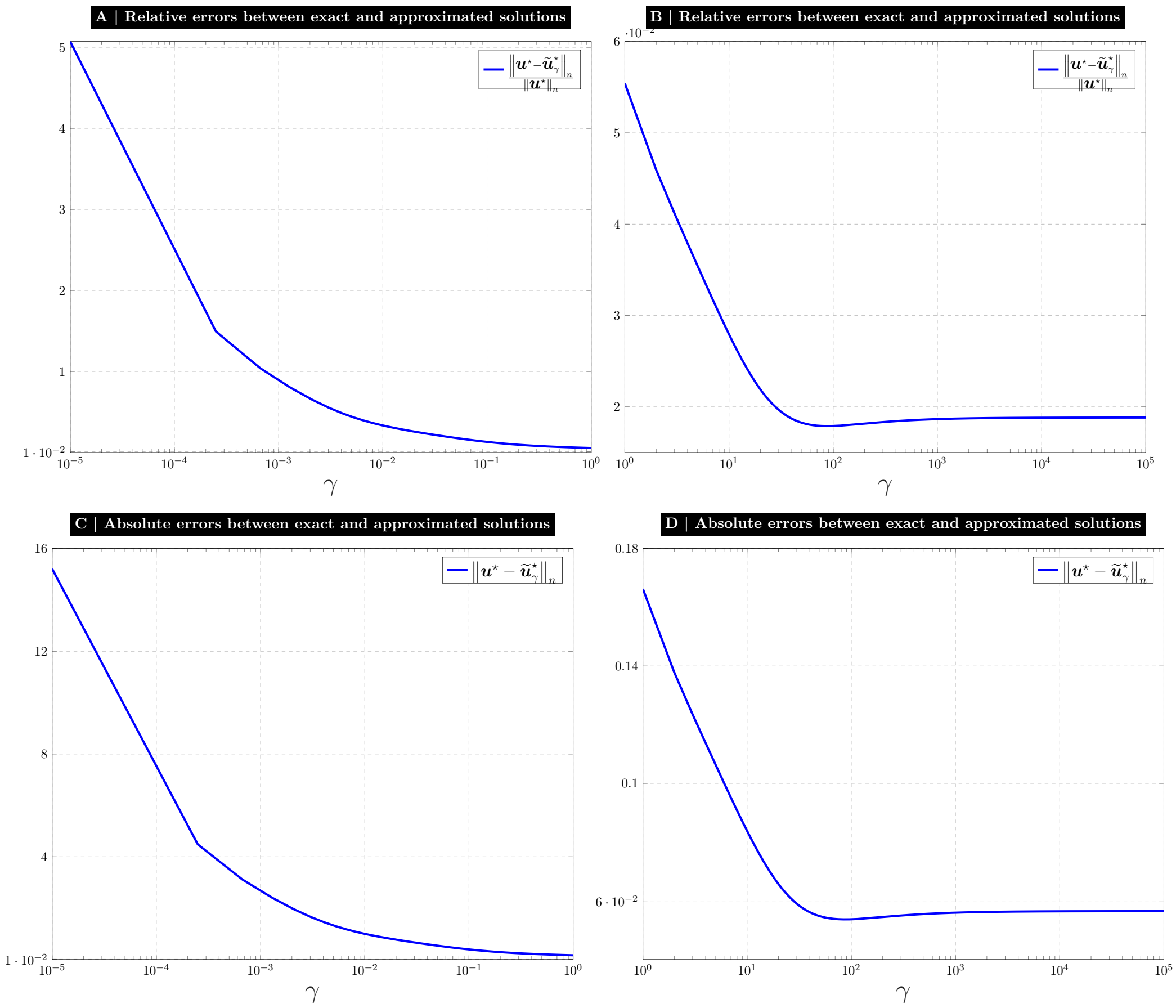

Figure 28: Relative and absolute errors between the exacte solution $\boldsymbol{u}^\star$ computed from (66) and our stabilized-regularized solution $\widetilde{\boldsymbol{u}}_\gamma^\star$ obtained from the solving of (22) by the Matlab command \, with respect to the parameter $\gamma \in [10^{-5},\, 1]$ (in A and C) and $\gamma \in [1,\, 10^5]$ (in B and D). The data $\mathbf{A}$, $\mathbf{L}$ and $\boldsymbol{b}$ are obtained from the Matlab codes *phillips(n)* and *get_l*$(n, 1)$ provided in *Regularization Tools* [22], the computations are made under perturbed data in (57), the right hand term $\gamma \mathbf{L}^t \boldsymbol{g}$ in (22) is not taken account since being considered as non available informations ($\boldsymbol{g} = \mathbf{0}$), and the dimensions are fixed to $m = n = 1\,000$ where $\mathbf{cond}(\mathbf{A}) = 2.8990e^{+10}$ and $\mathbf{rank}(\mathbf{A}) = 1\,000$.

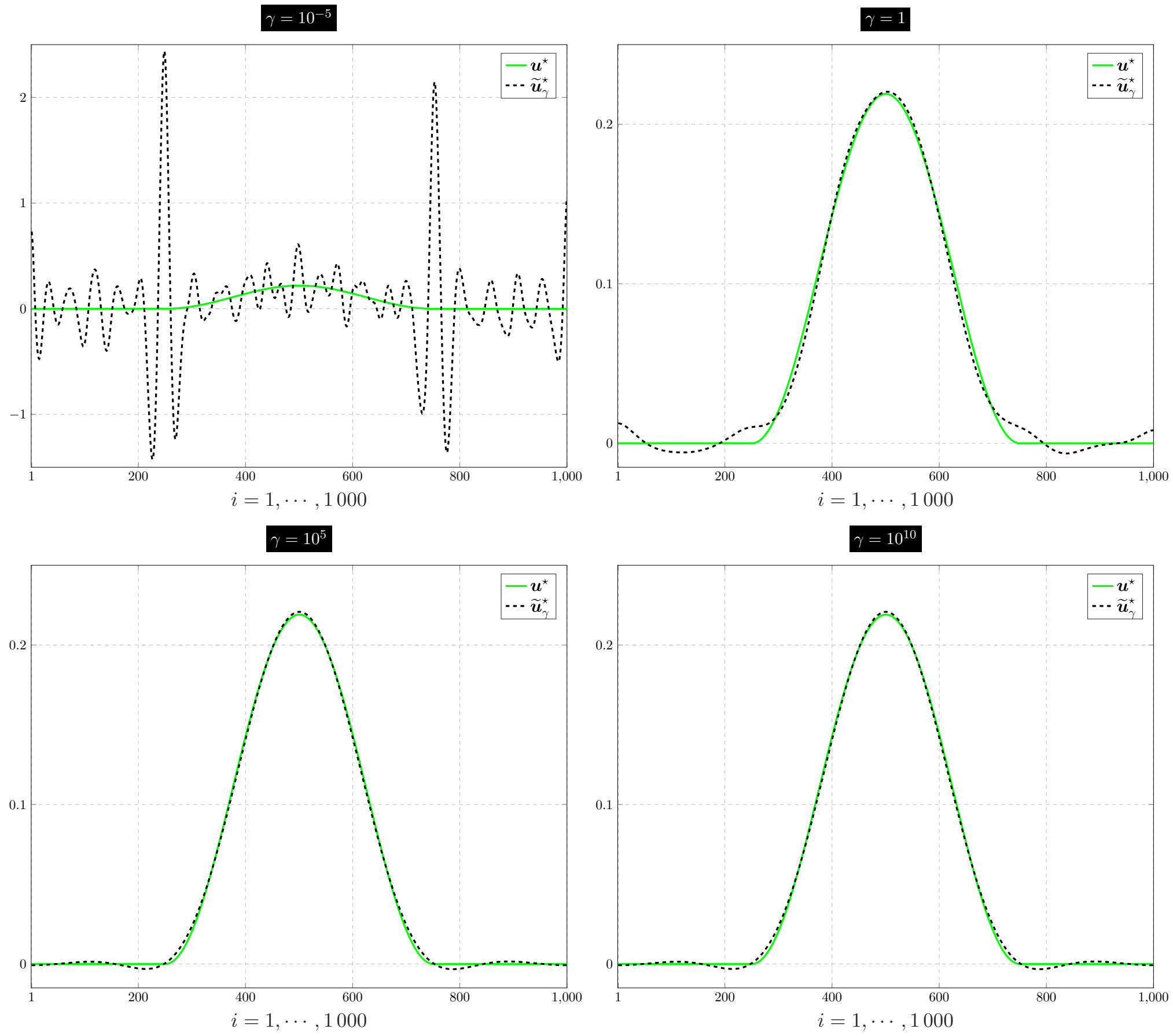

Figure 29: The exact solutions $\boldsymbol{u}^\star$ (in green solid line) and the stabilized-regularized solutions $\widetilde{\boldsymbol{u}}_\gamma^\star$ (in black dashed line) represented with respect to the discretization index $i = 1, \cdots, 500$, for different values of the regularization parameter $\gamma = 10^{-5}, 1, 10^5, 10^{10}$ (given in the top of the plots). The stabilized-regularized solution $\widetilde{\boldsymbol{u}}_\gamma^\star$ is computed from the solving of (22) by the Matlab command \ under the perturbed data in (57), and the right hand term $\gamma \mathbf{L}^t \boldsymbol{g}$ in (22) is not taken account since being considered as non available informations ($\boldsymbol{g} = \mathbf{0}$). The dimensions are fixed to $m = n = 1\,000$, where $\mathbf{cond}(\mathbf{A}) = 2.8990e^{+10}$ and $\mathbf{rank}(\mathbf{A}) = 1\,000$.

### *6.2.2.* Perturbations from smoothing of A and $\boldsymbol{b}$ in (58)

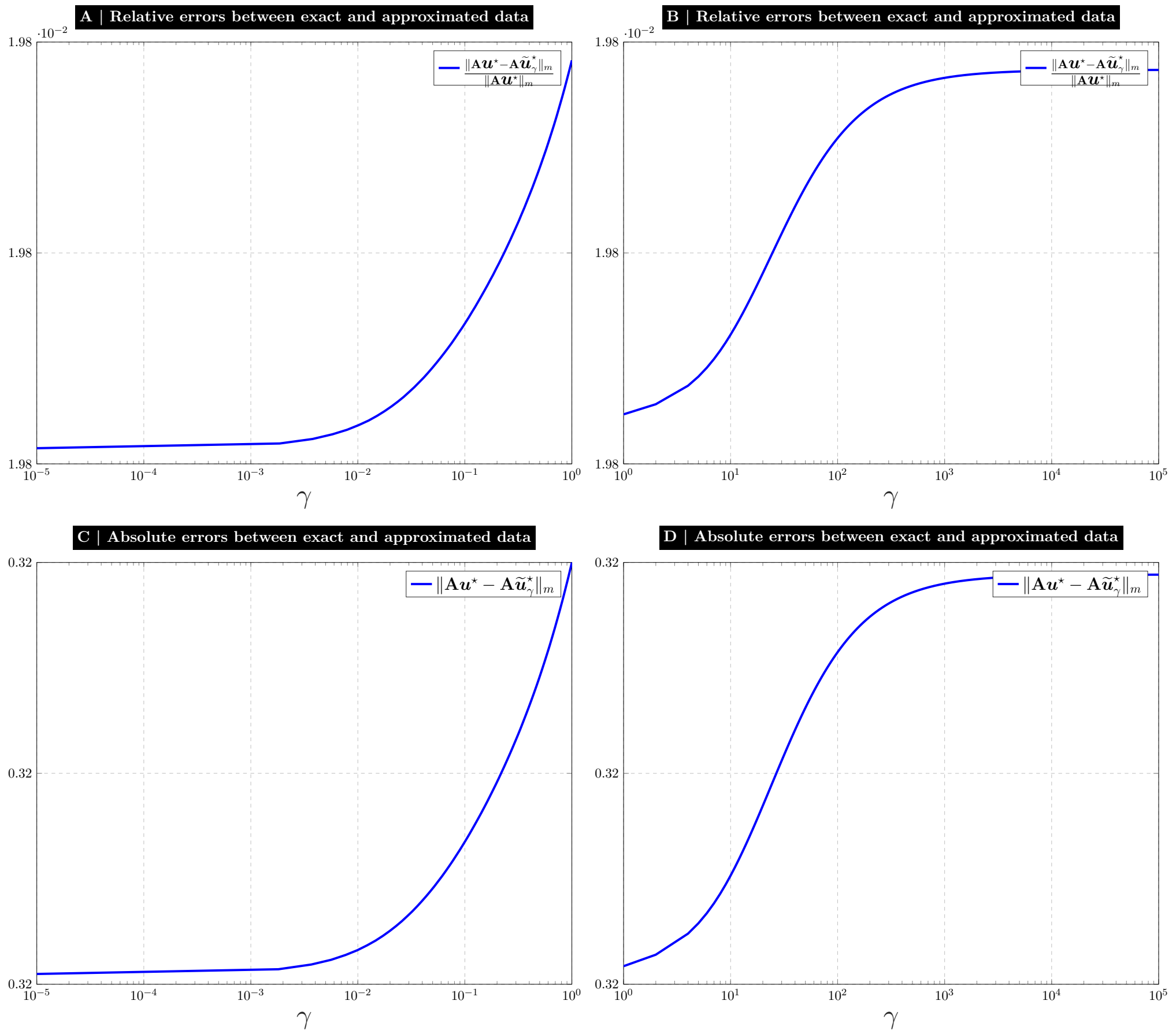

Figure 30: Absolute and relative errors between the exacte data $\mathbf{A}\boldsymbol{u}^\star$ obtained from (66) and our stabilized-regularized data $\mathbf{A}\widetilde{\boldsymbol{u}}^\star_\gamma$ obtained from the solving of (22) by the Matlab command \, with respect to the parameter $\gamma \in [10^{-5},\, 1]$ (in A and C) and $\gamma \in [1,\, 10^5]$ (in B and D). The data $\mathbf{A}$, $\mathbf{L}$ and $\boldsymbol{b}$ are obtained from the Matlab codes *phillips(n)* and *get_l*$(n, 1)$ provided in *Regularization Tools* [22], the computations are made under the perturbed data in (58), the right hand term $\gamma \mathbf{L}^t \boldsymbol{g}$ in (22) is not taken account since being considered as non available informations ($\boldsymbol{g} = \mathbf{0}$), and the dimensions are fixed to $m = n = 1\,000$, where $\mathbf{cond}(\mathbf{A}) = 2.7583e^{+11}$ and $\mathbf{rank}(\mathbf{A}) = 1\,000$.

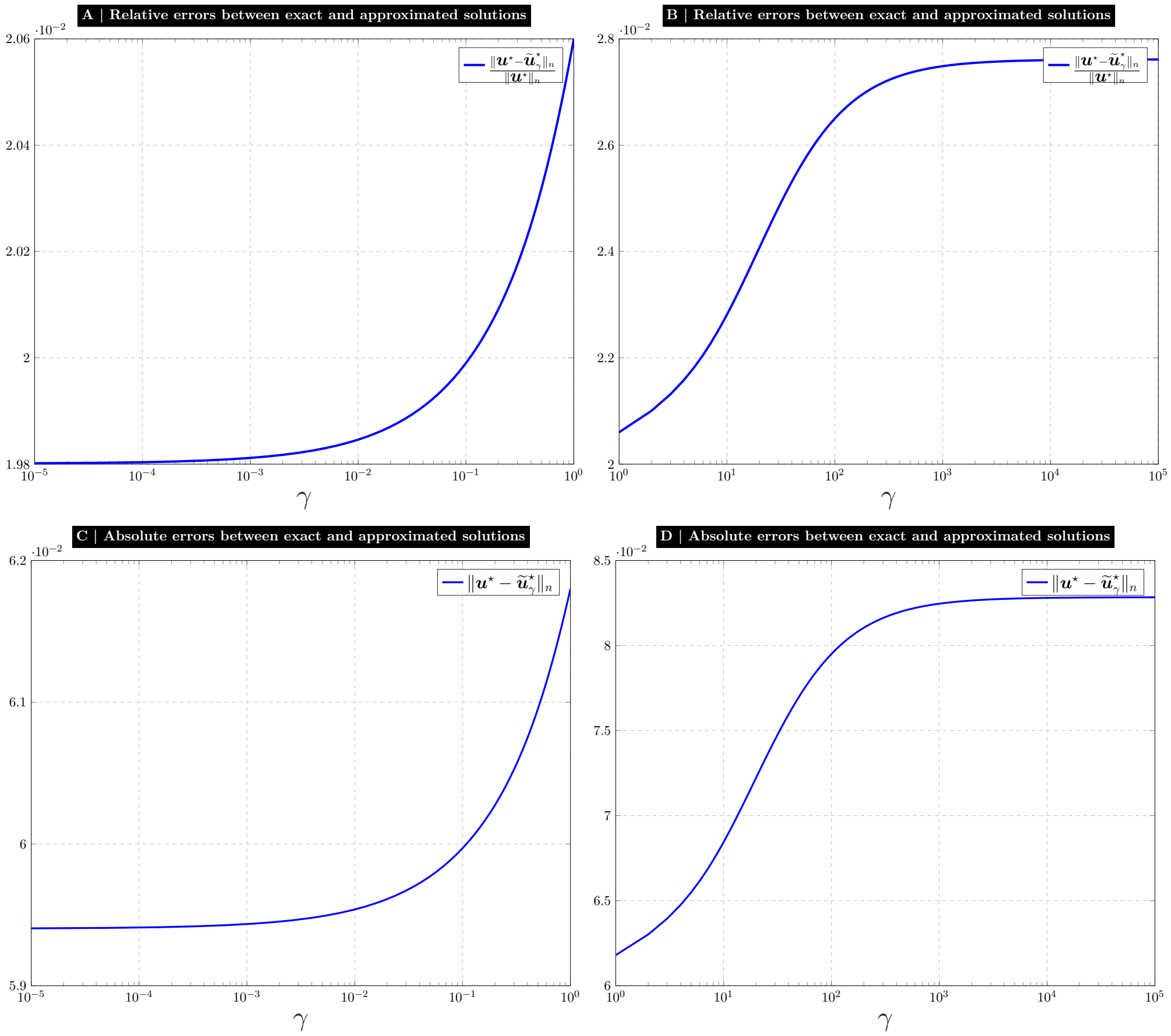

Figure 31: Relative and absolute errors between the exacte solution $\boldsymbol{u}^\star$ computed from (66) and our stabilized-regularized solution $\widetilde{\boldsymbol{u}}^\star_\gamma$ obtained from the solving of (22) by the Matlab command \, with respect to the parameter $\gamma \in [10^{-5},\, 1]$ (in A and C) and $\gamma \in [1,\, 10^5]$ (in B and D). The data $\mathbf{A}$, $\mathbf{L}$ and $\boldsymbol{b}$ are obtained from the Matlab codes *phillips(n)* and *get_l*$(n, 1)$ provided in *Regularization Tools* [22], the computations are made under perturbed data in (58), the right hand term $\gamma \mathbf{L}^t \boldsymbol{g}$ in (22) is not taken account since being considered as non available informations ($\boldsymbol{g} = \mathbf{0}$), and the dimensions are fixed to $m = n = 1\,000$, where $\mathbf{cond}(\mathbf{A}) = 2.7583e^{+11}$ and $\mathbf{rank}(\mathbf{A}) = 1\,000$.

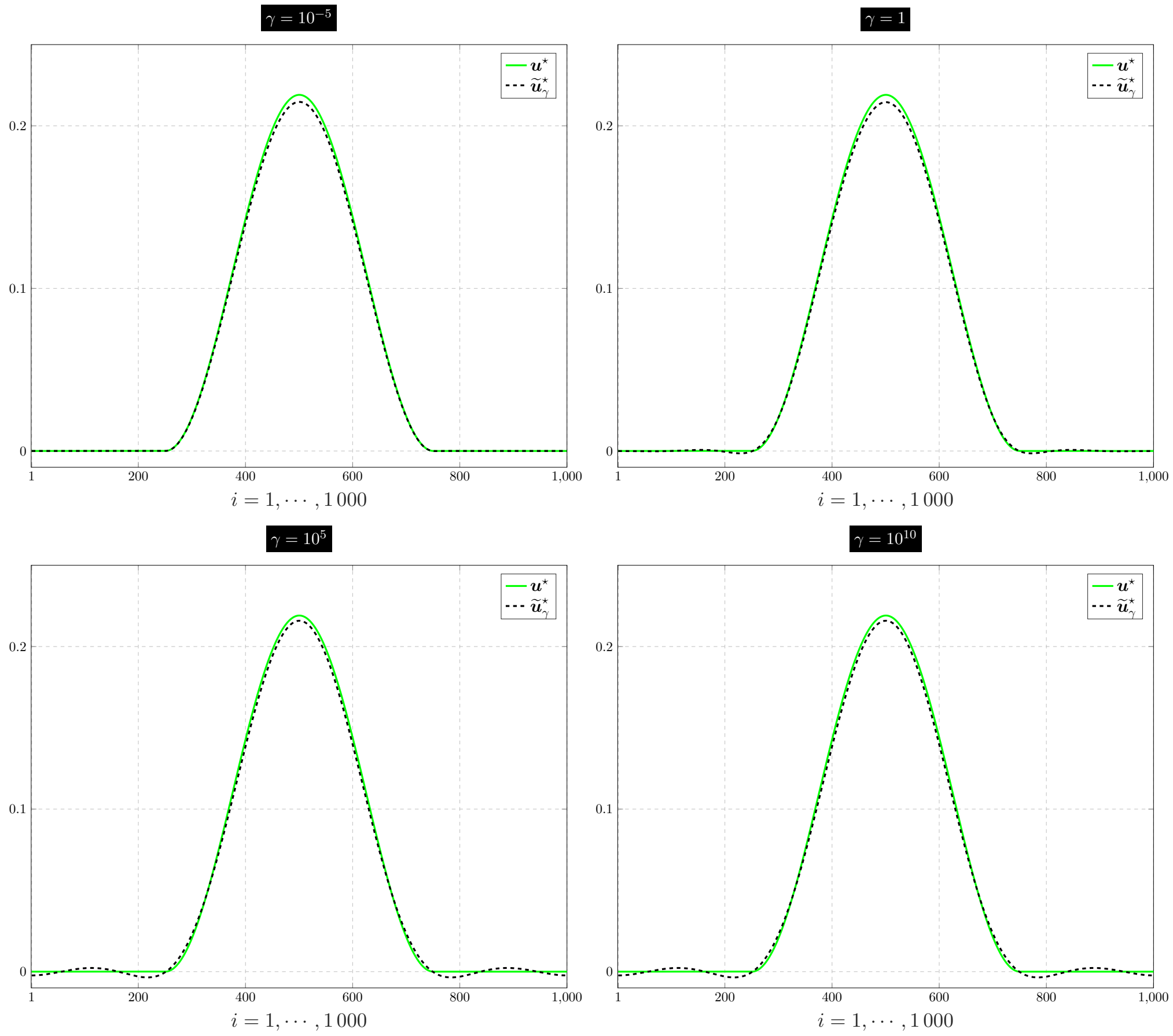

Figure 32: The exact solutions $\boldsymbol{u}^\star$ (in green solid line) and the stabilized-regularized solutions $\widetilde{\boldsymbol{u}}^\star_\gamma$ (in black dashed line) represented with respect to the discretization index $i = 1, \cdots, 500$, for different values of the regularization parameter $\gamma = 10^{-5}, 1, 10^5, 10^{10}$ (given in the top of the plots). The stabilized-regularized solution $\widetilde{\boldsymbol{u}}^\star_\gamma$ is computed from the solving of (22) by the Matlab command \ under the perturbed data in (58), and the right hand term $\gamma \mathbf{L}^t \boldsymbol{g}$ in (22) is not taken account since being considered as non available informations ($\boldsymbol{g} = \mathbf{0}$). The dimensions are fixed to $m = n = 1\,000$, where $\mathbf{cond}(\mathbf{A}) = 2.7583e^{+11}$ and $\mathbf{rank}(\mathbf{A}) = 1\,000$.

## 7. Supplement tests

### 7.1. *The gravity surveying model problem (see [21, 22, 27])*

The gravity surveying model problem is defined by the Fredholm integral equation

$$\int_0^1 K(s,x)u^\star(x)dx = g(s), \;\; a \le s \le b,$$

where the following analytical solution and the kernel $K(s,x)$ are considered

$$u^\star(x) = sin(\pi x) + \frac{1}{2}sin(2\pi x), \; 0 \le x \le 1, \tag{68}$$
$$K(s,x) = \frac{h}{(h^2+(s-x)^2)^{\frac{3}{2}}}, \;\; a \le s \le b, \; 0 \le x \le 1.$$

The problem is discretized by means of the midpoint quadrature rule with $n$ points, and the default integration interval $[0,1]$ leads to a symmetric Toeplitz matrix $\mathbf{A}$. Then the right-hand side is computed as $\boldsymbol{b} = \mathbf{A}\boldsymbol{u}^\star$, where the $n$-vector exact solution $\boldsymbol{u}^\star = (u_j^\star)_{1\le j\le n}$ is derived from the analytical formulas (68) by $u_j^\star = u^\star(x_j^\star)$, $\;i = 1,\cdots,n$. The parameter $h$ is the depth at which the magnetic deposit is located, and the default value is $h = 0.25$. The larger the depth $h$, the faster the decay of the singular values. In the forthcoming computations, we fix $a = 0$, $b = 1$ and choose the second derivative operator in (55) as regularization matrix.

#### 7.1.1. *Perturbations by white noises in (56)*

- $h = 0.25, \quad \eta = 10^{-3}$

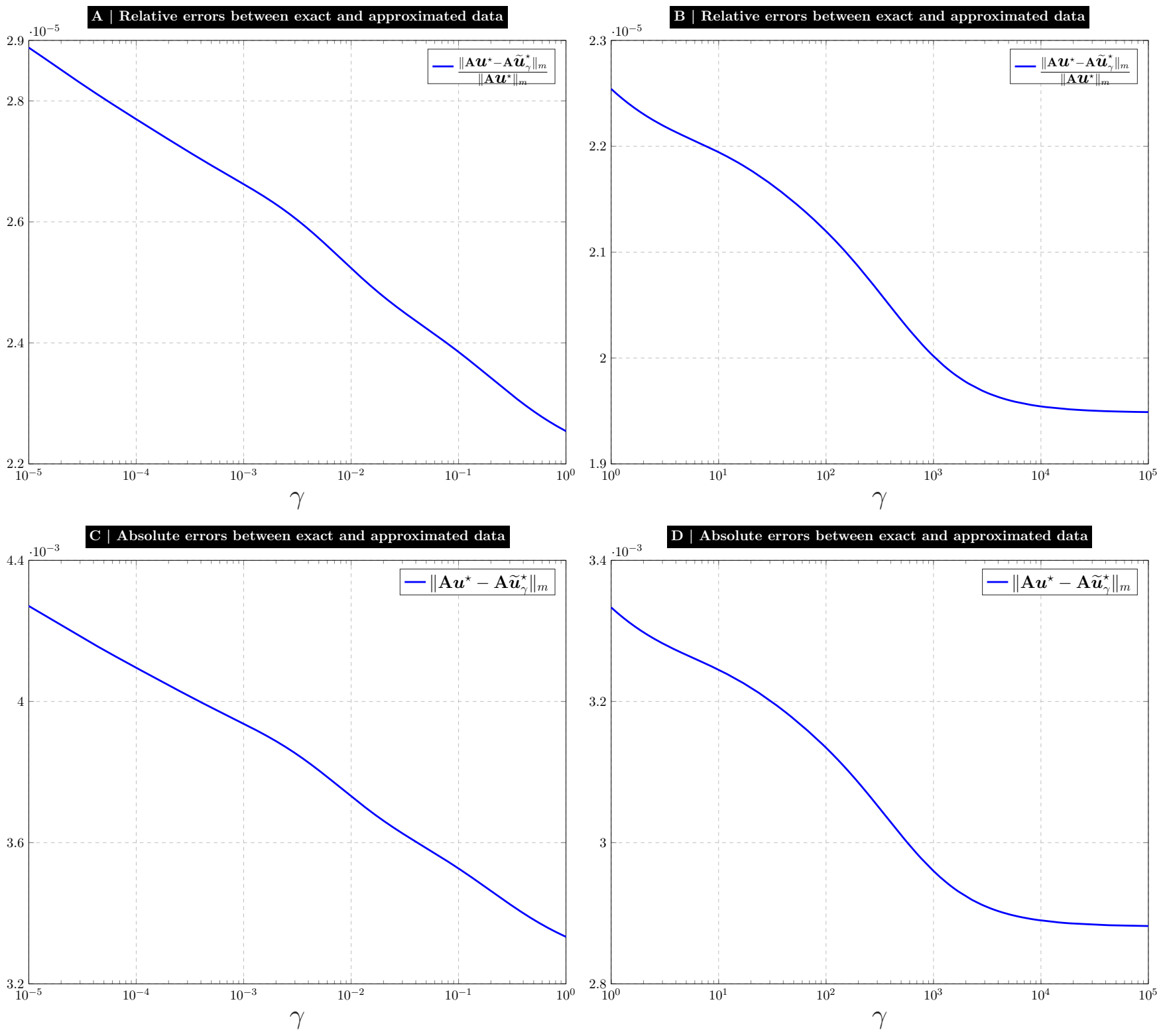

Figure 33: Absolute and relative errors between the exact data $\mathbf{A}\boldsymbol{u}^\star$ computed from (68), and the stabilized-regularized data $\mathbf{A}\widetilde{\boldsymbol{u}}_\gamma^\star$ obtained from the solving of (22) by the Matlab command \, with respect to the parameter $\gamma \in [10^{-5}, 1]$ (in A and C) and $\gamma \in [1, 10^5]$ (in B and D). The data $\mathbf{A}$, $\mathbf{L}$ and $\boldsymbol{b}$ are computed from the Matlab codes *gravity(n,1,a,b,h)* and *get_l*$(n,2)$ provided in *Regularization Tools* [22] under the dimensions $n = 1\,000$ (where $\mathbf{cond}(\mathbf{A}) = 2.5999e^{20+}$, $\mathbf{rank}(\mathbf{A}) = 45$), and the right hand term $\gamma\mathbf{L}^t\boldsymbol{g}$ in (22) is not taken into account since $\boldsymbol{g} = \mathbf{0}$. The computations are made under the perturbed data in (56) and the horizontal axis is displayed on a log-scale.

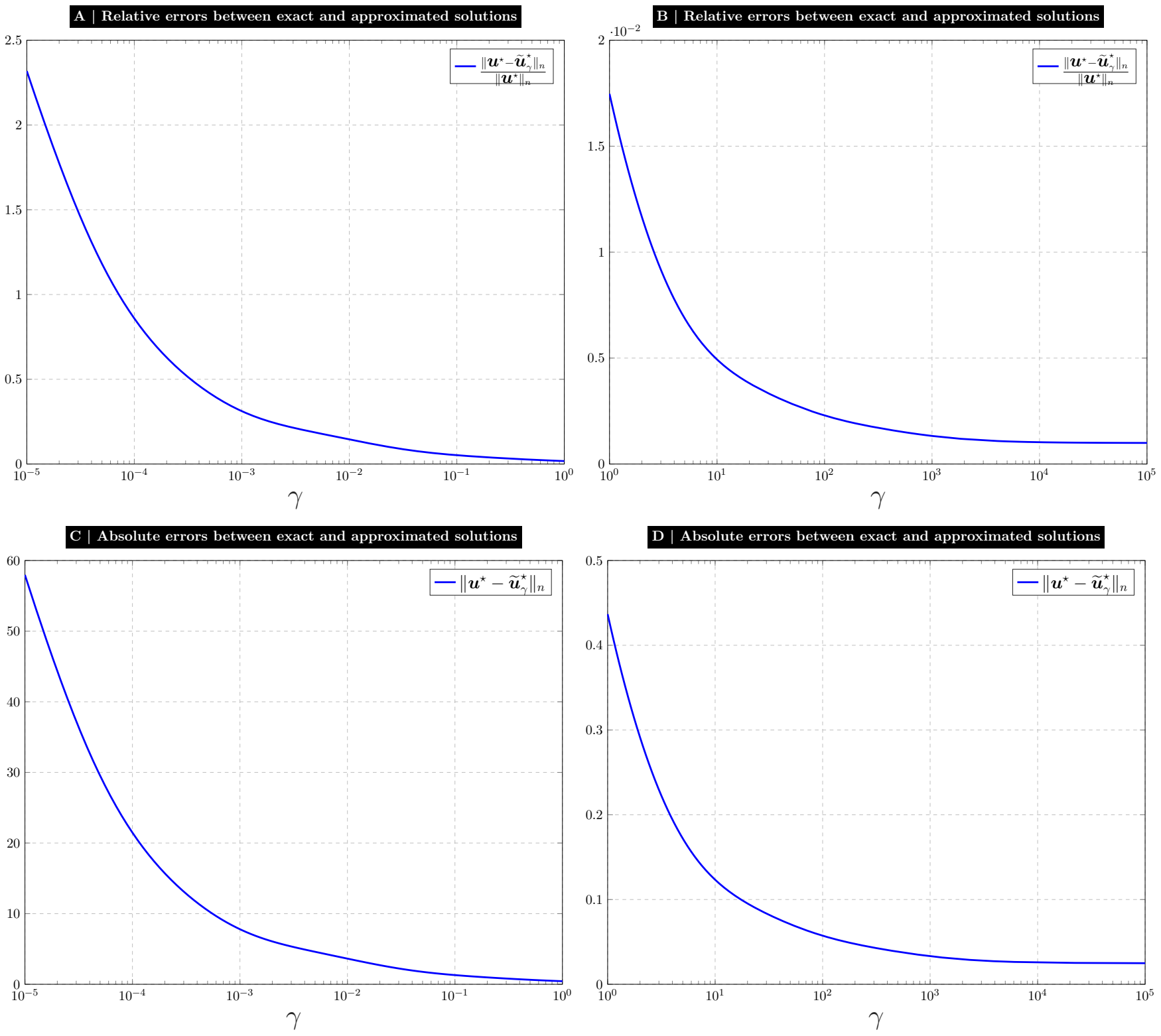

Figure 34: Absolute and relative errors between the exact solution $\boldsymbol{u}^\star$ computed from (68), and the stabilized-regularized solution $\widetilde{\boldsymbol{u}}^\star_\gamma$ obtained from the solving of (22) by the Matlab command \, with respect to the parameter $\gamma \in [10^{-5}, 1]$ (in A and C) and $\gamma \in [1, 10^5]$ (in B and D). The data $\mathbf{A}$, $\mathbf{L}$ and $\boldsymbol{b}$ are computed from the Matlab codes *gravity(n,1,a,b,h)* and *get_l*$(n, 2)$ provided in *Regularization Tools* [22] under the dimensions $n = 1\,000$ (where $\mathbf{cond}(\mathbf{A}) = 2.5999e^{20+}$, $\mathbf{rank}(\mathbf{A}) = 45$), and the right hand term $\gamma \mathbf{L}^t \boldsymbol{g}$ in (22) is not taken into account since $\boldsymbol{g} = \mathbf{0}$. The computations are made under the perturbed data in (56) and the horizontal axis is displayed on a log-scale.

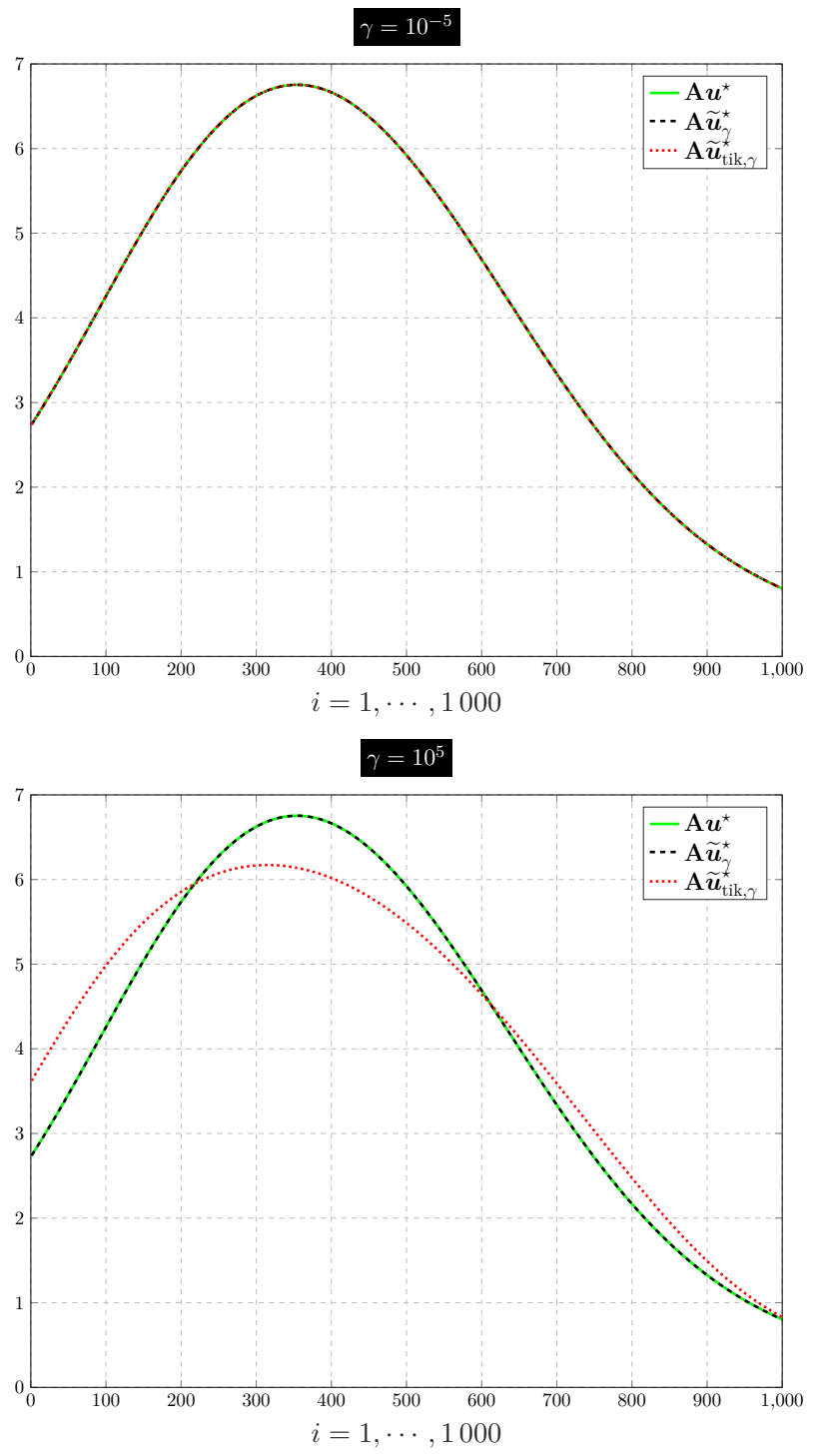

Figure 35: The exact data $\mathbf{A}\boldsymbol{u}^\star$ (in green solid line), the data reconstructed from the stabilized-regularized approach $\mathbf{A}\widetilde{\boldsymbol{u}}^\star_\gamma$ (in black dashed line), and the data made from the Tikhonov regularization $\mathbf{A}\widetilde{\boldsymbol{u}}^\star_{\mathrm{tik},\gamma}$ (in red dotted line), are presented with respect to the discretisation index $i = 1, \cdots, 1\,000$, for the different values of the regularization parameter $\gamma = 10^{-5}, 1, 10^5, 10^{10}$ (given in the top of the plots). The stabilized-regularized solution $\widetilde{\boldsymbol{u}}^\star_\gamma$ is obtained from the solving of (22) through the Matlab command \, and the Tikhonov regularization solution $\widetilde{\boldsymbol{u}}^\star_{\mathrm{tik},\gamma}$ is computed in the solving of (5) by the Matlab code *tikhonov( )* implemented in *Regularization Tools* [22]. The data $\mathbf{A}$, $\mathbf{L}$ and $\boldsymbol{b}$ are computed from the Matlab codes *gravity(n,1,a,b,h)* and *get_l*$(n, 2)$ in *Regularization Tools* [22] and the approximate derivative is fixed to $\boldsymbol{g} = \mathbf{0}$. The computations are made under the perturbed data in (56), with respect to the dimensions $n = 1\,000$ (where $\mathbf{cond}(\mathbf{A}) = 2.5999e^{20+}$, $\mathbf{rank}(\mathbf{A}) = 45$).

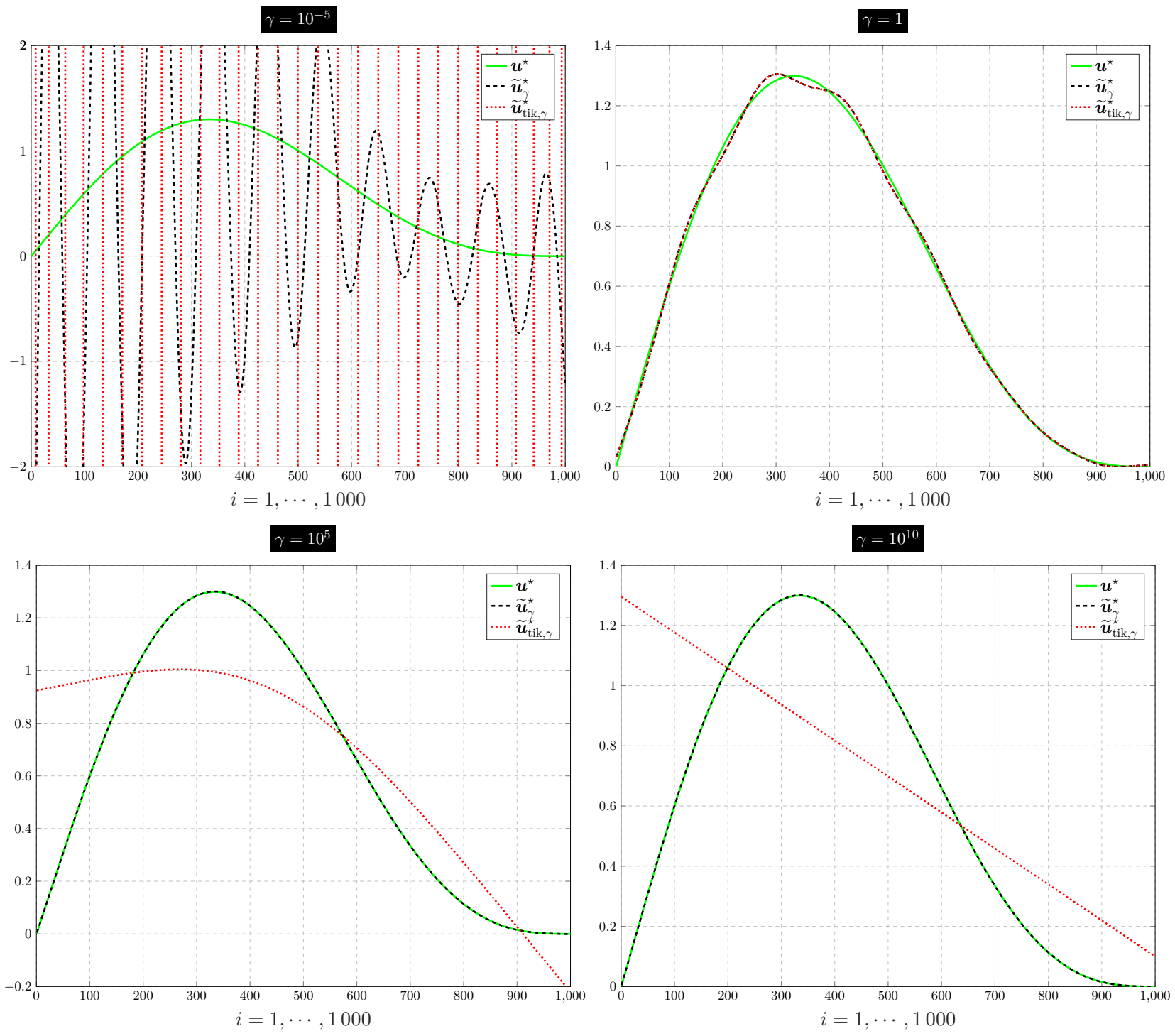

Figure 36: The exact solution $\boldsymbol{u}^\star$ (in green solid line) computed from (68), the stabilized-regularized solution $\widetilde{\boldsymbol{u}}^\star_\gamma$ (in black dashed line) obtained from the solving of (22) by the Matlab command \, and the Tikhonov regularization solution $\widetilde{\boldsymbol{u}}^\star_{\mathrm{tik},\gamma}$ (in red dotted line) obtained in the solving of (5) by the Matlab code *tikhonov( )* implemented in *Regularization Tools* [22], are presented with respect to the discretization index $i = 1, \cdots, 1\,000$, for different values of the parameter $\gamma = 10^{-5}, 1, 10^5, 10^{10}$. The data $\mathbf{A}$, $\mathbf{L}$ and $\boldsymbol{b}$ are computed from the Matlab codes *gravity(n,1,a,b,h)* and *get_l*$(n, 2)$ in *Regularization Tools* [22] and the term $\gamma \mathbf{L}^t \boldsymbol{g}$ in (5) and (22) is not taken into account ($\boldsymbol{g} = \mathbf{0}$). The computations are made under perturbed data in (56), with dimensions $n = 1\,000$ (where $\mathbf{cond}(\mathbf{A}) = 2.5999e^{20+}$, $\mathbf{rank}(\mathbf{A}) = 45$).

- $h = 0.25, \quad \eta = 10^{-2}$

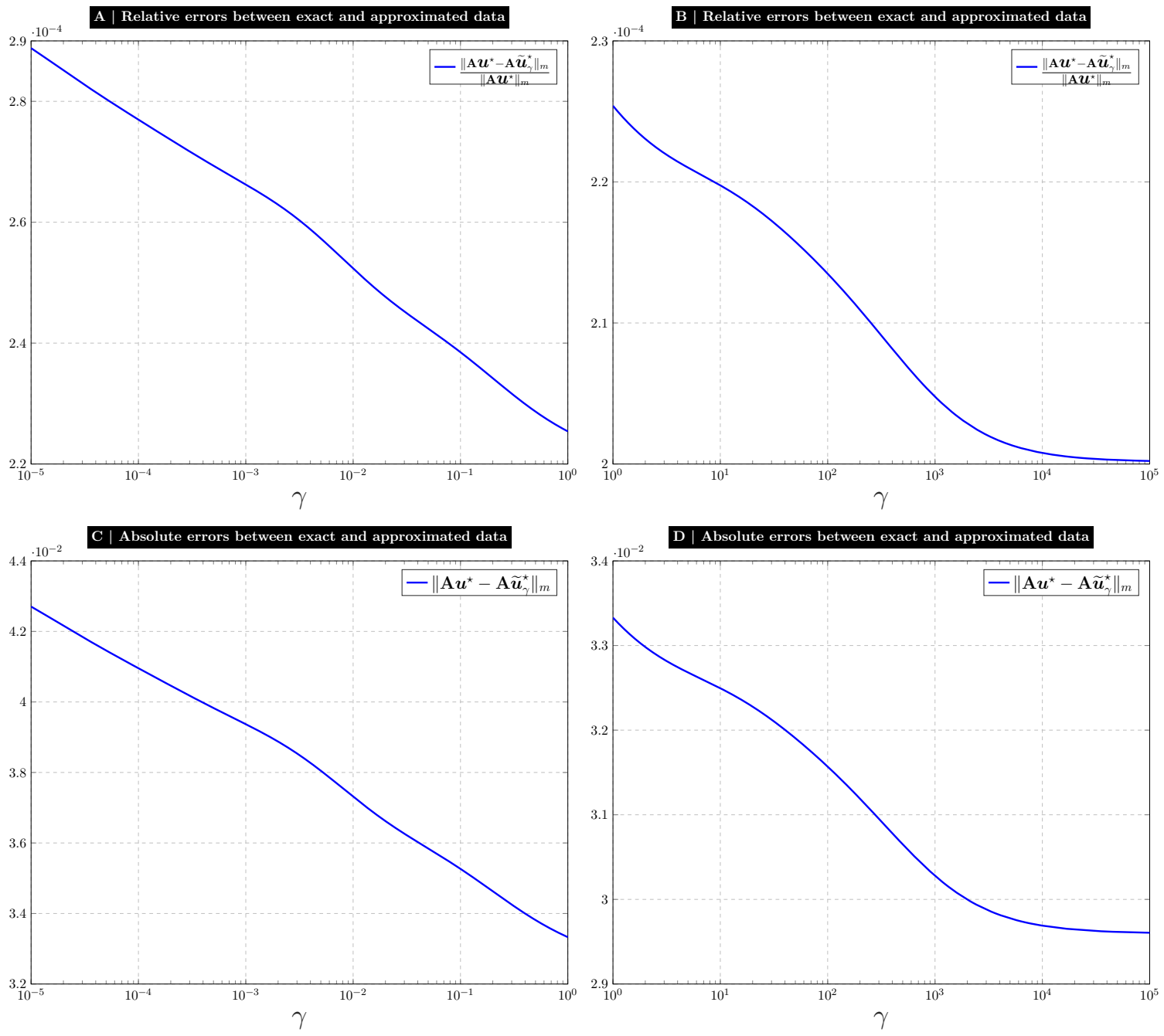

Figure 37: Absolute and relative errors between the exact data $\mathbf{A}\boldsymbol{u}^\star$ computed from (68), and the stabilized-regularized data $\mathbf{A}\widetilde{\boldsymbol{u}}^\star_\gamma$ obtained from the solving of (22) by the Matlab command \, with respect to the parameter $\gamma \in [10^{-5}, 1]$ (in A and C) and $\gamma \in [1, 10^5]$ (in B and D). The data $\mathbf{A}$, $\mathbf{L}$ and $\boldsymbol{b}$ are computed from the Matlab codes *gravity(n,1,a,b,h)* and *get_l*$(n, 2)$ provided in *Regularization Tools* [22] under the dimensions $n = 1\,000$ (where $\mathbf{cond}(\mathbf{A}) = 2.5999e^{20+}$, $\mathbf{rank}(\mathbf{A}) = 45$), and the right hand term $\gamma \mathbf{L}^t \boldsymbol{g}$ in (22) is not taken into account since $\boldsymbol{g} = \mathbf{0}$. The computations are made under the perturbed data in (56) and the horizontal axis is displayed on a log-scale.

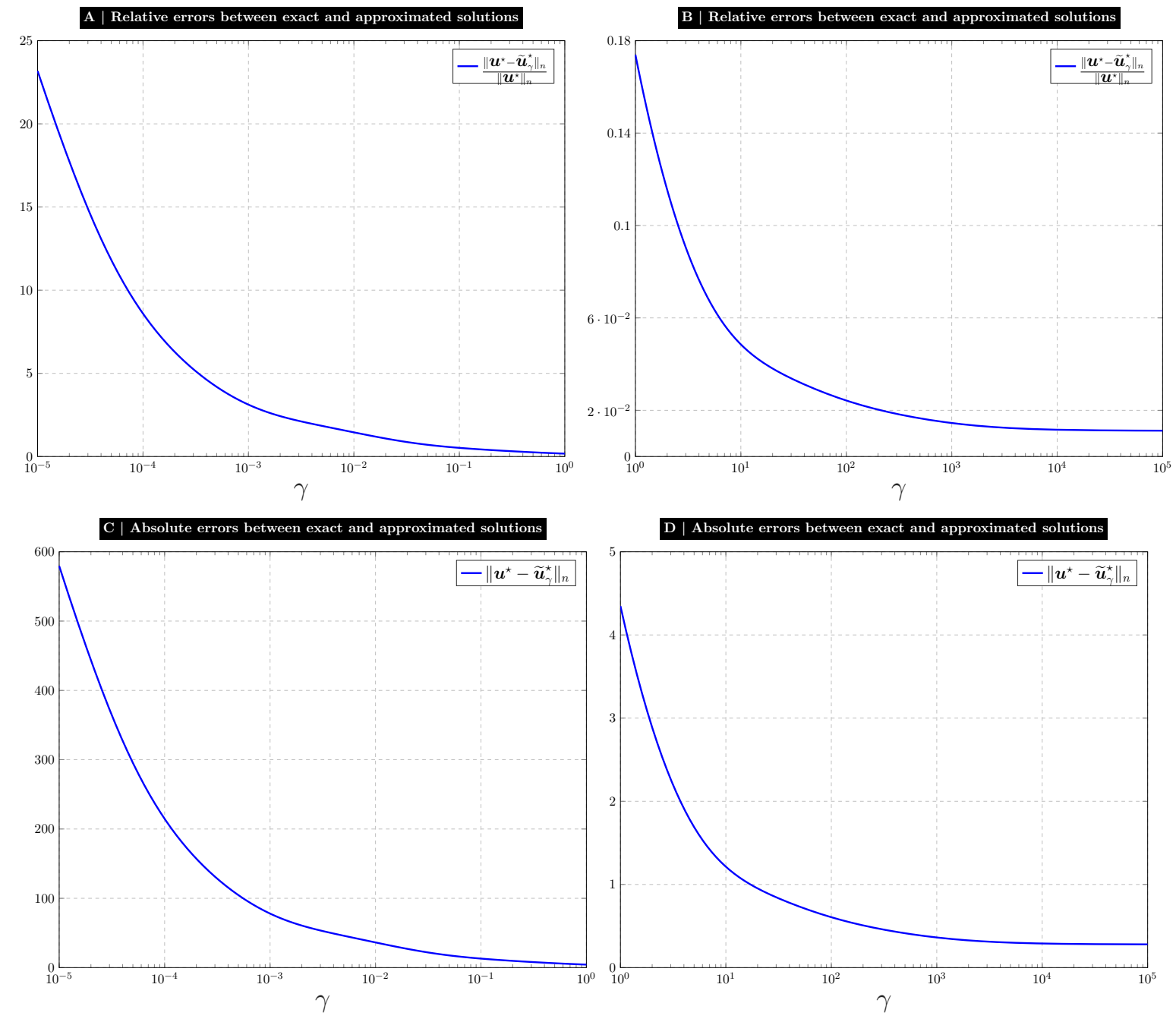

Figure 38: Absolute and relative errors between the exact solution $\boldsymbol{u}^\star$ computed from (68), and the stabilized-regularized solution $\widetilde{\boldsymbol{u}}^\star_\gamma$ obtained from the solving of (22) by the Matlab command \, with respect to the parameter $\gamma \in [10^{-5}, 1]$ (in A and C) and $\gamma \in [1, 10^5]$ (in B and D). The data $\mathbf{A}$, $\mathbf{L}$ and $\boldsymbol{b}$ are computed from the Matlab codes *gravity(n,1,a,b,h)* and *get_l*$(n, 2)$ provided in *Regularization Tools* [22] under the dimensions $n = 1\,000$ (where $\mathbf{cond}(\mathbf{A}) = 2.5999e^{20+}$, $\mathbf{rank}(\mathbf{A}) = 45$), and the right hand term $\gamma \mathbf{L}^t \boldsymbol{g}$ in (22) is not taken into account since $\boldsymbol{g} = \mathbf{0}$. The computations are made under the perturbed data in (56) and the horizontal axis is displayed on a log-scale.

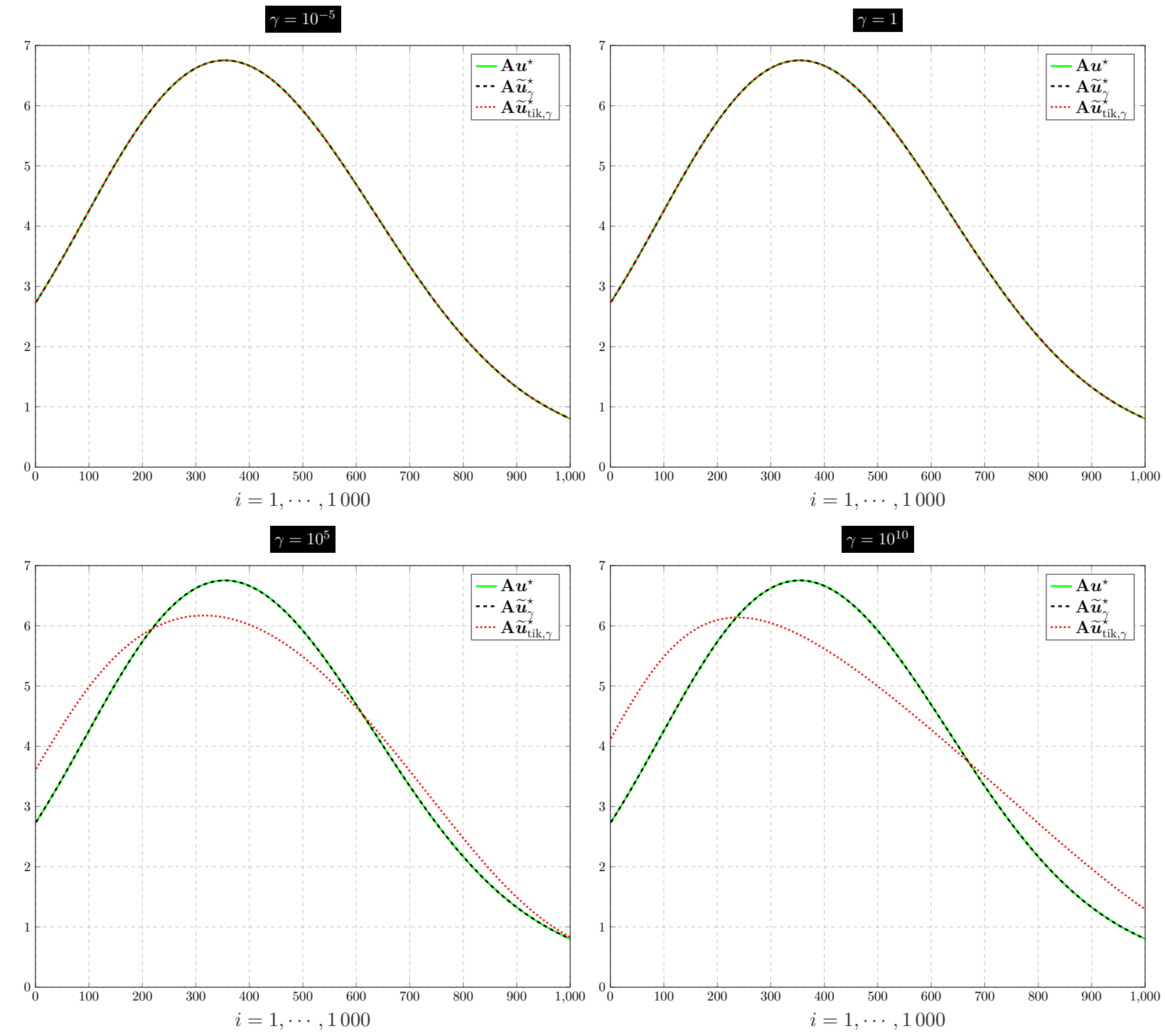

Figure 39: The exact data $\mathbf{A}\boldsymbol{u}^\star$ (in green solid line), the data reconstructed from the stabilized-regularized approach $\mathbf{A}\widetilde{\boldsymbol{u}}_\gamma^\star$ (in black dashed line), and the data made from the Tikhonov regularization $\mathbf{A}\widetilde{\boldsymbol{u}}_{\text{tik},\gamma}^\star$ (in red dotted line), are presented with respect to the discretisation index $i = 1, \cdots, 1\,000$, for the different values of the regularization parameter $\gamma = 10^{-5}, 1, 10^{5}, 10^{10}$ (given in the top of the plots). The stabilized-regularized solution $\widetilde{\boldsymbol{u}}_\gamma^\star$ is obtained from the solving of (22) through the Matlab command \, and the Tikhonov regularization solution $\widetilde{\boldsymbol{u}}_{\text{tik},\gamma}^\star$ is computed in the solving of (5) by the Matlab code *tikhonov( )* implemented in *Regularization Tools* [22]. The data $\mathbf{A}$, $\mathbf{L}$ and $\boldsymbol{b}$ are computed from the Matlab codes *gravity(n,1,a,b,h)* and *get_l*$(n, 2)$ in *Regularization Tools* [22] and the approximate derivative is fixed to $\boldsymbol{g} = \mathbf{0}$. The computations are made under the perturbed data in (56), with respect to the dimensions $n = 1\,000$ (where $\mathbf{cond}(\mathbf{A}) = 2.5999e^{20+}$, $\mathbf{rank}(\mathbf{A}) = 45$) and the parameters in (49).

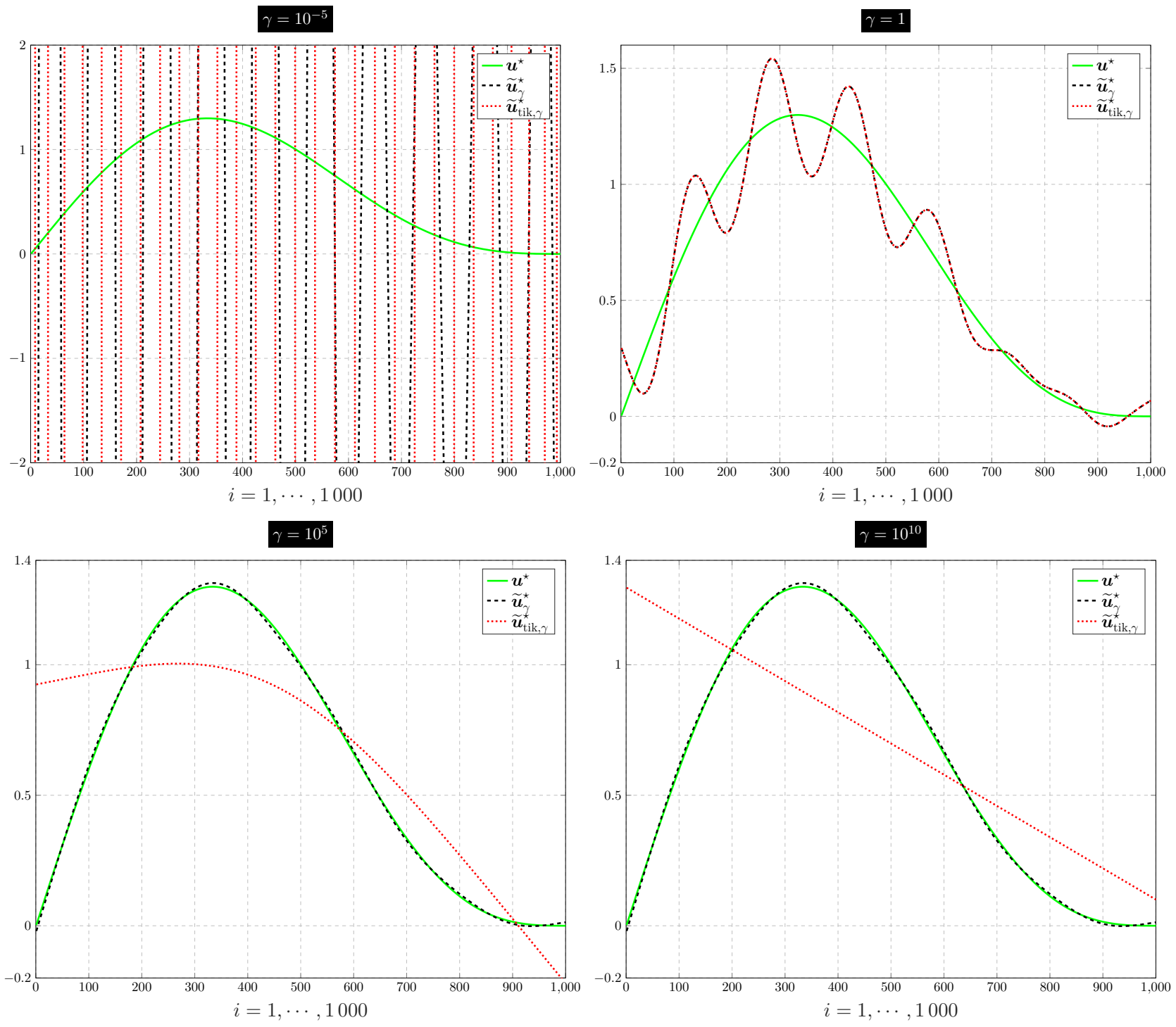

Figure 40: The exact solution $\boldsymbol{u}^\star$ (in green solid line) computed from (68), the stabilized-regularized solution $\widetilde{\boldsymbol{u}}_\gamma^\star$ (in black dashed line) obtained from the solving of (22) by the Matlab command \, and the Tikhonov regularization solution $\widetilde{\boldsymbol{u}}_{\text{tik},\gamma}^\star$ (in red dotted line) obtained in the solving of (5) by the Matlab code *tikhonov( )* implemented in *Regularization Tools* [22], are presented with respect to the discretization index $i = 1, \cdots, 1\,000$, for different values of the parameter $\gamma = 10^{-5}, 1, 10^{5}, 10^{10}$. The data $\mathbf{A}$, $\mathbf{L}$ and $\boldsymbol{b}$ are computed from the Matlab codes *gravity(n,1,a,b,h)* and *get_l*$(n, 2)$ in *Regularization Tools* [22] and the term $\gamma \mathbf{L}^t \boldsymbol{g}$ in (5) and (22) is not taken into account ($\boldsymbol{g} = \mathbf{0}$). The computations are made under perturbed data in (56), with dimensions $n = 1\,000$ (where $\mathbf{cond}(\mathbf{A}) = 2.5999e^{20+}$, $\mathbf{rank}(\mathbf{A}) = 45$) and parameters in (49).

- $h = 1, \quad \eta = 10^{-3}$

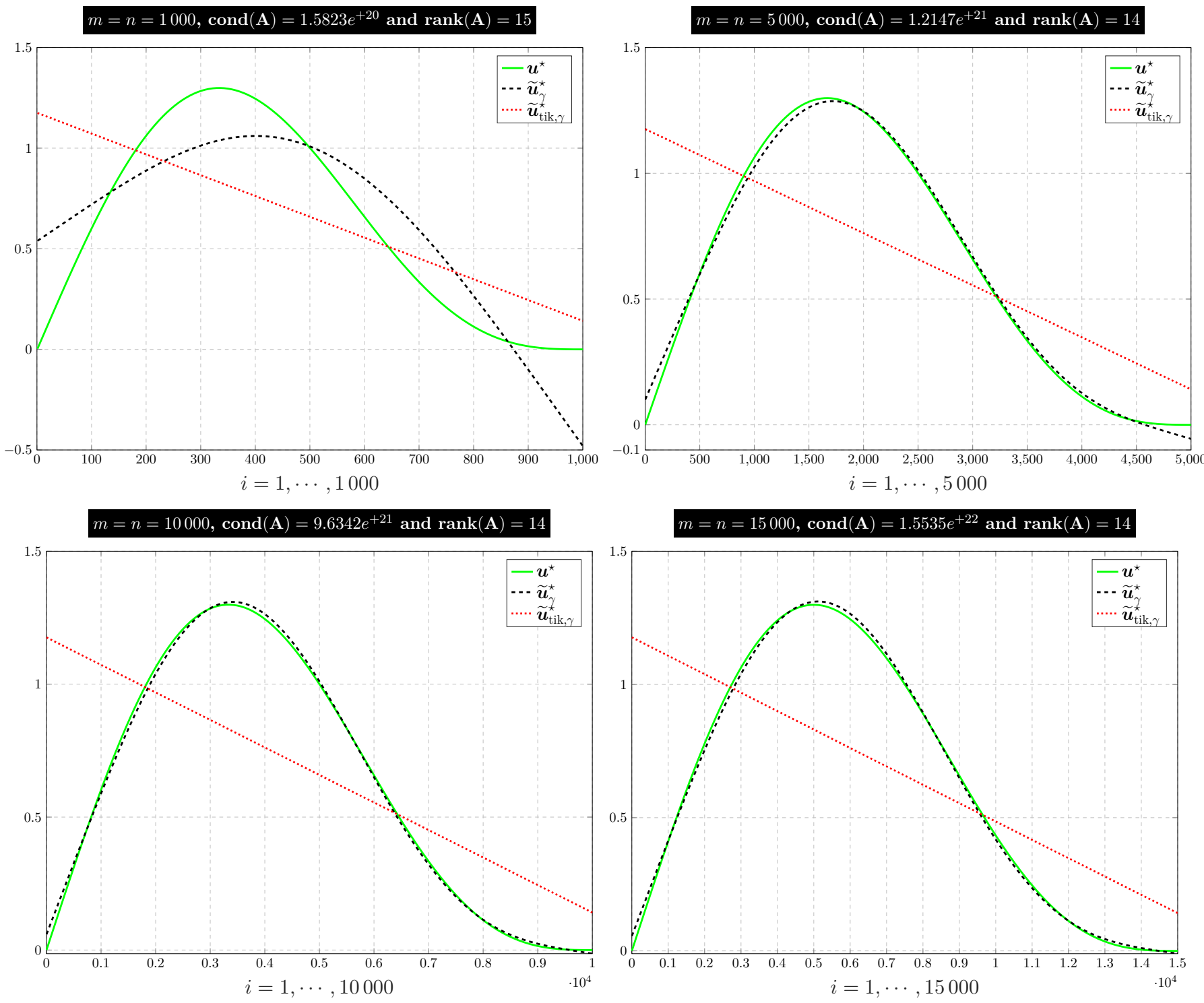

Figure 41: The exact solution $\boldsymbol{u}^\star$ (in green solid line) computed from (68), the stabilized-regularized solution $\widetilde{\boldsymbol{u}}^\star_\gamma$ (in black dashed line) obtained from the solving of (22) by the Matlab command \, and the Tikhonov regularization solution $\widetilde{\boldsymbol{u}}^\star_{\mathrm{tik},\gamma}$ (in red dotted line) obtained in the solving of (5) by the Matlab code *tikhonov( )* implemented in *Regularization Tools* [22], are presented with respect to the discretization index $i = 1, \cdots, m$, for different values of the dimensions $n = 1\,000, 5\,000, 10\,000, 15\,000$ and of the fixed parameter $\gamma = 10^{10}$. The data $\mathbf{A}$, $\mathbf{L}$ and $\boldsymbol{b}$ are computed from the Matlab codes *gravity(n,1,a,b,h)* and *get_l*$(n, 2)$ in *Regularization Tools* [22] and the term $\gamma \mathbf{L}^t \boldsymbol{g}$ in (5) and (22) is not taken into account ($\boldsymbol{g} = \mathbf{0}$). The computations are made under perturbed data in (56).

- $h = 1, \quad \eta = 10^{-2}$

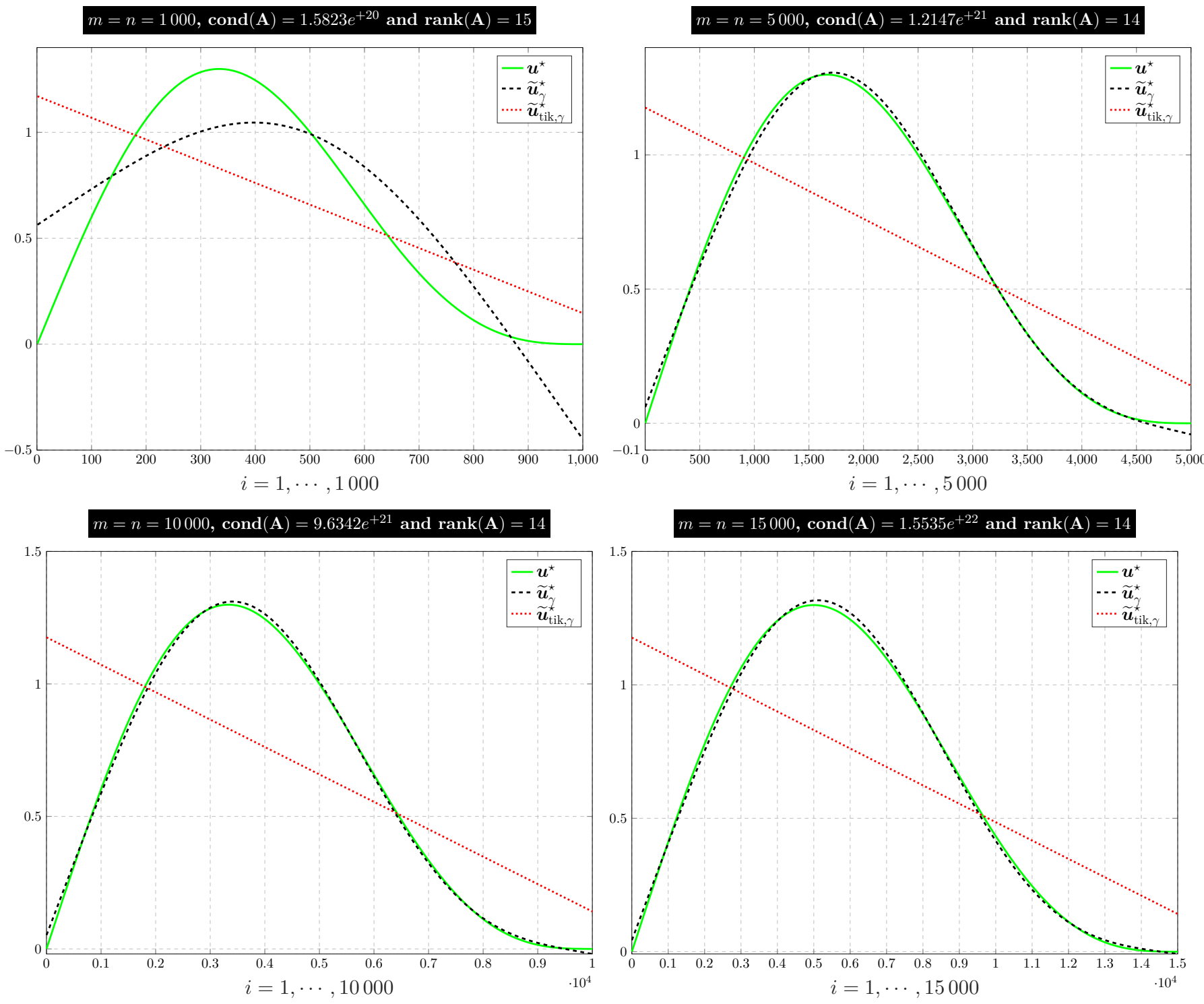

Figure 42: The exact solution $\boldsymbol{u}^\star$ (in green solid line) computed from (68), the stabilized-regularized solution $\widetilde{\boldsymbol{u}}^\star_\gamma$ (in black dashed line) obtained from the solving of (22) by the Matlab command \, and the Tikhonov regularization solution $\widetilde{\boldsymbol{u}}^\star_{\text{tik},\gamma}$ (in red dotted line) obtained in the solving of (5) by the Matlab code *tikhonov( )* implemented in *Regularization Tools* [22], are presented with respect to the discretization index $i = 1, \cdots, m$, for different values of the dimensions $n = 1\,000, 5\,000, 10\,000, 15\,000$ and of the fixed parameter $\gamma = 10^{10}$. The data $\mathbf{A}$, $\mathbf{L}$ and $\boldsymbol{b}$ are computed from the Matlab codes *gravity(n,1,a,b,h)* and *get_l*$(n, 2)$ in *Regularization Tools* [22] and the term $\gamma \mathbf{L}^t \boldsymbol{g}$ in (5) and (22) is not taken into account ($\boldsymbol{g} = \mathbf{0}$). The computations are made under perturbed data in (56).